%% file: paratuck2_algo.tex
\documentclass[a4paper,dvipsnames]{siamart250211}
\input{preamble.tex}
\input{macros.tex}

\newcommand{\TheTitle}{A lifting approach to \\ ParaTuck-2 tensor decompositions}
\newcommand{\TheShortTitle}{A lifting approach to ParaTuck-2}

\newcommand{\TheAuthors}{K. Usevich}
\headers{\TheShortTitle}{\TheAuthors}

\title{\TheTitle\thanks{Submitted to the editors DATE.
\funding{This work was supported by the ANR projects LeaFleT (ANR-19-CE23-0021-01) and AGDAM (ANR-23-CE94-0001).}}}

\author{
  Konstantin Usevich\thanks{Universit\'{e} de Lorraine and CNRS, CRAN (Centre de
    Recherche en Automatique en Nancy), UMR 7039, Campus Sciences, BP 70239,
    54506 Vand\oe{}uvre-l\`{e}s-Nancy cedex, France
    (\email{konstantin.usevich@cnrs.fr}).}%
}

\begin{document}
\maketitle

\begin{abstract}
The ParaTuck-2 decomposition (PT2D) of third-order tensor is a two-layer generalization of the well-known canonical polyadic decomposition (CPD).
While being more flexible than the CPD, the PT2D also possesses similar uniqueness properties.
In this paper, we show than under the best known uniqueness conditions, the exact PT2D can be computed by an algebraic algorithm (i.e., can the PT2D problems can be reduced to computing nullspaces and eigenvalues of certain matrices).
We do so by lifting the slices of the tensor to higher-dimensional space, which also allows for refining the existing uniqueness conditions.
The algorithms are developed for general PT2D and its symmetric version (DEDICOM), which leads to an algebraic algorithm for another generalization of the CPD, the PARAFAC2 decomposition.
Our methods are also applicable in the approximation scenario, as shown by the numerical experiments.
\end{abstract}



\section{Introduction}
ParaTuck-2 decomposition \cite{HarsL96:uniqueness} (or PT2D for short) decomposes a collection of $\edim{3}$ matrices $\matr{T}_{k} \in  \FF^{\edim{1} \times \edim{2}}$, $k=1,\ldots,\edim{3}$, (where  $\FF = \RR$ or $\CC$),
jointly  as a product
\begin{equation}\label{eq:paratuck_slices}
\matr{T}_{k} = \matr{A} \DiagA{k} \matr{\CoreFront}  \DiagB{k}  \matr{B}^{\T},
\end{equation}
where the three matrix factors
$\matr{A} \in \FF^{\edim{1} \times \rkRow}, \matr{B} \in \FF^{\edim{2} \times \rkCol}, \matr{\CoreFront} \in \FF^{\rkRow \times \rkCol}
$
are common for all $\matr{T}_{k} $ and the middle factors $\DiagA{k}$ and $\DiagB{k}$, $k=1,\ldots, \edim{3}$ are $R \times R$ and $S \times S$  diagonal matrices, respectively, that
depend on the index $k$.
The PT2D is a decomposition of an $\edim{1} \times \edim{2} \times \edim{3}$ 3rd-order tensor (whose frontal slices $\edim{1} \times \edim{2}$ are the matrices $\matr{T}_k$) and the decomposition problem is to recover the unknown factors given the $\matr{T}_k$.

The PT2D \eqref{eq:paratuck_slices} can be viewed as a ``2-layer'' generalization of the usual joint matrix factorization ($\matr{T}_k = \matr{A} \matr{D}^{(k)} \matr{B}^{\T}$, $\matr{D}^{(k)}$ diagonal, $k=1,\ldots,\edim{3}$), thus the PT2D is a generalization of the well-known CP (canonical polyadic) decomposition \cite{KoldB09:tensor}.
The name ParaTuck-2  was suggested in \cite{HarsL96:uniqueness} to show that PT2D shares the features of  PARAFAC\footnote{one of the older names for the CP decomposition} and Tucker tensor decompositions \cite{KoldB09:tensor}.


PT2D appears in several applications.
In multi-way analysis, the PT2D was proposed in  the context psychometrics \cite{HarsL96:uniqueness} and chemometrics \cite{Bro98:thesis} to account for additional variability in the joint  factor analysis of  several datasets.
In the context of telecommunications, PT2D was introduced as an encoding scheme  (as well as various generalizations \cite{FaviA14:constrained}) in the context of multiple-input multiple-output wireless communications \cite{AlmeFX13:MIMO,OlivFFB19:mimo,Nask20:thesis}.
Recently, PT2D appeared in block-structured system identification and neural network approximation (in the frame of so-called decoupling approach \cite{DreeIS15:decoupling}), where the CPD corresponds to models with 1 hidden layer  \cite{DreeIS15:decoupling} while the 2-layer case gives rise to the PT2D \cite{DeJoUDI23:CAMSAP,UsevZIDA23:EUSIPCO}.
The fully symmetric version of PT2D (case $\matr{A} =\matr{B}$, $\DiagA{k}=\DiagB{k}$ in \eqref{eq:paratuck_slices}), also known as the three-way DEDICOM\footnote{The name is an abbreviation of ``decomposition into directional components''} \cite{Hars78:dedicom}, is a way to describe asymmetric relationships  in multivariate data analysis (see \cite{KoldB09:tensor}). DEDICOM is also tightly linked to the popular PARAFAC-2 model \cite{RoalSCABCA22:AOADMM}; in particular,  DEDICOM is used as a backbone in so-called indirect approaches to PARAFAC-2 decomposition \cite{Kier93:ALS}.

One of the key features of the PT2D that makes it appealing in applications is its (essential) uniqueness properties (up to trivial permutations and scalings) (similarly to the CPD). As shown in \cite{HarsL96:uniqueness}, for the case $ \rkRow \le\edim{1}$, $\rkCol \le \edim{2}$ (full column rank $\matr{A},\matr{B}$), and $\edim{3}$ is sufficiently large, the PT2D is essentially unique under mild (generic) conditions on the factors $\matr{A}$, $\matr{B}$, $\matr{\CoreFront}$, $\DiagA{k}$,  $\DiagB{k}$.
A similar uniqueness result was also proved in \cite{HarsL96:uniqueness} for a partially symmetric version of PT2D (a generalized DEDICOM), which also serves as a base for the best known uniqueness results of the PARAFAC-2 decomposition \cite{BergK96:uniqueness}.
However, up to the author's knowledge, no significant progress has been made on uniqueness of PT2D and DEDICOM/PARAFAC2 since \cite{HarsL96:uniqueness}.

Despite the usefulness and the nice features of PT2D, reliable algorithms are lacking, especially in the non-symmetric case.
As for local optimization algorithms, many standard optimization tools, such as ALS (alternating least squares \cite[pp. 68--71]{Bro98:thesis})  suffer from slow convergence and local minima.
In many applications, a simplifying assumptions is introduced  that at least one of the factors (e.g., example $\matr{A}$) is known \cite{OlivFFB19:mimo}, \cite[p. 213]{Bro98:thesis}, which simplifies the problem.
Several recent efforts focused on improving local optimization strategies, such as double contractions in ALS in \cite{Nask20:thesis} or stochastic optimization \cite{ZniyA25:stochastic},but these methods still lack convergence proofs.
 PT2D can be also reformulated as a structured CPD \cite{FaviA14:constrained}, but this reformulation do not currently yield a reliable way to compute the decomposition (as the resulting tensor rank is too high). 
For DEDICOM (as well as for PARAFAC-2 \cite{RoalSCABCA22:AOADMM}), the available algorithms mainly employ (local)  alternating optimization minimization, such as ALS \cite{Kier89:ALS,Kier93:ALS} or ASALSAN\footnote{alternating simultaneous approximation, least squares and Newton} \cite{BadeHK07:temporal}, and are based on efficient techniques for computing updates.
Finally, the uniqueness results of  \cite{HarsL96:uniqueness} do not provide an explicit (algebraic) algorithm for  PT2D or DEDICOM (unlike  the algebraic algorithms to find the CPD under  Kruskal's uniqueness conditions \cite{KoldB09:tensor,DomaL14:gevd}).



In this paper, we show that, by lifting $\matr{T}_k$ to a higher-dimensional space, we can design  algebraic algorithms (relying on standard linear algebra tools) for the PT2D that work under the uniqueness conditions in  \cite{HarsL96:uniqueness} (i.e. $\rkRow \le \edim{1}, \rkCol \le \edim{2}$ and $\edim{3}$ sufficiently large).
We provide algebraic algorithms for the non-symmetric PT2D and DEDICOM (and as a consequence for PARAFAC2), which are first of this kind, up to the author's knowledge.
The lifting approach relies on properties of the kernel of a certain structured matrix constructed from $\matr{T}_k$.
This approach also  allows us to relax and clarify uniqueness results in \cite{HarsL96:uniqueness}.
The proofs are split into many intermediate steps and highlight the importance of each of the assumptions, as well as the structure of the core tensor. 

The paper is structured as follows.
Section 2 contains basic matrix and tensor notation used in the paper.
Section 3 recalls  known facts about PT2D, including ambiguities and a summary of the results of \cite{HarsL96:uniqueness}, and gives a preview of the results of the paper.
Section 4 provides an alternative definition of the PT2D through the core tensor and triple product of matrices and discusses the properties of the core tensor that are crucial for the following.
Sections 5 and 6 contain the  main results  and algorithms for the non-symmetric and the fully symmetric cases, respectively
In section 7, we discuss the advanced properties of the core tensor and its implications.
Section 8 contains a preview of the proof for a special case $R=S=2$.
Sections 9 and 10 contain the proofs for the nonsymmetric and symmetric case, respectively.
Finally, Section 11 contains details on implementation and the numerical experiments.

In the paper, for simplicity, we assume $\FF= \RR$ (however, all the results can be easily translated to the  case $\FF=\CC$).

\section{Matrix and tensor notation}
Vectors, matrices, and tensors,  will be respectively denoted by bold lowercase letters (\textit{e.g.} $\vect{u}$), with bold uppercase letters (\textit{e.g.} $\matr{M}$), and with bold calligraphic letters, (\textit{e.g.} $\tens{A}$).
For simplicity, here vectors, matrices, and tensors, are viewed as 1-way, 2-way and d-way arrays respectively. Their corresponding entries will be denoted by $\tenselem{A}_{ijk}$, $M_{ij}$, and $u_i$.
The vectors $\vect{v} \in \FF^{n}$ in this paper are always thought of as column vectors (i.e, matrices $\FF^{n \times 1}$).
The notation $\FF^{n_1 \times \ldots \times n_d} = (\FF^{n_1}) \otimes \cdots \otimes  (\FF^{n_d})$ denotes the space of $n_1 \times \cdots \times n_d$ tensors (where $\otimes$ stands for the tensor product).

For a tensor $\tens{T} \in \FF^{n_1 \times \ldots \times n_d}$  we denote by $\vecl{\tens{T}} \in \FF^{n_1  \cdots n_d}$ the standard column-major vectorization of the tensor.
For example, for the matrix  $\matr{V}  = \begin{bmatrix} \vect{v}_1 & \cdots & \vect{v}_{\edim{2}} \end{bmatrix}\in \FF^{\edim{1} \times \edim{2}}$, the column-major vectorization stacks its columns
\[
\vecl{\matr{V}}    =\begin{bmatrix} \vect{v}_1^{\T} & \cdots & \vect{v}_{\edim{2}}^{\T} \end{bmatrix}^{\T}
\]
We also allow the vectorization to be applied to tensor subspaces, e.g. $\vecl{\set{V}} \subseteq \FF^{IJK}$ for $\set{V} \subseteq \FF^{I\times J \times K}$
By $\matricize{\edim{1}}{\edim{2}}{\vect{v}}$ we denote $(\edim{1},\edim{2})$-matricization of a vector $\vect{v} \in \FF^{\edim{1}\edim{2}}$ the matrix $\matr{V}$ such that $\vect{v} = \vecl{\matr{V}}$.
We will equally allow the matricization to be applied directly on tensors, that is $\matricize{m}{n}{\tens{T}} = \matricize{m}{n}{\vecl{\tens{T}}}$.
We also denote by $\tensorize{\edim{1}}{\edim{2}}{\edim{3}}{\vect{v}} \in \RR^{{\edim{1}}\times{\edim{2}}\times{\edim{3}}}$ the tensorization of the vector $\RR^{{\edim{1}}{\edim{2}}{\edim{3}}}$

We use $\Span{\cdot}$ to denote  linear space spanned by a set of vectors (matrices, or tensors) and $\range{\matr{Z}}$ the range of a matrix.
We use $\lker{\matr{Z}}$ to denote the left kernel (i.e., orthogonal complement of $\range{\matr{Z}}$).
Operator $\Diag{\cdot}$ denotes the diagonal matrix constructed from a vector, or a block-diagonal matrix built from a sequence of matrices.

Operator $\contr{p}$ denotes contraction on the $p$th index of a tensor; when contracted with a matrix, it is understood that summation is always performed on the second index of the matrix. For instance, $(\tens{A}\contr{1}\matr{M})_{ijk}=\sum_\ell \tenselem{A}_{\ell jk} M_{i\ell}$. We denote by $\matr{T}^{(1)} \in \FF^{\edim{1}\times(\edim{2} \edim{3})}$ and $\matr{T}^{(2)}\in \FF^{\edim{2}\times(\edim{1} \edim{3})}$, $\matr{T}^{(3)}\in \FF^{\edim{3}\times(\edim{1} \edim{2})}$ the first, second and third unfoldings of a tensor $\tens{T} \in \FF^{\edim{1} \times \edim{2} \times \edim{3}}$, respectively. 
Note that $\matr{T}^{(1)} = \matricize{\edim{1}}{\edim{2}\edim{3}}{\tens{T}}$  and $\matr{T}^{(3)} = (\matricize{\edim{1}\edim{2}}{\edim{3}}{\tens{T}})^{\T}$.

The symbol $\tensp$ is used for tensor (outer) product of vectors (or tensors). That is, $\tens{T} = \vect{a} \tensp \vect{b} \tensp \vect{c}$ if $\tenselem{T}_{ijk} = a_ib_jc_k$.
To avoid confusion, we use the symbol $\kron$ for the Kronecker product of matrices. The usual identity for vectorization of products of matrices holds, $\vecl{\matr{A}\matr{X}\matr{B}^{\T}} = (\matr{B}\kron \matr{A}) \vecl{\matr{X}}$.
This implies, that the vectorization reverses the order of products for rank-one vectors, we have $\vecl{\vect{a} \tensp \vect{b}  \tensp \vect{c}} = \vect{c} \kron \vect{b} \kron \vect{a}$.
We use $\kr$ for the Khatri-Rao product of matrices, i.e., for two matrices 
\begin{equation}\label{eq:matr_AB}
\matr{A} = \bmx \vect{a}_1 &  \cdots & \vect{a}_R \emx, \matr{B} = \bmx \vect{b}_1 &  \cdots & \vect{b}_R \emx,
\end{equation}
we have $\matr{B} \kr \matr{A} = \bmx \vect{b}_1 \kron \vect{a}_1  &  \cdots &\vect{b}_R \kron \vect{a}_R \emx$. 

We will use the common notation 
\[
\cpd{\matr{A}}{\matr{B}}{\matr{C}} = \sum\limits_{r=1}^R \vect{a}_k \tensp \vect{b}_k \tensp \vect{c}_k
\]
 to denote an $R$-term (canonical\footnote{We will not make a distinction between polyadic and canonical polyadic decomposition (in some articles, the term ``canonical'' is used only for $\cpd{\matr{A}}{\matr{B}}{\matr{C}}$ where $R$ is the minimal possible such $R$).}) polyadic (CP) decomposition, or CPD, of an $\edim{1} \times \edim{2} \times \edim{3}$ tensor with $\matr{A},\matr{B}$ as in \eqref{eq:matr_AB} and $\matr{C} = \bmx \vect{c}_1 &  \cdots & \vect{c}_R \emx$.

We view symmetric tensors as subspace of order-$d$ $n\times \cdots \times n$ tensors, denoted by $S^{d}(\RR^{n}) \subset \RR^{n\times \cdots\times n}$.
For example $S^{2}(\RR^{n})$ denotes the vector space of symmetric $n\times n$ matrices.
We denote by $\sigma_d: \RR^{n^d} \to \vecl{S^{d}(\RR^{n})}$ the operation that symmetrizes of the tensor given in vectorized form (if is also an orthogonal projection on $\vecl{S^{d}(\RR^{n})}$).
For example, the symmetrization map $\symglobal_2: \RR^{n^2} \to \RR^{n^2}$ maps
$\vecl{\matr{X}} \mapsto \vecl{\frac{\matr{X} + \matr{X}^{\T}}{2}}$.
In particular, $\sigma_2(\vect{y} \kron \vect{z}) = \frac{1}{2}(\vect{y} \kron \vect{z} + \vect{z} \kron \vect{y})$ for $\vect{y}, \vect{z} \in \RR^{n}$.
Finally, for a linear  subspace $\set{V} \subset \RR^{n}$ we denote by $S^2(\set{V}) \subset \RR^{n\times n}$ the set of symmetric matrices with range in $\set{V}$.

\section{Background on  PT2D and its variants}
\subsection{ParaTuck-2 decomposition: notation and basic considerations}
We first introduce the following compact notation for the factors of the PT2D (that groups the diagonal factors $\DiagA{k}$ and $\DiagB{k}$ into two dense matrices $\matr{F}$ and $\matr{G}$).
\begin{definition}\label{def:PT2D}
A third-order tensor $\tens{T} \in \RR^{\edim{1} \times \edim{2} \times \edim{3}}$ is said to admit a rank-$(R,S)$ ParatTuck-2 decomposition with \emph{factors}
\[
\matr{A} \in \FF^{\edim{1} \times \rkRow}, \matr{B} \in \RR^{\edim{2} \times \rkCol}, \matr{\CoreFront} \in \RR^{\rkRow \times \rkCol}, 
\matr{\CoreVert} \in \RR^{\rkRow \times \edim{3}}, \matr{\CoreHoriz} \in \RR^{\rkCol \times \edim{3}},
\]
if its frontal slices  can be expressed as $\tens{T}_{:,:,k} = \matr{T}_k$ with $\matr{T}_k$ given in \eqref{eq:paratuck_slices}, where the  $\DiagA{k}$ and $\DiagB{k}$ are defined from the columns of
\begin{equation}\label{eq:columns_core_factors}
 \matr{\CoreVert} = \bmx \vect{g}_1 &  \cdots & \vect{g}_{\edim{3}} \emx, \quad\matr{\CoreHoriz} =  \bmx \vect{h}_1 &  \cdots & \vect{h}_{\edim{3}} \emx
 \end{equation}
 as follows:
\begin{equation}\label{eq:DiagAB}
\DiagA{k}  = \Diag{\vect{g}_k}, \quad \DiagB{k}  = \Diag{\vect{h}_k}.
\end{equation}
We will write use the notation 
\[
\tens{T} = \ptd{\matr{A}}{\matr{B}}{\matr{\CoreFront}}{\matr{\CoreVert}}{\matr{\CoreHoriz}},
\]
to denote that $\tens{T}$ has a PT2D with factors $\matr{A},\matr{B},\matr{F},\matr{G},\matr{H}$.
\end{definition}

\begin{example}
As mentioned in the introduction, CPD is a special case of PT2D for a particular choice of factors. Indeed, for $\rkRow=\rkCol$, $\boldsymbol{1}_{\rkRow\times \edim{3}}$ denoting an $R \times \edim{3}$ matrix of all ones, we have
\[
\cpd{\matr{A}}{\matr{B}}{\matr{C}}  = \ptd{\matr{A}}{\matr{B}}{\matr{I}}{\matr{C}}{\boldsymbol{1}_{R\times \edim{3}}}.
\]
Note that there is a plenty of ways to express a CPD as a special case of the PT2D. 
\end{example}
\begin{remark}
Unlike the tensor rank (minimal $R$ such that the tensor admits an $R$-term CPD), for PT2D, we cannot  define  the ``minimal'' rank in a unique way, because there is a pair $(R,S)$ of ``ranks''.
In what follows, we will use the term ``rank-(R,S)'' ParaTuck-2 decomposition (or $(R,S)$-PT2D) without assuming that this is a minimal rank.
\end{remark}

\subsection{Symmetric versions of PT2D and DEDICOM}
Several special cases of PT2D exist (see the introduction in \cite{HarsL96:uniqueness}) were proposed:
\begin{itemize}
\item symmetrically weighted PT2D \cite{HarsL96:uniqueness} : $\rkRow=\rkCol$, $\matr{G}=\matr{H}$ (which corresponds to $\DiagA{k} = \DiagB{k}$ in \eqref{eq:paratuck_slices});

\item fully symmetric case (also known as DEDICOM  \cite{HarsL96:uniqueness,KoldB09:tensor} ): $R=S$, $\matr{G} = \matr{H}$, $\matr{A} = \matr{B}$ and $\matr{F}$ symmetric.
We will  use the notation
\[
\ptdsym{\matr{A}}{\matr{F}}{\matr{G}} :=  \ptd{\matr{A}}{\matr{A}}{\matr{F}}{\matr{G}}{\matr{G}}.
\]
\end{itemize}
\begin{remark}
For any matrix $\matr{F}$, we have $\ptdsym{\matr{A}}{\matr{F}}{\matr{G}} = \ptdsym{\matr{A}}{{(\matr{F}+\matr{F}^{\T})}/{2}}{\matr{G}}$.
Therefore it makes no sense to consider nonsymmetric  $\matr{F}$ in DEDICOM.
\end{remark}

The fully symmetric is it tightly linked to the PARAFAC-2 decomposition, as explained below.
\begin{example}\label{ex:PT2_PF2}
Let $\matr{X}_k \in \RR^{\edim{1} \times J_k}$, $k =1,\ldots,\edim{3}$  be a collection of matrices that are decomposed as
\begin{equation}\label{eq:PARAFAC2}
\matr{X}_k = \matr{A} \matr{D}_k \matr{B}_k^{\T}, 
\end{equation}
where $\matr{A} \in \RR^{\edim{1} \times \rkRow}$, $\matr{D}_k \in \RR^{\rkRow \times \rkRow}$ and $\matr{B}_k \in \RR^{\rkRow \times J_k}$ satisfy $\matr{B}_k \matr{B}_k^{\T} = \matr{\Sigma}$.
Such a decomposition is known as  PARAFAC-2 decomposition (with common covariance constraint).

PARAFAC-2 decomposition is linked to DEDICOM as follows. Under the common covaraince constraint, the  matrices $\matr{M}_k = \matr{X}_k \matr{X}_k^{\T}$  can be expressed as 
\[
\matr{M}_k = \matr{A} \matr{D}_k \matr{\Sigma} \matr{D}_k \matr{A}^{\T}
\]
Thus the ``covariance'' tensor $\tens{M} \in \RR^{\edim{1} \times \edim{1} \times \edim{3}}$ with $\tens{M}_{:,:,k} = \matr{M}_k$ admits the rank-$R$ fully symmetric PT2D ($\tens{M} =  \ptdsym{\matr{A}}{\matr{\Sigma}}{\matr{G}}$).
\end{example}
The relation between DEDICOM and PARAFAC-2 decomposition in  \Cref{ex:PT2_PF2} is in fact used to establish the uniqueness properties of PARAFAC-2 \cite{BergK96:uniqueness}, based on uniqueness proofs for DEDICOM in \cite{HarsL96:uniqueness}.
Also, DEDICOM is  used to compute the PARAFAC-2 decomposition using the indirect approach \cite{Kier93:ALS}.

\subsection{Ambiguities and uniqueness}
In this section, we will recall the trivial ambiguities of PT2D (which are similar to the ones of the CPD) and define the notion of (essential) uniqueness.

As in the case of the CPD, there are permutation and scaling ambiguities. For example, permuting the columns of $\matr{A}$  together with the same permutation of rows $\matr{F}$ and $\matr{G}$, leads to an alternative PT2D. 
We will summarize these ambiguities in the following lemma.
\begin{lemma}\label{lem:ambiguities}
The following transformations (trivial ambiguities) give an alternative PT2D:
\begin{enumerate} 
\item 
$\ptd{\matr{A}}{\matr{B}}{\matr{\CoreFront}}{\matr{\CoreVert}}{\matr{\CoreHoriz}} = \ptd{\matr{A}\matr{\Pi}_{\matr{A}}}{\matr{B}}{\matr{\Pi}_{\matr{A}}^{\T}\matr{\CoreFront}}{\matr{\Pi}_{\matr{A}}^{\T} \matr{\CoreVert}}{\matr{\CoreHoriz}}$, for any permutation matrix $\matr{\Pi}_{A} \in \RR^{R\times R}$
\item 
$\ptd{\matr{A}}{\matr{B}}{\matr{\CoreFront}}{\matr{\CoreVert}}{\matr{\CoreHoriz}} = \ptd{\matr{A}}{\matr{B}\matr{\Pi}_{\matr{B}}}{\matr{\CoreFront}\matr{\Pi}_{\matr{B}}}{ \matr{\CoreVert}}{\matr{\Pi}_{\matr{B}}^{\T}\matr{\CoreHoriz}}$, for any permutation matrix $\matr{\Pi}_{B} \in \RR^{S\times S}$.
\item 
$\ptd{\matr{A}}{\matr{B}}{\matr{\CoreFront}}{\matr{\CoreVert}}{\matr{\CoreHoriz}} = \ptd{\matr{A}\matr{\Lambda}_{\matr{A}}}{\matr{B}}{(\matr{\Lambda}_{\matr{A}}\matr{\Lambda}_{\matr{G}})^{-1}\matr{\CoreFront}}{\matr{\Lambda}_{\matr{G}} \matr{\CoreVert}}{\matr{\CoreHoriz}}$, for all invertible diagonal  $\matr{\Lambda}_{\matr{A}}, \matr{\Lambda}_{\matr{G}}$.
\item 
$\ptd{\matr{A}}{\matr{B}}{\matr{\CoreFront}}{\matr{\CoreVert}}{\matr{\CoreHoriz}} = \ptd{\matr{A}}{\matr{B}\matr{\Lambda}_{\matr{B}}}{\matr{\CoreFront}(\matr{\Lambda}_{\matr{B}}\matr{\Lambda}_{\matr{H}})^{-1}}{\matr{\CoreVert}}{\matr{\Lambda}_{\matr{H}} \matr{\CoreHoriz}}$, for all invertible diagonal  $\matr{\Lambda}_{\matr{B}}, \matr{\Lambda}_{\matr{H}}$.
\item 
$\ptd{\matr{A}}{\matr{B}}{\matr{\CoreFront}}{\matr{\CoreVert}}{\matr{\CoreHoriz}} = \ptd{\matr{A}}{\matr{B}}{\matr{\CoreFront}}{\matr{\CoreVert}\matr{\Lambda}_{\matr{G}\matr{H}}}{\matr{\CoreHoriz}\matr{\Lambda}^{-1}_{\matr{G}\matr{H}}}$, for all invertible diagonal  $\matr{\Lambda}_{\matr{G}\matr{H}}$.
\end{enumerate}
\end{lemma}
\begin{proof}
Follows from explicit computations.
\end{proof}
We will call the PT2D unique if there are no other other PT2Ds rather than those obtained from \Cref{lem:ambiguities}.
\begin{definition}\label{def:uniqueness}
The decomposition $\tens{T} = \ptd{\matr{A}}{\matr{B}}{\matr{\CoreFront}}{\matr{\CoreVert}}{\matr{\CoreHoriz}}$ is called  (essentially)  unique if for any alternative decomposition
\[
\tens{T} =
\ptd{\widetilde{\matr{A}}}{\widetilde{\matr{B}}}{\widetilde{\matr{\CoreFront}}}{\widetilde{\matr{\CoreVert}}}{\widetilde{\matr{\CoreHoriz}}},
\]
can be obtained through a sequence of transformations described in \Cref{lem:ambiguities}, \emph{i.e.},
there exist permutation matrices $\matr{\Pi}_{\matr{A}} \in \mathbb{R}^{\rkRow \times \rkRow}$ and $\matr{\Pi}_{\matr{B}} \in \mathbb{R}^{\rkCol \times \rkCol}$ and nonsingular diagonal matrices
\[
\matr{\Lambda}_{\matr{A}},\matr{\Lambda}_{\matr{G}}  \in   \mathbb{R}^{\rkRow \times \rkRow}, \matr{\Lambda}_{\matr{B}},\matr{\Lambda}_{\matr{H}}  \in   \mathbb{R}^{\rkCol \times \rkCol}, \matr{\Lambda}_{\matr{GH}} \in \mathbb{R}^{\edim{3} \times \edim{3}} 
\]
such that
\begin{align*}
\widetilde{\matr{A}} &= {\matr{A}} \cdot (\matr{\Pi}_{\matr{A}} \cdot \matr{\Lambda}_{\matr{A}}),\quad
\widetilde{\matr{B}} = {\matr{B}} \cdot (\matr{\Pi}_{\matr{B}} \cdot \matr{\Lambda}_{\matr{B}}),\\
\widetilde{\matr{F}}&=  \matr{\Pi}_{\matr{A}}^{\T}  \cdot (\matr{\Lambda}_{\matr{G}}^{-1}  \cdot {\matr{F}} \cdot \matr{\Lambda}_{\matr{H}}^{-1}) \cdot \matr{\Pi}_{\matr{B}},\\
\widetilde{\matr{G}} &= ( \matr{\Lambda}_{\matr{A}}^{-1}  \cdot \matr{\Pi}_{\matr{A}}^{\T})  \cdot (\matr{\Lambda}_{\matr{G}} \cdot {\matr{G}} \cdot  \matr{\Lambda}_{\matr{GH}}),  \\
\widetilde{\matr{H}}&= (\matr{\Lambda}_{\matr{B}}^{-1} \cdot \matr{\Pi}_{\matr{B}}^{\T} ) \cdot  ( \matr{\Lambda}_{\matr{H}}\cdot  {\matr{H}} \cdot \matr{\Lambda}^{-1}_{\matr{GH}}).
\end{align*}
\end{definition}
In the next subsection, we recall the existing results on  uniqueness of PT2D.

\subsection{Known results: nonsymmetric case}\label{sec:nonsymmetric_uniqueness}
The uniqueness results of \cite{HarsL96:uniqueness} rely on the following assumptions.
\begin{assumption}\label{as:AB_full_rank}
The matrices $\matr{A}$ and $\matr{B}$ are of full column rank (i.e. $\rank{\matr{A}} = R, \rank{\matr{B}} = S$).
\end{assumption}
\begin{assumption}\label{as:nonzero_F}
The matrix $\matr{\CoreFront}$ does not have zero elements.
\end{assumption}

\begin{assumption}\label{as:maxrank_F}
The matrix $\matr{\CoreFront}$ is maximal possible rank (i.e $\rank{\matr{\CoreFront}} = \min(R,S)$).
\end{assumption}
\begin{assumption}\label{as:nonzero_GH}
All elements of $\matr{G}$ and $\matr{H}$ are non-zero  (equivalently, all   $\DiagA{k}$, $\DiagB{k}$ are nonsingular).
\end{assumption}

\begin{assumption}\label{as:full_rank}
The matrix $\matr{\CoreVert} \kr \matr{\CoreVert} \kr \matr{\CoreHoriz}  \kr  \matr{\CoreHoriz}$ has the maximal possible rank
\[
\rank{\matr{\CoreVert} \kr \matr{\CoreVert} \kr  \matr{\CoreHoriz}  \kr  \matr{\CoreHoriz}} = \frac{\rkRow(\rkRow+1)}2 \frac{\rkCol(\rkCol+1)}2.
\]
\end{assumption}
Then the following theorem holds true
\begin{theorem}[{\cite[Theorem 1]{HarsL96:uniqueness}}]\label{thm:HL_nonsymmetric}
Let $\tens{T} = \ptd{\matr{A}}{\matr{B}}{\matr{\CoreFront}}{\matr{\CoreVert}}{\matr{\CoreHoriz}}$ be an $(R,S)$-PT2D whose factors satisfy \Crefrange{as:AB_full_rank}{as:full_rank}.
Then the PT2D is unique in the sense of \Cref{def:uniqueness}.
\end{theorem}

The assumptions are quite natural (except probably \Cref{as:full_rank}); we will explain the intuition behind \Cref{as:full_rank} in the next sections.
We will also show that, using our approach, \Cref{as:maxrank_F,as:nonzero_GH} can be relaxed.
Moreover, all these assumptions are generic as long as the matrices have certain dimensions.

\begin{remark}[On genericity of assumptions]
\Crefrange{as:nonzero_F}{as:nonzero_GH} are satisfied for generic $\matr{F},\matr{G},\matr{H}$. i.e., the assumptions are satisfied except a set of $\matr{F},\matr{G},\matr{H}$ of Lebesgue measure $0$ \footnote{More precisely, the exceptional set is a semialgebraic subset of strictly smaller dimension} (equivalently all the assumptions are satisfied with probability $1$ for factors drawn from any absolutely continuous probability distribution.).
If $\edim{1} \ge \rkRow$, $\edim{2} \ge \rkCol$, then \Cref{as:AB_full_rank} is also satisfied for generic $\matr{A}$ and $\matr{B}$, while \Cref{as:full_rank} is satisfied generically if  $\edim{3} \ge \frac{\rkRow(\rkRow+1)}{2}\frac{\rkCol(\rkCol+1)}{2}$.
\end{remark}

\begin{remark}
Note that $\edim{1} \ge \rkRow$, $\edim{2} \ge \rkCol$ is a necessary condition for \Cref{as:AB_full_rank} while $\edim{3} 
\ge \frac{\rkRow(\rkRow+1)}{2}\frac{\rkCol(\rkCol+1)}{2}$ is necessary for
 \Cref{as:full_rank}.
 \end{remark}
 
\subsection{Known results: symmetrically weighted case}
The symmetrically weighted case is slightly different from the non-symmetric case: \Cref{thm:HL_nonsymmetric} does not apply, because \Cref{as:full_rank} is never satisfied.
Indeed, if $\matr{G} = \matr{H}$, then 
$\matr{\CoreVert} \kr \matr{\CoreVert} \kr  \matr{\CoreHoriz}  \kr  \matr{\CoreHoriz}  = \matr{\CoreVert}^{\kr4}$, whose rank is at most $\binom{R+3}{4}$ which is strictly less than $\binom{R+1}{2}^2$ (the bound in  \Cref{as:full_rank}).
Therefore,  \Cref{as:full_rank} is replaced by the following assumption.
\begin{assumption}\label{as:full_rank_symmetric}
It holds that $\matr{\CoreVert} = \matr{\CoreHoriz} $ and
the matrix $\matr{\CoreVert} \kr \matr{\CoreVert} \kr \matr{\CoreVert}  \kr  \matr{\CoreVert}$ has the maximal possible rank
\[
\rank{\matr{\CoreVert} \kr \matr{\CoreVert} \kr  \matr{\CoreVert}  \kr  \matr{\CoreVert}} = \binom{R+3}{4} = \frac{(R+3)(R+2)(R+1)R}{24}.
\]
\end{assumption}
Then the following result holds true.
\begin{theorem}[{\cite[Theorem 2]{HarsL96:uniqueness}}]\label{thm:HL_symmetric}
Let $\tens{T} = \ptd{\matr{A}}{\matr{B}}{\matr{\CoreFront}}{\matr{\CoreVert}}{\matr{\CoreVert}}$ be the symmetrically weighted PT2D (with $R= S$) whose factors satisfy \Crefrange{as:AB_full_rank}{as:nonzero_GH} and \Cref{as:full_rank_symmetric}.
Then the symmetrically weighted PT2D is unique in the sense of \Cref{def:uniqueness}, i.e. any alternative symmetrically weighted PT2D 
\[
\tens{T} =
\ptd{\widetilde{\matr{A}}}{\widetilde{\matr{B}}}{\widetilde{\matr{\CoreFront}}}{\widetilde{\matr{\CoreVert}}}{\widetilde{\matr{\CoreVert}}},
\]
is necessarily as in  \Cref{def:uniqueness} (with an extra constraint $\widetilde{\matr{\CoreVert}} = \widetilde{\matr{\CoreHoriz}}$ and $\matr{\Pi}_{\matr{A}} = \matr{\Pi}_{\matr{B}}$).
\end{theorem}

\begin{remark}[On genericity of the assumptions]
The   \cref{as:full_rank_symmetric} is again satisfied for generic $\matr{G}$ if  $\edim{3} \ge \frac{\rkRow(\rkRow+1)(\rkRow+2)(\rkRow+3)}{24}$. Note that  $\edim{3} \ge \frac{\rkRow(\rkRow+1)(\rkRow+2)(\rkRow+3)}{24}$ is a necessary condition for  \cref{as:full_rank_symmetric} to hold.
\end{remark}

In particular, \cref{thm:HL_symmetric} also applies to the fully symmetric P2TD (DEDICOM), and is used to establish uniqueness of PARAFAC-2 \cite{KierTB99:direct} (thanks to the correspondence in \Cref{ex:PT2_PF2}).
Up to the author's knowledge, \Cref{thm:HL_symmetric} yields the best best known results for uniqueness of  DEDICOM (and, as a consequence, of the PARAFAC-2 decomposition), , except a very particular case in \cite{BergK96:uniqueness}, where generic uniqueness is shown for $R=2$, $\edim{3}=4$ (see also references in \cite{RoalSCABCA22:AOADMM}).

\subsection{Preview of the results of the paper}
The drawback of the results of \cite{HarsL96:uniqueness} is that they do not come with an algorithm that can recover the factors of a given PT2D under the assumptions  of \Cref{thm:HL_nonsymmetric} and \Cref{thm:HL_symmetric}.
In this paper, we fill this gap by providing the algorithms for the non-symmetric PT2D and DEDICOM (the case of symmetrically weighted PT2D  is omitted, but can be also treated along the same lines).
Moreover, we show that some assumptions can be relaxed.

First of all \Cref{as:nonzero_GH} can be completely removed (if \Cref{as:full_rank} or \Cref{as:full_rank_symmetric} are satisfied). Second,  \Cref{as:maxrank_F} can be relaxed to the following assumption.
\begin{assumption}\label{as:krank_F}
 $\krank{{\matr{\CoreFront}}^{\T}} \ge 2$ and $\krank{{\matr{\CoreFront}}} \ge 2$.
\end{assumption}
\Cref{as:krank_F} actually means that  $\matr{\CoreFront}$ has neither proportional columns nor proportional rows.
Then, in particular, the following two results hold true.
\begin{proposition}[\Cref{thm:pt2d_uniqueness_nonsymmetric}]
Let $\tens{T}  = \ptd{\matr{A}}{\matr{B}}{\matr{F}}{\matr{G}}{\matr{H}} \in \RR^{\edim{1}\times \edim{2} \times \edim{3}}$,  whose factors satisfy  \Cref{as:AB_full_rank,as:nonzero_F,as:full_rank,as:krank_F}.
Then the $(R,S)$-PT2D of $\tens{T}$ is essentially unique.
Also, $(R,S)$ are minimal ranks (i.e., there does not exist $(R',S')$-PT2D with either $R' < R$ or $S' <S$).
\end{proposition}
\begin{proposition}[\Cref{thm:pt2d_uniqueness_symmetric}]
Let $\tens{T}  = \ptdsym{\matr{A}}{\matr{F}}{\matr{G}} \in \RR^{\edim{1}\times \edim{2} \times \edim{3}}$,  $\matr{F} \in \RR^{R \times R}$ symmetric, whose factors satisfy  \Cref{as:AB_full_rank,as:nonzero_F,as:full_rank_symmetric,as:krank_F}.
Then the  $R$-DEDICOM of $\tens{T}$ is essentially unique.
In addition, $R$ is the minimal possible rank for a DEDICOM of $\tens{T}$.
\end{proposition}

Our approach relies on constructing the so-called ``lifting'' matrix $\matr{\Phi}(\tens{T})$, whose left kernel has a particular structure. 
This enables us to retrieve the factors $\matr{A}$ and $\matr{B}$ from the left kernel of the lifting matrix (up to permutation and scaling of columns), uniquely, under the  assumptions mentioned above. 
The other factors ($\matr{F}$) can be retrieved uniquely (up to remaining ambiguities)  once $\matr{A}$ and $\matr{B}$ are fixed.

The lifting approach provides an algorithm to recover the factors  (${\matr{A}},{\matr{B}},{\matr{F}},{\matr{G}},{\matr{H}}$) for the PT2D and ${\matr{A}},{\matr{F}},{\matr{G}}$ for DEDICOM), up to trivial ambiguities, and also implies  uniqueness of the decompositions.
We also show that if  \Cref{as:krank_F} is not satisfied, this leads to additional ambiguities of the PT2D, namely  $\matr{A}$ (resp. $\matr{B}$) is non-unique if $\krank{{\matr{\CoreFront}}^{\T}} = 2$ (resp.  $\krank{{\matr{\CoreFront}}} \ge 2$).
However, in that case, the set of equivalent solutions can be described.

Although we do not cover the symmetrically weighted case in this paper,  it can be treated along the same lines by applying the lifting approach.


\section{PT2D and its core tensor}
In this section, we recall an alternative representation for the PT2D via the core tensor (see \cite{FaviA14:constrained}) and discuss the basic implications.
\subsection{Triple product of matrices}
First, we define the following special product of three matrices.
\begin{definition}[Triple product]\label{def:triprod}
For three matrices  $\matr{\CoreFront} \in \RR^{R \times S}$,
$\matr{\CoreVert} \in \RR^{R \times \edim{3}}$, $\matr{\CoreHoriz} \in \RR^{S \times \edim{3}}$, 
their triple product, denoted as $\tens{C} = \triprod{\matr{\CoreFront}}{\matr{\CoreVert}}{\matr{\CoreHoriz}}$,
 is the tensor $\tens{C} \in\RR^{R\times S\times K}$ whose elements are defined by
 \begin{equation}\label{eq:paratuck_core}
\tenselem{\Core}_{ijk} ={\CoreFront}_{ij}  {\CoreVert}_{ik} {\CoreHoriz}_{jk}.
\end{equation}
\end{definition}
Graphically, \Cref{def:triprod} can be visualised as a ``product'' of three matrices, see \Cref{fig:core_tensor}.
\begin{figure}[ht!]
\centering
{\begin{tikzpicture}
\begin{scope}[scale= 0.5]
\draw(-0.3,0.5) node (A) {\footnotesize $R$};
\draw(0.5,-0.3) node (A) {\footnotesize $S$};
\draw(-0.1,1.4) node (A) {\footnotesize $\edim{3}$};

\draw(0.5,0.5) node (A) {$\tens{C}$};
\draw(0,0) rectangle (1,1);
\draw(0,1) -- (0.5,1.5) -- (1.5,1.5) -- (1.5,0.5) --(1,0);
\draw(1,1) -- (1.5,1.5);

\draw(2,0.8) node (eq) {$=$};

\draw(3,0) -- (3,1) -- (3.5,1.5) -- (3.5,0.5) -- (3,0);
\draw(3.7,1.3) -- (4.7,1.3) -- (5.2,1.8) -- (4.2,1.8) -- (3.7,1.3);
\draw(3.8,-0.1) rectangle (4.8,0.9);

\draw(3.65,1.1) node (s22) {\small $\cdot$};

\draw(4.3,0.4) node (s22) {\small {\footnotesize$\matr{F}$}};
\draw(3.25,0.75) node (s23) {\small{\footnotesize$\matr{G}$}};
\draw(4.45,1.55) node (s23) {\small{\footnotesize$\matr{H}$}};
\end{scope}
\end{tikzpicture}}
\caption{Product of three matrices.}
\label{fig:core_tensor}
\end{figure}

\begin{remark}
If  $\tens{\Core} = \triprod{\matr{\CoreFront}}{\matr{\CoreVert}}{\matr{\CoreHoriz}}$, then its slices are equal to  
\begin{equation}\label{eq:hadamard_rank1}
\tens{\Core}_{:,:,k} =  ( \matr{\CoreFront} \hadam \vect{g}_k \vect{h}_k^{\T}), 
\end{equation}
where  $\vect{g}_k, \vect{h}_k$ are as in \eqref{eq:columns_core_factors}, and $\hadam$ is the Hadamard product (i.e., slices of $\tens{\Core}$ are Hadamard products of the same matrix $\matr{\CoreFront}$ with different rank-one matrices).
\end{remark}

\begin{remark}
The tensor factorization \eqref{eq:paratuck_core} is known in algebraic statistics as the ``no-three-way interaction'' model \cite{DiacS98:aos}.
It is a typical example of a tensor factorization  that cannot be realized as a graphical model. 
\end{remark}

\subsection{Factorizing the core tensor}\label{sec:factorizing_core}
We also remark that the factorization $\tens{\Core} = \triprod{\matr{\CoreFront}}{\matr{\CoreVert}}{\matr{\CoreHoriz}}$ admits trivial ambiguities
\begin{lemma}\label{lem:core_ambiguities}
A tensor $\tens{\Core} = \triprod{\matr{\CoreFront}}{\matr{\CoreVert}}{\matr{\CoreHoriz}}$ admits an alternative decomposition $\tens{\Core} = \triprod{\widetilde{\matr{\CoreFront}}}{\widetilde{\matr{\CoreVert}}}{\widetilde{\matr{\CoreHoriz}}}$ for all ${\widetilde{\matr{\CoreFront}}}, {\widetilde{\matr{\CoreVert}}}, {\widetilde{\matr{\CoreHoriz}}}$ as
\[
\widetilde{\matr{\CoreFront}} = \matr{\Lambda}^{-1}_{\matr{G}} \matr{\CoreFront}\matr{\Lambda}^{-1}_{\matr{H}},\quad
\widetilde{\matr{\CoreVert}} = \matr{\Lambda}_{\matr{G}} \matr{\CoreVert} \matr{\Lambda}_{\matr{GH}},\quad
\widetilde{\matr{\CoreHoriz}} = \matr{\Lambda}_{\matr{H}} \matr{\CoreHoriz}\matr{\Lambda}_{\matr{GH}}^{-1},
\]
where $\matr{\Lambda}_{\matr{G}} \in \RR^{R\times R}$, $\matr{\Lambda}_{\matr{H}} \in \RR^{S\times S}$, $\matr{\Lambda}_{\matr{GH}} \in \RR^{\edim{3}\times \edim{3}}$ are any diagonal nonsingular matrices.
\end{lemma}
\begin{proof}
Follows from straightforward calculation.
\end{proof}
The basic question is, whether for a  tensor $\tens{C}$ as in  \eqref{eq:paratuck_core}, we can recover the factors $\matr{F},\matr{G},\matr{H}$, given $\tens{C}$.
We show that $\tens{\Core}$ can be uniquely factorized (up to the ambiguities in \Cref{lem:core_ambiguities}) under some additional conditions.
We propose the following two conditions that lead to two different ways to compute the factors.

\begin{lemma}\label{lem:nonzero_slice}
Let $\tens{C}$ be of the form \eqref{eq:paratuck_core}.
If there  exists  one  index $r$ such that $\tenselem{\Core}_{ijr} \neq 0$ for any $i,j$, then $\matr{F}, \matr{G},\matr{H}$ can be uniquely recovered from $\tens{\Core}$ (up to scaling ambiguities in \Cref{lem:core_ambiguities} ).
\end{lemma}

\begin{remark}
In terms of assumptions  on factors, the condition of \Cref{lem:nonzero_slice} is equivalent to requiring that the \Cref{as:nonzero_F} ($F_{ij} \neq 0$) holds together with $\vect{g}_r$ and $\vect{h}_r$  not having zero elements for a fixed $r$ ( i.e., $\DiagA{r}$ and $\DiagB{r}$ are nonsingular for fixed $r$).
\end{remark}

However, the factors of the core tensors can be recovered even if none of   $\DiagA{r}$ and $\DiagB{r}$ are simultaneously nonsingular, but under an additional assumption.

\begin{proposition}\label{prop:core_factorization}
Let $\tens{C}$ be of the form \eqref{eq:paratuck_core}, whose factors $\matr{\CoreFront}$, $\matr{\CoreVert}$, $\matr{\CoreHoriz}$ satisfy \Cref{as:nonzero_F} and either of  \Crefrange{as:full_rank}{as:full_rank_symmetric}.
Then  $\matr{\CoreFront}$, $\matr{\CoreVert}$, $\matr{\CoreHoriz}$  can be uniquely recovered (up to the scaling ambiguities in \Cref{lem:core_ambiguities} ) from $\matr{C}$.
\end{proposition}

The proofs of \Cref{lem:nonzero_slice,prop:core_factorization} are constructive (i.e., provide  algorithms to compute the factorization) and are given in \Cref{sec:factorizing_core_proofs}.

\subsection{PT2D through the triple product}
It is easy to show that the PT2D can be defined via the triple product of matrices.
\begin{lemma}\label{lem:paratuck_contraction}
For any $\matr{A}$, $\matr{B}$, $\matr{F}$, $\matr{G}$, $\matr{H}$ as in  \Cref{def:PT2D} we have
\begin{equation}\label{eq:paratuck_contraction}
 \ptd{\matr{A}}{\matr{B}}{\matr{\CoreFront}}{\matr{\CoreVert}}{\matr{\CoreHoriz}} =
\Big(\triprod{\matr{\CoreFront}}{\matr{\CoreVert}}{\matr{\CoreHoriz}} \Big) \con_1 \matr{A} \con_2 \matr{B}.
\end{equation}
\end{lemma}
\begin{proof}
Indeed if  $\tens{\Core} = \triprod{\matr{\CoreFront}}{\matr{\CoreVert}}{\matr{\CoreHoriz}}$, then it is easy to see from \eqref{eq:hadamard_rank1} that its slices satisfy
\begin{equation}\label{eq:core_slices}
\tens{\Core}_{:,:,k} = \DiagA{k} \matr{\CoreFront}  \DiagB{k},
\end{equation}
with $\DiagA{k}, \DiagB{k}$  as in \cref{eq:DiagAB}.
Injecting  \eqref{eq:core_slices} in \eqref{eq:paratuck_contraction} and comparing with \eqref{eq:paratuck_slices} we see that \eqref{eq:paratuck_contraction} is verified.
\end{proof}
This is a building block in computing the PT2D, as this corresponds to the subproblem of finding the remaining factors of the PT2D when $\matr{A}$ and $\matr{B}$ are known (see the next subsection for the details).
\begin{remark}
The ambiguities in \Cref{lem:core_ambiguities} correspond exactly to ambiguities of PT2D in \Cref{def:uniqueness} when $\matr{A}$ and $\matr{B}$ are fixed.
\end{remark}
\Cref{lem:paratuck_contraction} shows that a PT2D tensor admits a Tucker-2 factorization $\tens{\Core}\con_1 \matr{A} \con_2 \matr{B}$ with the core tensor structured as \eqref{eq:paratuck_core}.
In the next subsection we show how to use the Tucker-2 factorization to reduce dimension of the problem.

\subsection{Multilinear ranks in the PT2D}\label{sec:ml_ranks}
Next, let us discuss how to reduce the PT2D decomposition problem from the case of rectangular factors $\matr{A}$ and $\matr{B}$  ($\edim{1} > R$ or $\edim{2} >S$) to the case when  $\matr{A}$ and $\matr{B}$  are square (and invertible).
For this, we use a standard strategy which performs Tucker compression first and relies on the following lemma.

\begin{lemma}\label{lem:pt2d_tucker2}
Let $\tens{T}$ have an $(R,S)$-PT2D with $R \le\edim{1} $ or $S \le \edim{2}$.
\begin{enumerate}
\item It holds that
\begin{equation}\label{eq:ml_rank_bound}
\rank{\matr{T}^{(1)}} \le R,\quad\rank{\matr{T}^{(2)}} \le S.
\end{equation}
\item Assume that $\rank{\matr{T}^{(1)}} = R,\rank{\matr{T}^{(2)}} = S$, and   $\matr{U} \in \RR^{\edim{1} \times R}$, $\matr{W} \in \RR^{\edim{2} \times S}$ be the bases of the corresponding column spaces (i.e., $\range{\matr{U} } = \range{\matr{T}^{(1)}}$,  $\range{\matr{W} } = \range{\matr{T}^{(2)}}$).

Then  $\tens{T} = \ptd{\matr{A}}{\matr{B}}{\matr{F}}{\matr{G}}{\matr{H}}$ is an $(R,S)$-PT2D of $\tens{T}$ if and only if the compressed tensor $\tens{T}_c = \tens{T} \con_1 \matr{U}^{\dagger} \con_2 \matr{V}^{\dagger}$ has the PT2D
\begin{equation}\label{eq:pt2d_compressed}
\tens{T}_c   = \ptd{\matr{A}_c}{\matr{B}_c}{\matr{F}}{\matr{G}}{\matr{H}};
\end{equation}
Also, in \eqref{eq:pt2d_compressed}, the factors of the compressed and full PT2D are linked as $\matr{A} = \matr{U}\matr{A}_c$, $\matr{B} = \matr{V}\matr{B}_c$.
\end{enumerate}
\end{lemma}
\begin{proof}
\begin{enumerate}
\item From \cref{eq:paratuck_slices}, we have that the $\range{\tens{T}_{:,:,k}}  \subseteq\range{\matr{A}}$ and $\range{(\tens{T}_{:,:,k})^{\T}}  \subseteq \range{\matr{B}}$, which implies that  $\range{\matr{T}^{(1)}} \subseteq\range{\matr{A}}$,  $\range{\matr{T}^{(2)}} \subseteq\range{\matr{B}}$, hence the rank condition is satisfied.

\item If $\rank{\matr{T}^{(1)}} = R,\rank{\matr{T}^{(2)}} = S$, then it must hold that  $\range{\matr{T}^{(1)}} =\range{\matr{A}}$ and  $\range{\matr{T}^{(2)}} = \range{\matr{B}}$, for any valid rank-$(R,S)$ PT2D.
In particular, there exist nonsingular $\matr{A}_c \in \RR^{R \times R}$ and $\matr{B}_c \in \RR^{S \times S}$, such that  $\matr{A} = \matr{U}\matr{A}_c$ and $\matr{B} = \matr{V}\matr{B}_c$.
Then by multiplying $\tens{T}_{:,:,k}$ by pseudoinverses of $\matr{U}$ and $\matr{V}$, we get
that the slices of the compressed tensor $\tens{T}_c = \tens{T} \con_1 \matr{U}^{\dagger} \con_2 \matr{V}^{\dagger}$  can be expressed as

\[
(\matr{U})^{\dagger} \tens{T}_{:,:,k} (\matr{V}^{\T})^{\dagger} =  \matr{A}_c \DiagA{k} \matr{\CoreFront}  \DiagB{k}  \matr{B}_c^{\T},
\]
hence $\tens{T}_c$ has a PT2D \eqref{eq:pt2d_compressed}.
Vice versa, any PT2D of $\tens{T}_c$ in \eqref{eq:pt2d_compressed} gives rise to the PT2D of $\tens{T}$ since $\tens{T} = \tens{T}_c \con_1 \matr{U} \con_2 \matr{V}$.
\end{enumerate}
\end{proof}

In the rest of the subsection, we will show when the first and second  unfoldings of  $\tens{T}$ are guaranteed to have maximal possible rank.
We first propose some remarks on the ranks of undoldings of $\tens{C}$.
\begin{proposition}\label{prop:core_multilinear_rank}
Let $\matr{F}$ satisfy  \Cref{as:nonzero_F}. Then the following holds true.
\begin{enumerate}
\item If  there exists $j$ such that $\rank{\matr{G} \kr \matr{H}_{j,:}} = R$, then $\rank{\matr{C}^{(1)}} = R$.
\item If  there exists $i$ such that $\rank{\matr{G}_{i,:} \kr \matr{H}} = S$, then $\rank{\matr{C}^{(2)}} = S$.
\end{enumerate}
\end{proposition}
\begin{proof}
See \Cref{sec:proofs_ml_rank} for the proofs.
\end{proof}

Next, we can formulate the following simple corollary  for the multilinear ranks of $\tens{T}$.

\begin{corollary}\label{prop:pt2d_multilinear_rank}
If both conditions of \Cref{prop:core_multilinear_rank} are satisfied together with \Cref{as:AB_full_rank}, then
\[
\rank{\matr{T}^{(1)}} = R, \quad \rank{\matr{T}^{(2)}} = S.
\]
\end{corollary}
\begin{proof}
Follows from the fact that $\matr{T}^{(1)} = \matr{A} \matr{C}^{(1)},\quad \matr{T}^{(2)} = \matr{B} \matr{C}^{(2)}$, where both matrices in the factorizations have rank $R$ and $S$ respectively by \Cref{prop:core_multilinear_rank}  together with \Cref{as:AB_full_rank}.
\end{proof}

Finally, we remark that the conditions of \Cref{prop:core_multilinear_rank} are, in fact, satisfied under   \Cref{as:full_rank} or \Cref{as:full_rank_symmetric}, as shown by the following corollary.
\begin{corollary}\label{cor:core_ml_rank_assumptions}
Under    \Cref{as:nonzero_F} and either of  \Cref{as:full_rank} or \Cref{as:full_rank_symmetric}, we have  
\[
\rank{\matr{C}^{(1)}} = R,\quad\rank{\matr{C}^{(2)}} = S.
\]
\end{corollary}
\begin{proof}
See \Cref{sec:proofs_ml_rank} for the proof.
\end{proof}

\subsection{Collinearities in $F$ and nonuniqueness}\label{sec:collinearities}
In this section, we explain why \Cref{as:krank_F} is necessary for uniqueness of the PT2D and why 
 additional indeterminacies can arise in  case \Cref{as:krank_F} is not satisfied.
This subsection can be skipped if the reader is not interested in nonunique cases.
 We start from a simple example to build an intuition.
 \begin{example}\label{ex:F_rank_1}
 Let $\rank{\matr{\CoreFront}} = 1$, i.e., all the rows are collinear and the columns as well.
 Then thanks to \Cref{lem:ambiguities}, we can assume, without loss of generality that 
 $\matr{\CoreFront} = \matr{1}\matr{1}^{\T}$.
 Then for the tensor $\tens{T} = \ptd{\matr{A}}{\matr{B}}{\matr{F}}{\matr{G}}{\matr{H}}$ decomposition of its slices (from \eqref{eq:hadamard_rank1} or \eqref{eq:paratuck_slices})
 reads simply
 \[
 \tens{T}_{:,:,k} =  \matr{A} (\vect{g}_k \vect{h}_k^{\T}) \matr{B}^{\T}.
 \]
 Let us take any nonsingular matrices $\matr{Y} \in \RR^{R\times R}$,  $\matr{Z} \in \RR^{S \times S}$. Then if we define
 \[
 \widetilde{\matr{A}} =  \matr{A} \matr{Y},\quad  \widetilde{\matr{B}} =  \matr{B} \matr{Z}, \quad
 \widetilde{\vect{g}}_k =  \matr{Y}^{-1} \vect{g}_k, \quad  \widetilde{\vect{h}}_k =  \matr{Z}^{-1} \vect{h}_k, 
 \] 
 we get that 
  \[
 \tens{T}_{:,:,k} = \widetilde{\matr{A}}  (\widetilde{\vect{g}}_k   \widetilde{\vect{h}}_k^{\T} )  \widetilde{\matr{B}}^{\T},
 \]
 hence $\tens{T} = \ptd{\widetilde{\matr{A}}}{\widetilde{\matr{B}}}{\matr{F}}{\widetilde{\matr{G}}}{\widetilde{\matr{H}}}$ with ${\widetilde{\matr{G}}}$ and ${\widetilde{\matr{H}}}$ defined from $\widetilde{\vect{g}}_k$, $\widetilde{\vect{h}}_k^{\T}$.
 \end{example}
 \Cref{ex:F_rank_1} means that for the rank-1 $\matr{F}$, the  factors $\matr{A}$ and $\matr{B}$ can be recovered uniquely in the PT2D only up to the subspaces spanned by their columns.
In the general case ($\krank{\matr{F}} =1$ or $\krank{\matr{F}^{\T}} =1$) we will show that a similar statement can be made for blocks of columns of  $\matr{A}$ and $\matr{B}$.

Let $\krank{\matr{F}} = 1$. Then, without loss of generality, we can assume that the columns of the matrix $\matr{F}$ consists of blocks of repeating columns of sizes $n_1, \ldots,n_{S'}$, i.e
\[
\matr{F} = [\overbrace{\vect{f}_{1}  \cdots \vect{f}_{1}}^{n_1 \text{ times}}  \cdots
 \overbrace{\vect{f}_{S'}  \cdots \vect{f}_{S'}}^{n_{S'} \text{ times}} ]
\] 
with $n_1 + \cdots + n_{S'} = R$ and all $\vect{f}_k$ are non-collinear ($\krank{\begin{bmatrix}
\vect{f}_{1} \ & \vect{f}_{2} & \cdots & \vect{f}_{S'}
\end{bmatrix}} \ge 2$).
Indeed, we can always reorder the rows to bring the collinear columns together and rescale the columns of $\matr{F}$ such the collinear columns are represented by the same vector.
This can be summarized in the following lemma.
\begin{lemma}\label{lem:F_repeating_columns}
In PT2D with  $\matr{F}$ satisfying \Cref{as:nonzero_F} (and $\rank{\matr{F}} > 1$),  we can assume that there exists $\matr{F}' \in \matr{F}^{R' \times S'}$, $R'\le R,S'\le S$, and $\krank{\matr{F}'} \ge 2$, $\krank{(\matr{F}')^{\T}} \ge 2$,  so that $\matr{F}$ can be represented\footnote{i.e.,  any $\matr{F}$ satisfying \Cref{as:nonzero_F} can be brought to the form \eqref{eq:F_repeating_columns} by permutations and rescalings of columns and rows.} as
\begin{equation}\label{eq:F_repeating_columns}
\matr{F} = 
\begin{bmatrix}
\matr{1}_{m_1} &  & & \\
 &  \matr{1}_{m_2} & & \\
 &&\ddots&\\
 &  & & \matr{1}_{m_R'} 
\end{bmatrix}
\matr{F}' 
\begin{bmatrix}
\matr{1}^{\T}_{n_1} &  & & \\
 &  \matr{1}^{\T}_{n_2} & & \\
 &&\ddots&\\
 &  & & \matr{1}^{\T}_{n_S'} 
\end{bmatrix}
\end{equation}
where $m_1+ \cdots + m_{R'} = R$ and $n_1+ \cdots + n_{S'} = S$.
\end{lemma}
Then for $\matr{F}$ in a such form, we can describe the ambiguities that arise from the collinearities.
\begin{proposition}
Let $\matr{F}$ be as in \Cref{lem:F_repeating_columns}, and $\tens{T} = \ptd{\matr{A}}{\matr{B}}{\matr{F}}{\matr{G}}{\matr{H}}$.
Then $\tens{T}$ admits the family of following alternative PT2D (in addition to those that can be obtained in  \Cref{lem:ambiguities}):
\[
\tens{T} = \ptd{\matr{A}\matr{Y}}{\matr{B}\matr{Z}}{\matr{F}}{\matr{Y}^{-1}\matr{G}}{\matr{Z}^{-1}\matr{H}},
\]
where $\matr{Y}$ and $\matr{Z}$ are block-diagonal matrices
\[
\matr{Y} = \Diag{\matr{Y}_{1} ,\ldots,\matr{Y}_{R'}} = \begin{bmatrix}
\matr{Y}_{1} &  & & \\
 &  \matr{Y}_{2} & & \\
 &&\ddots&\\
 &  & & \matr{Y}_{R'} 
\end{bmatrix},\quad \matr{Z} =  \Diag{\matr{Z}_{1} ,\ldots,\matr{Z}_{S'}} =\begin{bmatrix}
\matr{Z}_{1} &  & & \\
 &  \matr{Z}_{2} & & \\
 &&\ddots&\\
 &  & & \matr{Z}_{S'} 
\end{bmatrix},
\]
with $\matr{Y}_k \in \RR^{m_k \times m_k}$ and $\matr{Z}_k \in \RR^{n_k \times n_k}$ are all invertible.
\end{proposition}
\begin{proof}
The proof is analogous to the derivations in \Cref{ex:F_rank_1}.
\end{proof}


\section{Lifting approach to PT2D: the non-symmetric case}
In the next two subsections, we describe the main tool used in this paper, and present the main results for the non-symmetric PT2D.
Our goal will be to determine $\matr{A}$ and $\matr{B}$, since when these matrices are known, then we can find the other PT2D factors thanks to the results of the previous subsection.
In the rest of the section, in view of \cref{sec:ml_ranks}, we assume that $\edim{1} = R$ and $\edim{2} =S$, hence the matrices $\matr{A}$ and $\matr{B}$ are square invertible.

\subsection{Lifting and the structured matrix: basic properties}
\begin{definition}\label{def:PhiT}
For $\tens{T} \in \RR^{R \times S \times \edim{3}}$, we define the following $(RS)^2 \times \edim{3}$ structured matrix:
\[
\matr{\Phi}(\tens{T})
= 
\begin{bmatrix}
\vecl{\tens{T}_{:,:,1}} \kron \vecl{\tens{T}_{:,:,1}} & \cdots &
\vecl{\tens{T}_{:,:,\edim{3}}} \kron \vecl{\tens{T}_{:,:,\edim{3}}}
\end{bmatrix} \in \RR^{(RS)^2 \times \edim{3}},
\]
or equivalently, $\matr{\Phi}(\tens{T})$ can be expressed as $\matr{\Phi}(\tens{T}) = (\unfold{T}{3})^{\T}\kr(\unfold{T}{3})^{\T}$.
\end{definition}
\begin{remark}\label{rem:Phi_symmetry}
Note that the columns of $\matr{\Phi}(\tens{T})$, in fact, correspond 4D tensors $\{\tens{T}_{:,:,k} \otimes \tens{T}_{:,:,k}\}_{k=1}^{\edim{3}} \subset \RR^{R \times S \times R \times S}$ reshaped in a special way. This implies that $\range{\matr{\Phi}(\tens{T})}$ consists of  vectorized symmetric $RS \times RS$ matrices, hence
 $\range{\matr{\Phi}(\tens{T})} \subset \vecl{S^2(\RR^{RS})}$ for any $\tens{T}$.
\end{remark}

Now let us look at the properties of  $\matr{\Phi}(\tens{T})$ when $\tens{T}$ admits a PT2D, and particularly at the  the left nullspace of $\matr{\Phi}(\tens{T})$.
Due to \Cref{rem:Phi_symmetry} only  the symmetric part of the left kernel is informative.
To this end, we will consider the symmetrization  of the left nullspace $\symlker{\Phi{(\tens{T})}}$, which is, by \Cref{rem:Phi_symmetry}, equal to
\[
\symlker{\Phi{(\tens{T})}} = \lker{\Phi{(\tens{T})}} \cap \vecl{S^{2}(\RR^{RS})}.
\]
Then the following key lemma summarizes the properties of rank and symmetric left nullspace properties for ${\matr{\Phi}(\tens{T})}$ and also explains the importance of  \Cref{as:full_rank}.

\begin{lemma}\label{lem:Phi_rank_bound_nonsymmetric}
\begin{enumerate}
\item 
 For $\tens{T} = \ptd{\matr{A}}{\matr{B}}{\matr{F}}{\matr{G}}{\matr{H}} \in \RR^{R\times S\times \edim{3}}$ with $\matr{A} \in \RR^{R\times R}$, $\matr{B} \in \RR^{S\times S}$ we have
\begin{equation}\label{eq:Phi_rank_bound}
\rank{\matr{\Phi}(\tens{T})} \le \binom{R+1}{2} \binom{S+1}{2}.
\end{equation}
\item  If
all \Cref{as:AB_full_rank,as:nonzero_F,as:full_rank} are satisfied, then equality is achieved in \eqref{eq:Phi_rank_bound}.
In this case, $\symlker{\matr{\Phi}(\tens{T})}$ spanned by $M =  \binom{R}{2} \binom{S}{2}$ vectors
\[
\symlker{\matr{\Phi}(\tens{T})} = \Span{\vect{p}_1,\ldots,\vect{p}_M}.
\]
\item If one of \Cref{as:nonzero_F,as:full_rank} is not satisfied, we have a strict inequality in \eqref{eq:Phi_rank_bound}.
\end{enumerate}
 \end{lemma}
\begin{proof}
The proof  is presented later in \Cref{sec:nonsymmetric_proofs}.
\end{proof}

\subsection{The main result for non-symmetric case}
The main result  relies on \Cref{lem:Phi_rank_bound_nonsymmetric} and needs an extra bit of notation.
Let $\symperm: \vecl{\RR^{R \times S \times R\times S}} \to \RR^{R^2 \times S^2}$ be the permutation-symmetrization-matricization map that is defined on rank-one tensors as follows:
\begin{equation}\label{eq:symperm}
\symperm: \vecl{\vect{a} \tensp \vect{b} \tensp  \vect{y} \tensp \vect{z}} \mapsto 
\frac{1}{2}(\vect{a} \kron \vect{y} + \vect{y} \kron \vect{a}) \cdot \frac{1}{2}(\vect{b} \kron \vect{z} + \vect{z} \kron \vect{b})^{\T}, \quad \mbox{for all } \vect{y}, \vect{a} \in \RR^{R}, \vect{z}, \vect{b} \in \RR^{S},
\end{equation}
that is $\symperm$ permutes modes $2$ and $3$ of the 4-way tensor, performs symmetrization along the modes of the same size, and flattens the tensor into an $R^2 \times S^2$ matrix.
Then the key result for the nonsymmetric PT2D can be formulated as follows.

\begin{theorem}\label{thm:pt2d_nonsymmetric}
Let  $\tens{T} = \ptd{\matr{A}}{\matr{B}}{\matr{F}}{\matr{G}}{\matr{H}} \in \RR^{R\times S\times \edim{3}}$, whose factors satisfy  \Cref{as:AB_full_rank,as:nonzero_F,as:full_rank} and let
$\matr{P}_k = \symperm(\matr{p}_k)$ where $\vect{p}_k$ be any basis of $\symlker{\matr{\Phi}(\tens{T})}$ (as in \Cref{lem:Phi_rank_bound_nonsymmetric}).
\begin{enumerate}
\item Let $\matr{P}_A = \begin{bmatrix}\matr{P}_1 &\cdots & \matr{P}_M\end{bmatrix} \in \RR^{R^2 \times S^2K}$. Then $\dim(\symlker{\matr{P}_A})  = R$ if and only if $\krank{\matr{F}^{\T}} \ge 2$.
In that case
\[
\symlker{\matr{P}_A} = \Span{ \vect{a}_1 \kron \vect{a}_1, \ldots,\vect{a}_R \kron \vect{a}_R } = \range{\matr{A}\kr\matr{A}},
\]
where $\vect{a}_k$ are the columns of $\matr{A}$.
\item Let $\matr{P}_B = \begin{bmatrix}\matr{P}_1^{\T} &\cdots & \matr{P}_M^{\T}\end{bmatrix}\in \RR^{S^2 \times R^2K}$. Then $\dim(\symlker{\matr{P}_B})  = S$ if and only if $\krank{\matr{F}} \ge 2$.
In that case
\[
\symlker{\matr{P}_B} = \Span{ \vect{b}_1 \kron \vect{b}_1, \ldots,\vect{b}_S \kron \vect{b}_S } = \range{\matr{B}\kr\matr{B}},
\]
where $\vect{b}_k$ are the columns of $\matr{B}$.
\end{enumerate}
\end{theorem}
\begin{proof}
The proof will be split into several lemmas and presented later in a  \Cref{sec:nonsymmetric_proofs}.
\end{proof}
Then \Cref{thm:pt2d_nonsymmetric} implies the following improved result on identifiability of PT2D.

\begin{proposition}\label{thm:pt2d_uniqueness_nonsymmetric}
Let $\tens{T}  = \ptd{\matr{A}}{\matr{B}}{\matr{F}}{\matr{G}}{\matr{H}} \in \RR^{\edim{1}\times \edim{2} \times \edim{3}}$,  whose factors satisfy  \Cref{as:AB_full_rank,as:nonzero_F,as:full_rank,as:krank_F}.
Then the $(R,S)$-PT2D of $\tens{T}$ is essentially unique.
Also, $(R,S)$ are minimal ranks (i.e., there does not exist $(R',S')$-PT2D with either $R' < R$ or $S' <S$).
\end{proposition} 
\begin{proof}
First, note that by \Cref{as:AB_full_rank,as:nonzero_F,as:full_rank} we have that $\rank{\matr{T}^{(1)}} = R$ and $\rank{\matr{T}^{(2)}} = S$, due to \Cref{lem:pt2d_tucker2}.
This implies that there is no $(R',S')$-PT2D with either $R' < R$ or $S' <S$ (otherwise we would have $\rank{\matr{T}^{(1)}} = R' < R$ or $\rank{\matr{T}^{(2)}} = S' <S$, a contradiction).
Next, due to \Cref{lem:pt2d_tucker2} we can assume, without loss of generality that $\edim{1} = R$ and $\edim{2} = S$.
If $\edim{1} > R$ and $\edim{2} > S$, then by  \Cref{as:AB_full_rank,as:nonzero_F,as:full_rank} we have $\rank{\matr{T}^{(1)}} = R$ and $\rank{\matr{T}^{(2)}} = S$, thus by \Cref{lem:pt2d_tucker2}, the  $(R,S)$-PT2D of $\tens{T}$ is unique if and only if the $(R,S)$-PT2D of the compressed tensor $\tens{T}_c\in \RR^{R\times S\times K}$ is unique.

Now, let us assume $\edim{1} = R, \edim{2} =S$. Under \Cref{as:AB_full_rank,as:nonzero_F,as:full_rank}, by  \Cref{lem:Phi_rank_bound_nonsymmetric} we have $\rank{\matr{\Phi}(\tens{T})} = \binom{R+1}{2} \binom{S+1}{2}$.
If \Cref{as:krank_F} is also satisfied, we have  $\dim(\symlker{\matr{P}_A})  = R$ and $\dim(\symlker{\matr{P}_B})  = S$ by \Cref{thm:pt2d_nonsymmetric}.

Let $\matr{q}_{A,1}, \ldots, \matr{q}_{A,R} \in \RR^{R^2}$ be some basis of $\symlker{\matr{P}_A}$. Then, by \Cref{thm:pt2d_nonsymmetric} there is a nonsingular (change of basis) matrix $\matr{W}_A \in \RR^{R\times R}$
\[
\begin{bmatrix}\matr{q}_{A,1} & \ldots & \matr{q}_{A,R} \end{bmatrix} = (\matr{A} \kr \matr{A}) (\matr{W}_A)^{\T},
\]
This implies that the tensor $\tens{Q}_{A} \in \RR^{R\times R \times R}$ built as $(\tens{Q}_{A})_{:,:,k} = \matricize{R}{R}{\matr{q}_{A,k}}$ admits a CPD
\[
\tens{Q}_{A} = \cpd{\matr{A}}{\matr{A}}{\matr{W}_A}.
\]
Note that, by our assumptions, $\krank{\matr{A}} = \krank{\matr{W}_A} = R$, and therefore the Kruskal's uniqueness condition is satisfied \cite{KoldB09:tensor}, and the matrix $\matr{A}$ is unique up to permutation and scaling of columns.

Similarly, we can construct $\tens{Q}_{B} \in \RR^{S\times S \times S}$ built as ${(\tens{Q}_{B})_{:,:,k}} = \matricize{S}{S}{\matr{q}_{B,k}}$ (with $\{\matr{q}_{B,k}\}_{k=1}^S$ be any basis of  $\symlker{\matr{P}_B}$). For this tensor,  we have that, for some nonsingular $\matr{W}_2 \in \RR^{S\times S}$,
\[
\tens{Q}_{B} = \cpd{\matr{B}}{\matr{B}}{\matr{W}_B},
\]
whose CPD is unique by Kruskal's theorem and $\matr{B}$ is unique up to permutation and scaling of columns.
Thus, we can find ${\tens{C}} = \tens{T} \con_1 \matr{A}^{-1} \con \matr{B}^{-1}$ which must be of the form \eqref{eq:paratuck_core}.
The factors $\matr{F},\matr{G},\matr{H}$ can be recovered from $\tens{C}$ up to scaling indeterminacies, thanks to \Cref{prop:core_factorization}.

Finally, note that for any alternative $(R,S)$-PT2D,  $\tens{T}  = \ptd{\widetilde{\matr{A}}}{\widetilde{\matr{B}}}{\widetilde{\matr{F}}}{\widetilde{\matr{G}}}{\widetilde{\matr{H}}}$ we have that:
we have
\begin{itemize}
\item \Cref{as:AB_full_rank} for $(\widetilde{\matr{A}},\widetilde{\matr{B}})$ must be satisfied since $\rank{\matr{T}^{(1)}} = R$ and $\rank{\matr{T}^{(2)}} = S$.
\item \Cref{as:nonzero_F,as:full_rank} must be satisfied as well from statement 3) of  \Cref{lem:Phi_rank_bound_nonsymmetric};

\item  and  \Cref{as:krank_F} must be satisfied due to ``if and only if'' statements in   and \Cref{thm:pt2d_nonsymmetric}.
\end{itemize}
Therefore, (essential) uniqueness of the factors $\matr{A}$ and $\matr{B}$ in the CPDs of  $\tens{Q}_{A}$ and $\tens{Q}_{B}$ implies uniqueness of $\matr{A}$ and $\matr{B}$ in $(R,S)$-PT2D.
\end{proof}

\subsection{Algorithm for nonsymmetric case}
The proof of \Cref{thm:pt2d_uniqueness_nonsymmetric,thm:pt2d_nonsymmetric} suggest a way to compute the (unique under \Cref{as:AB_full_rank,as:nonzero_F,as:krank_F,as:full_rank}) PT2D of a given tensor, and we summarize it in \Cref{alg:paratuck_2x2_complete}.
We also show that the assumptions on factors  (\Cref{as:AB_full_rank,as:nonzero_F,as:krank_F,as:full_rank}) can be replaced by the assumptions on  the tensor $\tens{T}$ which can be checked \emph{a posteriori}.
\begin{algorithm}[hbt!]
\caption{Algebraic decomposition algorithm for nonsymmetric ParaTuck-2}\label{alg:paratuck_2x2_complete}
\begin{algorithmic}[1]
\STATE{{\bf Input:}$\tens{T} \in \RR^{\edim{1} \times \edim{2}\times \edim{3}}$ ranks $(R,S)$}
\STATE{Compute Tucker-2 decomposition: $\tens{T} = \tens{T}_c\con_1 \matr{U}\con_2 \matr{V}$ with $\matr{U} \in \RR^{\edim{1} \times R}$, $\matr{V} \in \RR^{\edim{2} \times S}$}
\STATE{Find the $M = \binom{R}{2} \binom{S}{2}$ vectors $\vect{p}_1,\ldots,\vect{p}_{M}$ vector in  symmetric left kernel of $\symlker{\matr{\Phi}(\tens{T}_c)}$}
\STATE{Find the projections $\matr{P}_k = \pi(\vect{p}_k)$, build the matrices ${\Mreshaped{A}}$ and ${\Mreshaped{B}}$.}
\STATE{Compute   bases $\{\matr{q}_{A,1},\ldots,\matr{q}_{A,R}\}$ and $\{\matr{q}_{B,1},\ldots,\matr{q}_{B,R}\}$ of the  symmetric left kernels $\symlker{\Mreshaped{A}}$ and $\symlker{\Mreshaped{B}}$, respectively.}
\STATE{Stack these matrices in tensors $\widehat{\tens{Q}}_A$ and $\widehat{\tens{Q}}_B$ and compute their CPD

$\widehat{\tens{Q}}_A = \cpd{\widehat{\matr{A}}_c}{\widehat{\matr{A}}_c}{{\matr{W}_A}} \quad \text{and} \quad 
\widehat{\tens{Q}}_B = \cpd{\widehat{\matr{B}}_c}{\widehat{\matr{B}}_c}{{\matr{W}_B}}$.
}
\STATE{Set  $\widehat{\matr{A}} = \matr{U} \widehat{\matr{A}}_c, \widehat{\matr{B}} = \matr{V}\widehat{\matr{B}}_c $. }
\STATE{Find the core tensor $\widehat{\tens{\Core}} = \tens{T}_c \contr{1} \widehat{\matr{A}}_c^{\dagger} \contr{2} \widehat{\matr{B}}_c^{\dagger}$. }
\STATE{Determine the factors $\widehat{\matr{\CoreFront}}$. $\widehat{\matr{\CoreVert}}$ and $\widehat{\matr{\CoreHoriz}}$ using the methods from \Cref{sec:factorizing_core} (namely, \Cref{prop:core_factorization}).}
\RETURN{$\widehat{\matr{A}},\widehat{\matr{B}},\widehat{\matr{F}},\widehat{\matr{G}},\widehat{\matr{H}}$, such that $\tens{T} = \ptd{\widehat{\matr{A}}}{\widehat{\matr{B}}}{\widehat{\matr{F}}}{\widehat{\matr{G}}}{\widehat{\matr{H}}}$}
\end{algorithmic}
\end{algorithm}

\begin{remark}
\Cref{alg:paratuck_2x2_complete} can be also applied in the approximation scenario, where all the exact decomposition steps are replaced by approximations.
For example, the nullspace computations can be performed by using the SVD, and the CPD can be replace by CP approximation.
\end{remark}

\begin{remark}
The algorithm can be  modified to handle the case when \Cref{as:full_rank} is not satisfied. 
In this case, the nullspaces $\symlker{\Mreshaped{A}}$, $\symlker{\Mreshaped{B}}$ will become larger and instead of the CPD of $\widehat{\tens{Q}}_A$, $\widehat{\tens{Q}}_B$ we would need to perform the symmetric joint block diagonalization. 
See \Cref{sec:collinear_jbd} for more details.
\end{remark}

\begin{remark}[From assumptions on factors to assumptions on the tensor]
The proof of \Cref{thm:pt2d_uniqueness_nonsymmetric} shows that \Cref{as:AB_full_rank,as:nonzero_F,as:krank_F,as:full_rank} hold if and only if
$\rank{\matr{T}^{(1)}} = R$ and $\rank{\matr{T}^{(2)}} = S$ and $\rank{\matr{\Phi}(\tens{T}_c)} = \binom{R+1}{2} \binom{S+1}{2}$, and $\tens{T}$ admits an $(R,S)$-PT2D.
Therefore \Cref{alg:paratuck_2x2_complete} returns the unique decomposition if the ``implicit'' conditions on $\tens{T}$ are satisfied.
\end{remark}

\paragraph{SVD- and EVD-based specialization of the algorithm}
Let us comment on how different steps of \Cref{alg:paratuck_2x2_complete} can be performed, by using only standard linear algebra operations, in the view of performing possibly the  approximate ParaTuck-2.
\begin{itemize}
\item Step 2 can by computed by the HOSVD: we take $\matr{U} \in \RR^{\edim{1} \times R}$, $\matr{V} \in \RR^{\edim{2} \times S}$  the matrices of the leading $R,S$ left singular vectors of $\matr{T}^{(1)}$ and $\matr{T}^{(2)}$, respectively. In this case, the core (compressed) tensor $\tens{T}_c \in \RR^{R \times S \times \edim{3}}$  is found as $\tens{T}_c = \tens{T} \con_1 \matr{U}^{\T}\con_2 \matr{V}^{\T}$
\item In particular, $R$ and $S$ can be estimated at the step $2$ in case they are not given \textit{a priori} (\emph{i.e.,} they can be estimated from the singular values of $\matr{T}^{(1)}$ and $\matr{T}^{(2)}$).
\item Step 3 can be performed by using the SVD of the matrix $\symlker{\matr{\Phi}(\tens{T}_c)}$. There are two options.
\begin{itemize}
\item The simplest is to compute the last $(RS)^2 - \binom{R+1}{2}\binom{S+1}{2}$ left singular vectors of $\symlker{\matr{\Phi}(\tens{T}_c)}$, symmetrize them with $\sigma_2(\cdot)$, and extract a basis of the symmetrized singular vectors.
\item Alternatively, we can form the matrix $\matr{\Phi}'(\tens{T}_c) \in \RR^{\binom{RS+1}{2} \times \edim{3}}$ by removing the repeating rows (which can only repeat twice) and multiplying the rows which had a repetition by a factor of $\sqrt{2}$, computing the last $M = \binom{R}{2} \binom{S}{2}$ left singular vectors  of  $\matr{\Phi}'(\tens{T}_c)$ (denoted as $\vect{p}'_1,\ldots,\vect{p}'_{M}$)
and expanding $\vect{p}'_k$ to $\vect{p}_k$ by repeating the corresponding elements and dividing them by  $\sqrt{2}$.

\end{itemize}
\item In steps $5$ and $6$: we can avoid extra symmetrization in $\pi$. In fact, let $\pi'$ be just the map that maps $\vecl{\vecl{\vect{a} \tensp \vect{b} \tensp  \vect{y} \tensp \vect{z}} \mapsto 
(\vect{y} \kron \vect{a}) \cdot \frac{1}{2}(\vect{z} \kron \vect{b})^{\T}}$.
Then  for any $\vect{p} \in \vecl{S^{2}(\RR^{RS})}$
\[
\range{\pi{(\vect{p})}} = \symglobal_2(\range{{\pi'(\vect{p})}}) \text{ and } \range{\pi{(\vect{p})}^{\T}} = \symglobal_2(\range{{\pi'(\vect{p})}^{\T}}),
\]
hence we can replace $\pi$ with $\pi'$ and  perform the symmetrization at the moment of computing $\symlker{\Mreshaped{A}}$ and $\symlker{\Mreshaped{B}}$.

\item For $\widehat{\tens{Q}}_A$, $\widehat{\tens{Q}}_B$ so that we can find the CPD (approximation) by a joint eigenvalue decomposition of several matrices.
Indeed, we have for any index $k$ 
\[
(\widehat{\tens{Q}}_A)_{k,:,:} = \widehat{\matr{A}}_c\Diag{({\widehat{\matr{A}}_c})_{k,:}}{\matr{W}_A}^{\T},
\]
and therefore
\[
 \matr{Y}_{k,\ell} =  (\widehat{\tens{Q}}_A)_{k,:,:}^{-1}(\widehat{\tens{Q}}_A)_{\ell,:,:} =
{\matr{W}_A}^{-\T} \matr{D}_{k,\ell} {\matr{W}_A}^{\T}
\]
with $\matr{D}_{k,\ell}  = \Diag{({\widehat{\matr{A}}_c})_{k,:}}^{-1}\Diag{({\widehat{\matr{A}}_c})_{\ell,:}}$.
Therefore, we can fix $k$ (say, to $k = 1$, which will choose a normalization $({\widehat{\matr{A}}_c})_{1,:} = 1$), and take $\ell = 2,\ldots,R$, to perform joint eigenvalue decomposition of $\matr{Y}_{1,2}, \ldots,\matr{Y}_{1,R}$ (for example, by the matrix pencil method).
This will give the remaining $(\widehat{\matr{A}}_c)_{\ell,:}$.
To avoid ill-conditioning, we can premultiply the first and second modes of the tensor  $\widehat{\tens{Q}}_A$ by a random rotation.

\end{itemize}


\section{Symmetric case (DEDICOM)}
As in the non-symmetric case,  we apply the lifting approach and the algorithm relies on the finding  the left nullspace of $\matr{\Phi}(\tens{T})$. However, there are important differences with respect to the nonsymmetric case.
\subsection{Properties of the structured matrix}
Note that the column span of  $\matr{\Phi}(\tens{T})$ has additional symmetry structure, in case when all slices  $\tens{T}_{:,:,k}$  are symmetric matrices.
Recall the definition of the subspace $\vecl{S^2(\vecl{S^2(\mathbb{R}^{R})})} \subset \RR^{R\times R\times R\times R}$, which is
\[
\vecl{S^2(\vecl{S^2(\mathbb{R}^{R})})} =  \Span{ \vecl{\matr{X}} \kron \vecl{\matr{X}}, \matr{X} \in \matr{S}^{2} (\mathbb{R}^{R})} .
\]
In fact, $\vecl{S^2(\vecl{S^2(\mathbb{R}^{R})})}$ is the space spanned by  the vectorizations of tensors $\tens{Y} \in \mathbb{R}^{R \times R\times R\times R}$ satisfying the following sets of symmetries
\[
\tenselem{Y}_{ijk\ell} = \tenselem{Y}_{jik\ell} = \tenselem{Y}_{ij\ell k} = \tenselem{Y}_{k\ell ij}.
\]
This implies the following remark.
\begin{remark}\label{lem:kernel_symmetric_symmetries}
For $\tens{T} = \ptdsym{\matr{A}}{\matr{F}}{\matr{G}} \in \RR^{R \times R \times \edim{3}}$, the column space of  $\matr{\Phi}(\tens{T})$ is a subspace of the space of partially symmetric tensors defined above.
\[
\range{\{ \matr{\Phi}(\tens{T})\}} \subseteq \vecl{S^2(\vecl{S^2(\mathbb{R}^{R})})} 
\]
\end{remark}
We will introduce the following useful definition. 
\begin{definition}
Let $\sigma_{(2,2)}: \RR^{R^4} \to \RR^{R^4} $ denote the orthogonal projection on $\vecl{S^2(\vecl{S^2(\RR^S)})}$, i.e. it is the map that acts on rank-one tensors as:
\[
\sigma_{(2,2)}:  \vect{z} \kron \vect{y} \kron \vect{b} \kron \vect{a} \mapsto 
\sigma(\sigma(\vect{z} \kron \vect{y}) \kron \sigma(\vect{b}\kron\vect{a})), \quad \mbox{for all } \vect{z}, \vect{y}, \vect{b} , \vect{a} \in \RR^{R}.
\]
Then for any matrix $\matr{Z} \in \RR^{R^4\times K}$, we define the symmetrized left kernel as $\symlkerd{\matr{Z}}$
\end{definition}
Note that due \Cref{lem:kernel_symmetric_symmetries}, we have (for any $\tens{T} = \ptdsym{\matr{A}}{\matr{F}}{\matr{G}} $)
\[
\symlkerd{ \matr{\Phi}(\tens{T})} = \lker{ \matr{\Phi}(\tens{T})} \cap \vecl{S^2(\vecl{S^2(\RR^R)})}
\]
Armed with the  definition of the symmetrized kernel, we can formulate the following analogue of \Cref{lem:Phi_rank_bound_nonsymmetric}.
\begin{lemma}\label{lem:Phi_rank_bound_symmetric}
\begin{enumerate}
\item 
 For $\tens{T} = \ptdsym{\matr{A}}{\matr{F}}{\matr{G}} \in \RR^{R\times S\times \edim{3}}$ with $\matr{A} \in \RR^{R\times R}$ and $\matr{F}$ symmetric we have
\begin{equation}\label{eq:Phi_rank_bound_symmetric}
\rank{\matr{\Phi}(\tens{T})} \le \binom{R+3}{4},
\end{equation}
\item If all \Cref{as:AB_full_rank,as:nonzero_F,as:full_rank_symmetric} are satisfied, then the equality in \eqref{eq:Phi_rank_bound_symmetric} is achieved.
In this case, $\symlkerd{\matr{\Phi}(\tens{T})}$ spanned by $N =  \frac{R^2(R-1)(R+1)}{12}$ vectors
\[
\symlkerd{\matr{\Phi}(\tens{T})} = \Span{\vect{p}_1,\ldots,\vect{p}_N}.
\]
\item If either  of \Cref{as:nonzero_F,as:full_rank_symmetric} is not satisfied, we have a strict inequality in \eqref{eq:Phi_rank_bound_symmetric}.
\end{enumerate}
\end{lemma}
\begin{proof}
The proof is given in \Cref{sec:symmetric_proofs}.
\end{proof}

\subsection{Main result for the symmetric case}
Recall that  for any matrix $\matr{Z} \in \RR^{R^3\times N}$, the symmetrization of the
left kernel is defined as $\symlkert{\matr{Z}}$.
\begin{theorem}\label{thm:pt2d_symmetric}
Let  $\tens{T} = \ptdsym{\matr{A}}{\matr{F}}{\matr{G}} \in \RR^{R\times S\times \edim{3}}$, whose factors satisfy  \Cref{as:AB_full_rank,as:nonzero_F,as:full_rank_symmetric} and let
$\matr{P}_k = \matricize{S^3}{S}{\sigma_{4}(\vect{p}_n)} $ where $\{\vect{p}_n\}_{n=1}^{N}$ be any basis of $\symlkerd{\matr{\Phi}(\tens{T})}$ (as in \Cref{lem:Phi_rank_bound_symmetric}).
Consider the matrix
\begin{equation}\label{eq:Psym}
\matr{P}_{sym} = \begin{bmatrix}\matr{P}_1 &\cdots & \matr{P}_N\end{bmatrix}
\end{equation}
Then,  $\rank{\matr{P}_{sym}} \le \binom{R+2}{3} - R$ and  the equality holds if and only if \Cref{as:krank_F} holds. In the latter case, we have
\[
\symlkert{\matr{P}_{sym}} = \Span{ \vect{a}_1 \kron \vect{a}_1 \kron \vect{a}_1, \ldots,\vect{a}_R \kron \vect{a}_R  \kron \vect{a}_R} = \range{\matr{A}\kr\matr{A}\kr\matr{A}}.
\]
\end{theorem}
\begin{proof}
The proof will be split into several lemmas and presented later in \Cref{sec:symmetric_proofs}.
\end{proof}
Then \Cref{thm:pt2d_symmetric} implies the following improved result on identifiability of PT2D.

\begin{proposition}\label{thm:pt2d_uniqueness_symmetric}
Let $\tens{T}  = \ptdsym{\matr{A}}{\matr{F}}{\matr{G}} \in \RR^{\edim{1}\times \edim{1} \times \edim{3}}$, $\matr{F} \in \RR^{R \times R}$ symmetric, whose factors satisfy  \Cref{as:AB_full_rank,as:nonzero_F,as:full_rank_symmetric,as:krank_F}.
Then the $R$-DEDICOM of $\tens{T}$ is essentially unique.
In addition, $R$ is the minimal possible rank for a DEDICOM of $\tens{T}$.
\end{proposition} 
\begin{proof}
Again, due to \Cref{lem:pt2d_tucker2} and  \Cref{as:AB_full_rank,as:nonzero_F,as:full_rank_symmetric}, we have that $\rank{\matr{T}^{(1)}} = R$ which implies that there is no $R'$-DEDICOM for $R' < R$.
Also,  we can assume, without loss of generality that $\edim{1} = R$.
If $\edim{1} > R$, then we can first perform Tucker compression on $\tens{T}$  (see \Cref{lem:pt2d_tucker2}) and compute the PT2D on the compressed tensor, as int the proof of \Cref{thm:pt2d_uniqueness_symmetric}.

Under \Cref{as:AB_full_rank,as:nonzero_F,as:full_rank_symmetric}, by  \Cref{lem:Phi_rank_bound_nonsymmetric} we have $\rank{\matr{\Phi}(\tens{T})} = \binom{R+3}{4}$.
If \Cref{as:krank_F} is also satisfied, we have  $\dim(\symlkert{\matr{P}_{sym}})  = R$ by \Cref{thm:pt2d_symmetric}.

Let $\matr{q}_{1}, \ldots, \matr{q}_{R} \in \RR^{R^3}$ be some basis of $\symlkert{\matr{P}_{sym}}$. Then, by \Cref{thm:pt2d_nonsymmetric} there is a nonsingular (change of basis) matrix $\matr{W} \in \RR^{R\times R}$
\[
\begin{bmatrix}\matr{q}_{1} & \ldots & \matr{q}_{R} \end{bmatrix} = (\matr{A} \kr \matr{A} \kr \matr{A}) \matr{W}^{\T},
\]
Which implies that the tensor $\tens{Q} \in \RR^{R\times R \times R \times R}$ built as ${(\tens{Q})_{:,:,:,k}} =\tensorize{R}{R}{R}{\vect{q}_{k}}$ has a CPD
\[
\tens{Q} = \cpdfour{\matr{A}}{\matr{A}}{\matr{A}}{\matr{W}}.
\]
Note that, by our assumptions, $\krank{\matr{A}} = \krank{\matr{W}} = R$. Thus  the Kruskal's uniqueness conditions are satisfied and the matrix $\matr{A}$ is unique up to permutation and scaling of columns.
The factors $\matr{F},\matr{G}$ can be recovered from $\tens{C}$ up to scaling indeterminacies, thanks to \Cref{prop:core_factorization}.

Finally, due to ``if and only if'' statements in  \Cref{lem:Phi_rank_bound_symmetric} and \Cref{thm:pt2d_symmetric},  for any alternative  $(R)$-DEDICOM,  $\tens{T}  = \ptdsym{\widetilde{\matr{A}}}{\widetilde{\matr{F}}}{\widetilde{\matr{G}}}$,
all \Cref{as:AB_full_rank,as:nonzero_F,as:full_rank_symmetric,as:krank_F} must be satisfied as well.
\end{proof}

\subsection{Algorithm for the symmetric case}
The proofs of  \Cref{thm:pt2d_uniqueness_symmetric,thm:pt2d_symmetric} are also constructive.
We summarize in \Cref{alg:paratuck_2x2_symmetric} the algorithms for  computing the (unique) R-DEDICOM of a given tensor, under   \Cref{as:AB_full_rank,as:nonzero_F,as:krank_F} and  \Cref{as:full_rank_symmetric}, respectively. 
\begin{algorithm}[hbt!]
\caption{Algebraic decomposition algorithm for DEDICOM}\label{alg:paratuck_2x2_symmetric}
\begin{algorithmic}[1]
\STATE{{\bf Input:}$\tens{T} \in \RR^{\edim{1},\edim{1},\edim{3}}$ rank $(R)$}
\STATE{Compute the symmetric Tucker-2 decomposition (approximation) $\tens{T} = \tens{T}_c \con_1 \matr{U}\con_2 \matr{U}$  with $\matr{U} \in \RR^{\edim{1} \times R}$}
\STATE{Find the $N$ vectors $\vect{p}_1,\ldots,\vect{p}_{N}$ in (approximate) symmetric left kernel  $\symlkerd{\matr{\Phi}(\tens{T}_c)}$}
\STATE{Find the symmetrizations $\matr{P}_k =   \matricize{R^3}{R}{\sigma_{4}(\vect{p}_n)}$, build $\matr{P}_{sym}$ from $\matr{P}_k$.}
\STATE{Compute a  basis $\{\matr{q}_{1},\ldots,\matr{q}_{R}\}$  of the (approximate) symmetric left kernel $\symlkert{\matr{P}_{sym}}$.}
\STATE{Stack these matrices in tensor $\widehat{\tens{Q}}_{sym} \in \RR^{R \times R\times R\times R}$ and compute its CPD (approximation)

$\widehat{\tens{Q}}_{sym} = \cpdfour{\widehat{\matr{A}}_c}{\widehat{\matr{A}}_c}{\widehat{\matr{A}}_c}{{\matr{W}_A}}$.}
\STATE{Set  $\widehat{\matr{A}} = \matr{U} \widehat{\matr{A}}_c$ and find the core tensor $\widehat{\tens{\Core}} = \tens{T}_c \contr{1} \widehat{\matr{A}}^{\dagger}_c \contr{2} \widehat{\matr{A}}^{\dagger}_c$. }
\STATE{Determine the factors $\widehat{\matr{\CoreFront}}$, $\widehat{\matr{\CoreVert}}$ using the methods from \Cref{sec:factorizing_core}  (namely, \Cref{prop:core_factorization}).}
\RETURN{$\widehat{\matr{A}},\widehat{\matr{F}},\widehat{\matr{G}}$, such that $\tens{T} = \ptdsym{\widehat{\matr{A}}}{\widehat{\matr{F}}}{\widehat{\matr{G}}}$}
\end{algorithmic}
\end{algorithm}

As in the nonsymmetric algorithm, we make remarks on each of the individual steps of the algorithm. 

\paragraph{SVD- and EVD- based specializations of the algorithm}

\begin{itemize}
\item In step 2, the HOSVD again is used: we can compute $\matr{U}^{I \times R}$ from singular vectors and also estimate $R$ if needed.
\item The symmetrized kernel $\symlkerd{\matr{\Phi}(\tens{T}_c)}$ can be computed as in the nonsymmetric case.

\item The symmetrization $\sigma_4$ can be skipped and combined with computing $\symlkert{\matr{P}_{sym}}$.
 
 \item The CPD of $\widehat{\tens{Q}}_{sym}$ can be computed by joint eigenvalue decomposition of $R\times R$ matrices
 \[
 \matr{Y}_{k,\ell} = 
  (\matricize{R^2}{R}{(\widehat{\tens{Q}}_A)_{k,:,:,:}})^{\dagger}\matricize{R^2}{R}{(\widehat{\tens{Q}}_A)_{\ell,:,:,:}}.
 \] 
\end{itemize}



\section{Properties of the core tensor}
In this section, we review advanced properties of the core tensor that serve as the base for the proofs.
\subsection{Implicit equations for the core tensor and their determinantal representation}\label{sec:core_equations}
We begin by recalling some equations that must be satisfied by a core tensor
\begin{lemma}
Let $\tens{C} \in \RR^{R \times S \times \edim{3}}$ be a tensor of the form \eqref{eq:paratuck_core}, i.e., $\tens{C} = \triprod{\matr{\CoreFront}}{\matr{\CoreVert}}{\matr{\CoreHoriz}}$.
Then for all $1 \le i,r \le R$, $1 \le j,s \le S$, $1 \le k,t \le \edim{3}$ it holds that
\begin{equation}\label{eq:equations_core}
\tenselem{\Core}_{ijk}\tenselem{\Core}_{rsk}\tenselem{\Core}_{ist}\tenselem{\Core}_{rjt} - 
\tenselem{\Core}_{isk}\tenselem{\Core}_{rjk}\tenselem{\Core}_{ijt}\tenselem{\Core}_{rst} = 0.
\end{equation}
\end{lemma}
\begin{proof}
Follows from straightforward substitution, as
\[
\tenselem{\Core}_{ijk}\tenselem{\Core}_{rsk}\tenselem{\Core}_{ist}\tenselem{\Core}_{rjt} = 
{\CoreFront}_{ij}{\CoreFront}_{rs} {\CoreFront}_{is}{\CoreFront}_{rj}
{\CoreVert}_{ik}{\CoreVert}_{rk} {\CoreVert}_{it}{\CoreVert}_{rt} 
{\CoreHoriz}_{jk}{\CoreHoriz}_{sk} {\CoreHoriz}_{jt}{\CoreHoriz}_{st}
 = \tenselem{\Core}_{isk}\tenselem{\Core}_{rjk}\tenselem{\Core}_{ijt}\tenselem{\Core}_{rst}
\]
\end{proof}
The equation \eqref{eq:equations_core} can be viewed as vanishing of a ``generalized determinant'' for  $2\times 2 \times 2$ tensors (see  \eqref{fig:core_equations}). 
Note that there are other polynomial equations for tensors satisfying \eqref{eq:equations_core}, and describing them is a difficult task \cite[end of \S6]{GarcSS05:AG}.
However, for our purposes, equations \eqref{eq:equations_core} will be largely sufficient.


\begin{figure}[ht!]
\centering
$\begin{array}{l}
\color{blue}{\tenselem{\Core}_{ijk}\tenselem{\Core}_{rsk}\tenselem{\Core}_{ist}\tenselem{\Core}_{rjt}}\color{black} - \\
\color{red}{\tenselem{\Core}_{isk}\tenselem{\Core}_{rjk}\tenselem{\Core}_{ijt}\tenselem{\Core}_{rst}}\color{black} = 0
\end{array}
\quad\quad\quad
\raisebox{-0.5\height}{\begin{tikzpicture}[square/.style={regular polygon,regular polygon sides=4}]
\begin{scope}[xscale= 0.75,yscale=0.75]
\tikzstyle{p1} = [color = red,draw,square,minimum size=1.5mm,fill,inner sep=0mm]
\tikzstyle{p0} = [color = blue,draw,circle,minimum size=1.5mm,fill,inner sep=0mm]

\draw(-0.6,1) node (i) {\small $i$};
\draw(-0.6,0) node (p) {\small $r$};
\draw(-0.5, -0.5) node (j) {\small $j$};
\draw(0.5, -0.5) node (q) {\small $s$};
\draw(0, 1.6) node (k) {\small $k$};
\draw(0.5, 2.1) node (r) {\small $t$};
\draw(0,0) node[p1]          (v0001) {};
\draw(1,0) node[p0]          (v0011) {};
\draw(0,1) node[p0]          (v0000) {};
\draw(1,1) node[p1]          (v0010) {};
\draw(0.5,0.5) node[p0]      (v0101) {};
\draw(1.5,0.5) node[p1]      (v0111) {};
\draw(0.5,1.5) node[p1]      (v0100) {};
\draw(1.5,1.5) node[p0]      (v0110) {};

\draw[->] (i) -- (v0000);
\draw[->] (p) -- (v0001);
\draw[->] (j) -- (v0001);
\draw[->] (q) -- (v0011);
\draw[->] (k) -- (v0001);
\draw[->] (r) -- (v0100);

\draw (v0000) -- (v0010) -- (v0011) -- (v0001) -- (v0000);
\draw (v0100) -- (v0110) -- (v0111) -- (v0101) -- (v0100);
\draw (v0000) -- (v0100);
\draw (v0010) -- (v0110) ;
\draw (v0001) -- (v0101);
\draw (v0011) -- (v0111) ;
\end{scope}
\end{tikzpicture}}$
\caption{Equations for core tensor via in terms of the elements of the $2\times 2\times 2$ subtensor $\tens{C}_{(i,r),(j,s),(k,t)}$.}
\label{fig:core_equations}
\end{figure}
Note that the equation \eqref{eq:equations_core} can be also rewritten as
\begin{equation}\label{eq:equations_core_2x2}
\det 
\begin{bmatrix}
\tenselem{\Core}_{ijk}\tenselem{\Core}_{rsk} & \tenselem{\Core}_{ijt}\tenselem{\Core}_{rst} \\
\tenselem{\Core}_{isk}\tenselem{\Core}_{rjk} & \tenselem{\Core}_{ist}\tenselem{\Core}_{rjt}
\end{bmatrix} = 0
\end{equation}
Then equations \eqref{eq:equations_core_2x2} for all $1 \le k < t \le \edim{3}$ can be grouped to 
obtain the following result.
\begin{lemma}\label{lem:determinantal_equations_C}
Let $\tens{C} \in \RR^{R \times S \times \edim{3}}$ be a tensor such that $\tens{C} = \triprod{\matr{\CoreFront}}{\matr{\CoreVert}}{\matr{\CoreHoriz}}$.
\begin{enumerate}
\item Then for all $1 \le i, r \le R$, $1 \le j, s \le S$, the matrix
\[
\matr{\Psi}^{(ijrs)} := 
\begin{bmatrix}
\tenselem{\Core}_{ij1}\tenselem{\Core}_{rs1}  & \cdots & \tenselem{\Core}_{ij\edim{3}}\tenselem{\Core}_{rs\edim{3}} \\
\tenselem{\Core}_{is1}\tenselem{\Core}_{rj1} & \cdots & \tenselem{\Core}_{ij\edim{3}}\tenselem{\Core}_{rs\edim{3}}
\end{bmatrix} \in \RR^{2 \times \edim{3}}
\]
must be at most rank $1$.
\item In particular, the vector 
\begin{equation}\label{eq:theta_expr}
\begin{bmatrix}
\theta^{(ijrs)}_1 &
\theta^{(ijrs)}_2
\end{bmatrix} :=
\begin{bmatrix}
-{\CoreFront}_{is}{\CoreFront}_{rj} &
{\CoreFront}_{ij}{\CoreFront}_{rs} 
\end{bmatrix},
\end{equation}
lies in the left kernel of $\matr{\Psi}^{(ijrs)}$, i.e.,
\begin{equation}\label{eq:left_kernel}
\begin{bmatrix}
\theta^{(ijrs)}_1 &
\theta^{(ijrs)}_2
\end{bmatrix} \matr{\Psi}^{(ijrs)} = 0.
\end{equation}
\end{enumerate}
\end{lemma}
\begin{proof}
\begin{enumerate}
\item
Let us fix $1 \le i, r \le R$, $1 \le j, s \le S$. Then, thanks to \eqref{eq:equations_core_2x2}, vanishing of all $2\times2$ minors of $\matr{\Psi}^{(ijrs)}$ is equivalent to the inequality $\rank{\matr{\Psi}^{(ijrs)}} \le 1$.

\item The straightforward computation leads to 
\[
(\matr{\Psi}^{(ijrs)})_{:,k} = (\CoreVert_{ik}\CoreVert_{rk} \CoreHoriz_{jk}\CoreHoriz_{sk})
\vect{z}^{(ijrs)}, \quad \text{with }
\vect{z}^{(ijrs)} = 
\begin{bmatrix}
{\CoreFront}_{ij}{\CoreFront}_{rs} \\
{\CoreFront}_{is}{\CoreFront}_{rj}
\end{bmatrix},
\]
or, equivalently
\begin{equation}\label{eq:Psi_ijrs_rank_one}
\matr{\Psi}^{(ijrs)} = \vect{z}^{(ijrs)} (\matr{\CoreVert}_{i,:} \kr \matr{\CoreVert}_{r,:} \kr \matr{\CoreHoriz}_{j:} \kr\matr{\CoreHoriz}_{s,:})
\end{equation}
Note that the vector $\begin{bmatrix}
\theta^{(ijrs)}_1 &
\theta^{(ijrs)}_2
\end{bmatrix}^{\T}$ is obviously orthogonal to $\vect{z}^{(ijrs)}$ and the proof is complete.
\end{enumerate}
\end{proof}

\subsection{Factorizing the core tensor: proofs}\label{sec:factorizing_core_proofs}


\begin{proof}[Proof \Cref{lem:nonzero_slice}]
Indeed, the condition of lemma implies ${\CoreFront}_{ij} \neq 0$ and hence we can take  $\widetilde{\matr{\CoreFront}} = \tens{\Core}_{:,:,r} $ and $\widetilde{\matr{g}}_{r} = \matr{1}$ and  $\widetilde{\matr{h}}_{r} = \matr{1}$.
Now, for fixed $k$ let $\matr{Z}_k \in \RR^{R\times S}$ be defined by elementwise division $(\matr{Z}_k)_{i,j} =  \tenselem{\Core}_{i,j,r} / \widetilde{{\CoreFront}}_{ij}$.
Then  we can find the k-th columns of $\widetilde{\matr{\CoreVert}}$ and $\widetilde{\matr{\CoreHoriz}}$ (i.e., $\widetilde{\matr{g}}_{k}$ and  $\widetilde{\matr{h}}_{k}$) from the rank-one factorization of $\matr{Z}_k$:
\[
\widetilde{\matr{g}}_{k} (\widetilde{\matr{h}}_{k})^{\T} = \matr{Z}_k.
\]
It is easy to see that for such matrices we have $\tens{\Core} = \triprod{\widetilde{\matr{\CoreFront}}}{\widetilde{\matr{\CoreVert}}}{\widetilde{\matr{\CoreHoriz}}}$, thus it gives a valid factorization of $\tens{\Core}$.
\end{proof}

\begin{proof}[Proof of  \Cref{prop:core_factorization}]
By the scaling ambiguities, we can take always take $\matr{F}$ to have the first row and column consisting of all ones.
Then we have that 
\begin{equation}\label{eq:F_pattern}
\matr{F} = 
\begin{bmatrix}
1 & 1 & \cdots & 1 \\
1 & \ast & \cdots & \ast \\
\vdots &\vdots && \vdots\\
1 & \ast & \cdots & \ast 
\end{bmatrix}.
\end{equation}
We will show that we can complete the $\ast$ elements uniquely. 

If \Cref{as:full_rank} or  \Cref{as:full_rank_symmetric} is satisfied, neither of rows of $\matr{G}\kr\matr{G}\kr\matr{H}\kr \matr{H}$ is identically zero.
That is, for arbitrary $1  <  r \le R$, $1 <s \le S$, 
\[
(\matr{\CoreVert}_{1,:} \kr \matr{\CoreVert}_{r,:} \kr \matr{\CoreHoriz}_{1,:} \kr\matr{\CoreHoriz}_{s,:})\neq \vect{0}
\]
and hence $\rank{\matr{\Psi}^{(11rs)}} =1$ in view of \eqref{eq:Psi_ijrs_rank_one}.
Therefore, we can determine $\theta^{(11rs)}_1,\theta^{(11rs)}_2$ in \eqref{eq:left_kernel} up to a constant from $\lker{\matr{\Psi}^{(11rs)}}$, and thus
\[
 F_{rs} = \frac{F_{rs} F_{11}}{F_{1s} F_{r1}}  =   -\frac{\theta^{(11rs)}_2}{\theta^{(11rs)}_1},
\]
can be uniquely recovered from $\Psi^{(11rs)}$.
\end{proof}

\begin{remark}
Note that we can relax  \Cref{as:full_rank,as:full_rank_symmetric} to the assumption of  existence of $i,j$ such that $(\matr{\CoreVert}_{1,:} \kr \matr{\CoreVert}_{r,:} \kr \matr{\CoreHoriz}_{1,:} \kr\matr{\CoreHoriz}_{s,:})\neq \vect{0}$
for all $r\neq i$ and $s \neq j$.
\end{remark}

\subsection{Multilinear ranks of the core tensor}\label{sec:proofs_ml_rank}

\begin{proof}[Proof of \Cref{prop:core_multilinear_rank}]
We will just prove the first statement as the second is analogous.
The first unfolding of $\tens{C}$ has the following block structure
\[
\matr{C}^{(1)} = 
\begin{bmatrix}
\DiagA{1} \matr{\CoreFront} \DiagB{1} & \DiagA{2} \matr{\CoreFront} \DiagB{2} & \cdots&
\DiagA{\edim{3}} \matr{\CoreFront} \DiagB{\edim{3}}
\end{bmatrix}.
\]
Let us select every $j$-th column of each block of $\matr{C}^{(1)}$ and stack them into a submatrix  $\matr{Y}\in \RR^{R \times K}$.
This submatrix has the form
\begin{align*}
\matr{Y}&  = \begin{bmatrix}
\DiagA{1} \matr{\CoreFront}_{:,j} H_{j,1} & \DiagA{2} \matr{\CoreFront}_{:,j} H_{j,2} & \cdots&
\DiagA{\edim{3}}\matr{\CoreFront}_{:,j} H_{j,\edim{3}} 
\end{bmatrix} \\
& =
\Diag{\matr{\CoreFront}_{:,j} }
\begin{bmatrix}
\vect{g}_1 H_{j,1} & \vect{g}_2 H_{j,2} & \cdots & \vect{g}_{\edim{3}} H_{j,\edim{3}}
\end{bmatrix} = \Diag{\matr{\CoreFront}_{:,j} } (\matr{G} \kr \matr{H}_{j,:})
\end{align*}
where $\Diag{\matr{\CoreFront}_{:,j} }$ is nonsingular, hence $\rank{\matr{Y}} = R$ by assumption.
This implies that $\rank{\matr{C}^{(1)}} = R$.
\end{proof}

\begin{proof}[{Proof of \Cref{cor:core_ml_rank_assumptions}}]
\begin{enumerate}
\item
Let \Cref{as:full_rank} be satisfied. Then  the matrix  $\matr{G} \kr \matr{G}\kr\matr{H} \kr \matr{H}$ has $\binom{R+1}{2}\binom{S+1}{2}$ distinct rows, and, in particular, any collection of the distinct rows must be linearly independent. This implies that for any $i,j$
\begin{equation}\label{eq:ranks_subsets_GGHH}
\rank{\matr{G} \kr \matr{G}_{i,:}\kr\matr{H}_{j,:} \kr \matr{H}_{j,:}} = R \mbox{ and }
\rank{\matr{G}_{i,:} \kr \matr{G}_{i,:}\kr\matr{H}_{j,:} \kr \matr{H}} = S.
\end{equation}
Note that the matrices in \eqref{eq:ranks_subsets_GGHH} are just diagonally scaled versions of 
$\matr{G}  \kr \matr{H}_{j,:}$ and $\matr{G}_{i,:} \kr \matr{H}$, hence   $\rank{\matr{G}  \kr \matr{H}_{j,:}} = R$, $\rank{\matr{G}_{i,:} \kr \matr{H}} = S$ and thus the conditions of \Cref{prop:core_multilinear_rank} are satisfied.
\item If  \Cref{as:full_rank_symmetric} is satisfied, then the proof is very similar.
In this case, the matrix   $\matr{G} \kr \matr{G}\kr\matr{G} \kr \matr{G}$ has $\binom{R+3}{4}$ distinct rows, and, as in the previous case, we have that for any $i$
\begin{equation}\label{eq:ranks_subsets_GGGG}
\rank{\matr{G} \kr \matr{G}_{i,:}\kr\matr{G}_{i,:} \kr \matr{G}_{i,:}} = R,
\end{equation}
hence $\rank{\matr{G} \kr \matr{G}_{i,:}} = R$, and  the conditions \Cref{prop:core_multilinear_rank} are also satisfied (since $\matr{G} = \matr{H}$).
\end{enumerate}
\end{proof}

\section{Proofs for the $2\times 2$ case}
\subsection{Building intuition: the $2\times 2$ nonsymmetric case}\label{sec:2x2_nonsymmetric}
We  provide a sketch proof of \Cref{thm:pt2d_nonsymmetric} for a particular case $R=S=2$.
Assume that \cref{as:AB_full_rank,as:nonzero_F,as:full_rank,as:krank_F}.
Note that for $R=2$, $S=2$ there is essentially one nontrivial matrix $\matr{\Psi}^{(ijrs)}$ that can be obtained, by choosing $(i,r)= (j,s)  = (1,2)$.
It is easy to see that $\matr{\Psi}^{(1122)}$ is, in fact, a submatrix $\Phi(\tens{C})$
then the condition \eqref{eq:left_kernel} on $\matr{\Psi}^{(1122)}$ can be rewritten  as
\[
(\theta_1 \vecl{\vect{e}_1 \tensp \vect{e}_1 \tensp \vect{e}_2 \tensp \vect{e}_2}+ 
\theta_2 \vecl{\vect{e}_1 \tensp \vect{e}_2 \tensp \vect{e}_2 \tensp \vect{e}_1})^{\T} \Phi(\tens{C}) = 0,
\]
where $\theta_1 = \theta^{(1122)}_1, \theta_2 = \theta^{(1122)}_2$ are the coefficients defined in \Cref{lem:determinantal_equations_C}. 
Now let $\widetilde{\vect{a}}_{k}$ and  $\widetilde{\vect{b}}_{k}$ be the rows of the inverses of $\matr{A}$ and   $\matr{B}$ respectively\footnote{These are the dual vectors, $\langle\widetilde{\vect{a}}_{i},{\vect{a}}_{j}\rangle = \langle\widetilde{\vect{b}}_{i},{\vect{b}}_{j}\rangle  = \delta_{ij}$}
\[
\matr{A}^{-\T} = \begin{bmatrix} \widetilde{\vect{a}}_1 &  \widetilde{\vect{a}}_2  \end{bmatrix},
\quad
\matr{B}^{-\T} = \begin{bmatrix} \widetilde{\vect{b}}_1 &  \widetilde{\vect{b}}_2  \end{bmatrix}.
\]
Then by duality (i.e., using the fact that $\Phi(\tens{T}) = (\matr{B} \kron\matr{A}\kron \matr{B}\kron \matr{A}) \Phi(\tens{C})$), we have that the vector
\begin{equation}\label{eq:v_1122}
\vect{v} := (\theta_1 \vecl{\widetilde{\vect{a}}_1 \tensp \widetilde{\vect{b}}_1 \tensp \widetilde{\vect{a}}_2 \tensp \widetilde{\vect{b}}_2}+ 
\theta_2 \vecl{\widetilde{\vect{a}}_1 \tensp \widetilde{\vect{b}}_2 \tensp \widetilde{\vect{a}}_2 \tensp \widetilde{\vect{b}}_1})
\end{equation}
lies in the left kernel of $\matr{\Phi}(\tens{T})$ (i.e. $\vect{v} ^{\T}\matr{\Phi}(\tens{T}) = 0$).
Then  the symmetrized vector is nonzero 
\begin{align*}
\symglobal_2(\vect{v} ) ={\theta_1}\sigma_2(\vecl{\widetilde{\vect{a}}_1 \tensp \widetilde{\vect{a}}_1 \tensp \widetilde{\vect{a}}_2 \tensp \widetilde{\vect{a}}_2})+{\theta_2} \sigma(\vecl{\widetilde{\vect{a}}_1 \tensp \widetilde{\vect{a}}_2 \tensp \widetilde{\vect{a}}_2 \tensp \widetilde{\vect{a}}_1}) \neq 0
\end{align*}
 and it lies in $\symlker{\matr{\Phi}(\tens{T})}$. By \Cref{lem:Phi_rank_bound_nonsymmetric}, this is, in fact, the only generator of the symmetrized kernel:
\[
\symlker{\matr{\Phi}(\tens{T})} = \Span{\symglobal_{2}(\vect{v})}.
\]
Thus any  vector $\vect{p}$ generating $\symlker{\matr{\Phi}(\tens{T})}$ has the form $c\symglobal_2(\vect{v})$ for come nonzero constant.
Now, applying the symmetrization and permutation map \eqref{eq:symperm} leads to:
\[
\matr{P}  = \symperm(c \sigma_2(\vect{v}))  =
c(\det{\matr{F}}) \symglobal_2(\widetilde{\vect{a}}_1 \kron \widetilde{\vect{a}}_2) 
(\symglobal_2(\widetilde{\vect{b}}_1 \kron \widetilde{\vect{b}}_2))^{\T},
\]
where we have used the fact that $ (\theta_1+\theta_2) = \det{\matr{F}}$ by \Cref{lem:determinantal_equations_C} (which is nonzero by \Cref{as:krank_F}).
This means that $\matr{P}$ is a rank-$1$ matrix whose column and row spaces are 
\[
\range{\matr{P}} = \Span{ \symglobal(\widetilde{\vect{a}}_1 \kron \widetilde{\vect{a}}_2)}, \quad
\range{\matr{P}^{\T}} = \Span{ \symglobal(\widetilde{\vect{b}}_1 \kron \widetilde{\vect{b}}_2)} 
\]
Next, observe that
\[
\langle \vect{a}\kron \vect{a},\symglobal(\widetilde{\vect{a}}_1 \kron \widetilde{\vect{a}}_2) \rangle= \langle \vect{a},\widetilde{\vect{a}}_1\rangle\langle \vect{a},\widetilde{\vect{a}}_2\rangle,
\]
which shows that both $\vect{a}_1\kron \vect{a}_1$ and  $\vect{a}_2\kron \vect{a}_2$ are orthogonal to $\symglobal(\widetilde{\vect{a}}_1 \kron \widetilde{\vect{a}}_2)$.
These two vectors are linearly independent, and since $\dim(\symlker{\matr{Q}}) = 2$, we get that
\[
\symlker{\matr{P}} = \Span{\vect{a}_1\kron \vect{a}_1,\vect{a}_2\kron \vect{a}_2}.
\]
By  similar manipulations, we can show that
\[
\symlker{\matr{P}^{\T}} = \Span{\vect{b}_1\kron \vect{b}_1,\vect{b}_2\kron \vect{b}_2},
\]
which completes the proof of \Cref{thm:pt2d_nonsymmetric} for a particular case $R=S=2$.

Finally, we comment on how to simply obtain $\vect{a}_1,\vect{a}_2$ from $c\symglobal(\widetilde{\vect{a}}_1 \kron \widetilde{\vect{a}}_2)$.
Assume for simplicity, that 
\[
\vect{a}_1 = \begin{bmatrix}1 & x\end{bmatrix}^{\T},\vect{a}_2 =  \begin{bmatrix}1 & y\end{bmatrix}^{\T}.
\]
then $\widetilde{\vect{a}}_1 = c_1  \begin{bmatrix}-x & 1\end{bmatrix}^{\T}$,  $\widetilde{\vect{a}}_2  =  c_2  \begin{bmatrix}-y & 1\end{bmatrix}^{\T}$  and thus
\[
c\symglobal(\widetilde{\vect{a}}_1 \kron \widetilde{\vect{a}}_2) = [z_0,z_1,z_1,z_2]^{\T} := c   \begin{bmatrix}xy &-(x+y)/2 & -(x+y)/2& 1\end{bmatrix}^{\T}.
\]
Now, given $z_0,z_1,z_2$, we can see that $x,y$ must be the roots of the univariate polynomial $z(t) = z_0 + 2z_1t+z_2t^2$, thus the vectors $\vect{a}_1,\vect{a}_2$ can be easily found from $\symglobal(\widetilde{\vect{a}}_1 \kron \widetilde{\vect{a}}_2)$.

\subsection{Building intuition: the $2\times 2$ symmetric case}
Now let us consider the $2\times 2$ symmetric case, i.e. $R =S=2$, $\matr{A} = \matr{B}$, $\matr{G} = \matr{H}$.
Note that $\vect{v}$ still lies in the left nullspace of  $\matr{\Phi}(\tens{T})$. Let us denote by $\vect{w}$ the projection this vector on $\vecl{S^2(\vecl{S^2(\RR^2)})}$ by double symmetrization, which has the form
\[
\vect{w} = {\theta_1}\sigma_2(\widetilde{\vect{a}}_1 \kron \widetilde{\vect{a}}_1 \kron \widetilde{\vect{a}}_2 \kron \widetilde{\vect{a}}_2)+{\theta_2} \sigma_2(\widetilde{\vect{a}}_1 \kron \widetilde{\vect{a}}_2) \kron \sigma_2(\widetilde{\vect{a}}_1 \kron \widetilde{\vect{a}}_2).
\]
Then $\vecl{S^2(\vecl{S^2(\RR^R)})} \cap \lker{\matr{\Phi}(\tens{T})}$ is spanned by a single vector $\vect{q} = c\vect{w}$.
Now take a projection of this vector on $S^4(\RR^2)$, which is 
\[
\sigma_4(\vect{w}) = \sigma_4(\vect{v}) = (\det{\matr{F}})
\sigma_4 (\widetilde{\vect{a}}_1 \kron \widetilde{\vect{a}}_1 \kron \widetilde{\vect{a}}_2 \kron \widetilde{\vect{a}}_2).
\]
Let $\matr{Q} \in \RR^{2^3 \times 2}$ be the reshaping of $\sigma_4 (\widetilde{\vect{a}}_1 \kron \widetilde{\vect{a}}_1 \kron \widetilde{\vect{a}}_2 \kron \widetilde{\vect{a}}_2)$, then the matix $\matr{Q}$ is of rank $2$ because it has the following low-rank decomposition
\[
\matr{Q} = \frac{1}{2} \left( \sigma_3 (\widetilde{\vect{a}}_1 \kron \widetilde{\vect{a}}_1 \kron \widetilde{\vect{a}}_2 ) \vect{a}_2^{\T} +  \frac{1}{2} \sigma_3 (\widetilde{\vect{a}}_1 \kron \widetilde{\vect{a}}_2 \kron \widetilde{\vect{a}}_2 ) \vect{a}_1^{\T}\right). 
\]
Note that, one the one hand
\[
\dim(\vecl{S^3(\RR^{3})}\cap\lker{\matr{Q}} = \dim(\vecl{S^3(\RR^{3})}) - \rank{\matr{Q}} = 2,
\]
and on the other hand, the two vectors  ${\vect{a}}_1 \kron {\vect{a}}_1 \kron {\vect{a}}_1$, ${\vect{a}}_2 \kron {\vect{a}}_2 \kron {\vect{a}}_2$ are in the left nullspace of $\matr{Q}$.
Therefore the symmetric left nullspace is spanned by those vectors.

\begin{remark}
The $2\times2$ symmetric case has also interpretation in terms of univariate polynomials. Let $\vect{a}_1, \vect{a}_2$ and $\widetilde{\vect{a}}_1, \widetilde{\vect{a}}_2$ be as in the previous subsection. Then the polynomial corresponding to the vector $\sigma_4 (\widetilde{\vect{a}}_1 \kron \widetilde{\vect{a}}_1 \kron \widetilde{\vect{a}}_2 \kron \widetilde{\vect{a}}_2 )$
\[
z(t) = \left\langle \sigma_4 (\widetilde{\vect{a}}_1 \kron \widetilde{\vect{a}}_1 \kron \widetilde{\vect{a}}_2 \kron \widetilde{\vect{a}}_2 ) , \begin{bmatrix} 1 \\ t\end{bmatrix}^{\kron 4} \right\rangle  = \left(\left\langle \widetilde{\vect{a}}_1, \begin{bmatrix} 1 \\ t\end{bmatrix}\right\rangle\right)^2\left( \left\langle \widetilde{\vect{a}}_2, \begin{bmatrix} 1 \\ t\end{bmatrix}\right\rangle\right)^2 =  c_1^2c_2^2(t-x)^2 (t-y)^2
\]
thus it has two double roots. Therefore $x,y$ can be found just by finding the  roots of the polynomial $\sqrt{z(t)}$.
\end{remark}

\section{Nonsymmetric case: proof of the main theorem}\label{sec:nonsymmetric_proofs}
\subsection{Basis of the left kernel}
Similarly to \Cref{sec:2x2_nonsymmetric}, we can get the basis of the symmetric kernel of the structured matrix $\matr{\Phi}(\tens{T})$ in the general case $R,S \ge 2$.
Let
\begin{equation}\label{eq:matr_ABtilde}
\widetilde{\matr{A}} = \matr{A}^{-\T} = \begin{bmatrix} \widetilde{\vect{a}}_1 & \cdots  & \widetilde{\vect{a}}_{\rkRow} \end{bmatrix}, \quad
\widetilde{\matr{B}} = \matr{B}^{-\T} = \begin{bmatrix} \widetilde{\vect{b}}_1 & \cdots  & \widetilde{\vect{b}}_{\rkCol} \end{bmatrix},
\end{equation}
so that $\{\widetilde{\vect{a}}_k\}^R_{k=1}$ $\{\widetilde{\vect{b}}_k\}^S_{k=1}$ are  the dual vectors for  $\{{\vect{a}}_k\}^R_{k=1}$ and $\{{\vect{b}}_k\}^S_{k=1}$, respectively.
Then the following lemma generalized  the derivations in  \Cref{sec:2x2_nonsymmetric}.
\begin{proposition}\label{prop:lker_unstructured}
Let $\tens{T}$ have a PT2D with invertible $\matr{A}$, $\matr{B}$.
For each pair $1 \le i, r \le \rkRow$, $1 \le j, s \le \rkCol$, define the following vector $\vect{v}^{(ijrs)} \in \sigma(\RR^{SRSR})$:
\begin{equation}\label{eq:lker_unstructured}
\vect{v}^{(ijrs)}:= \left(\theta^{(ijrs)}_1 \vecl{\widetilde{\vect{a}}_i \tensp \widetilde{\vect{b}}_j \tensp \widetilde{\vect{a}}_r \tensp \widetilde{\vect{b}}_s}+ 
\theta^{(ijrs)}_2 \vecl{\widetilde{\vect{a}}_i \tensp \widetilde{\vect{b}}_s \tensp \widetilde{\vect{a}}_r \tensp \widetilde{\vect{b}}_j}\right),
\end{equation}
where $\theta^{(ijrs)}_1,\theta^{(ijrs)}_2$ are defined in \eqref{eq:theta_expr}.
Then we have that:
\begin{enumerate}
\item  Each vector $\vect{v}^{(ijrs)}$ belongs to the left nullspace of $\matr{\Phi}(\tens{T})$;
\item Under  \Cref{as:AB_full_rank,as:nonzero_F,as:full_rank},  the symmetrizations of $\binom{R}{2}\binom{S}{2}$ vectors corresponding to distinct pairs of $(i,r)$ and $(j,s)$ form a basis of the symmetric left kernel of $\matr{\Phi}(\tens{T})$, i.e.,
\[
\symlker{\matr{\Phi}(\tens{T}}) = \Span{ \left\{\sigma_2 (\vect{v}^{(ijrs)} ) \right\}_{\substack{1 \le i < r \le R\\1 \le j < s \le S}}}
\]
\item The permutation/symmetrization map \eqref{eq:symperm} applied to $\vect{v}^{(ijrs)}$ gives the $R^2 \times S^2$ matrix of rank at most $1$
\[
\pi(\sigma_2 (\vect{v}^{(ijrs)} )) = \pi(\vect{v}^{(ijrs)}) = \left(\det{\matr{F}_{(i,r),(j,s)}} \right) \cdot
\sigma_2(\widetilde{\vect{a}}_i \kron \widetilde{\vect{a}}_r) (\sigma_2(\widetilde{\vect{b}}_j \kron \widetilde{\vect{b}}_s) )^{\T}
\]
\end{enumerate}
\end{proposition}
\begin{proof}
\begin{enumerate}
\item The proof repeats the derivations in \Cref{sec:2x2_nonsymmetric}, by substituting $\widetilde{\vect{a}}_i, \widetilde{\vect{a}}_r, \widetilde{\vect{b}}_j, \widetilde{\vect{b}}_s$ instead of $\widetilde{\vect{a}}_1, \widetilde{\vect{a}}_1, \widetilde{\vect{b}}_1, \widetilde{\vect{b}}_2$ in \eqref{eq:v_1122}.
In particular we use the fact that
\begin{equation}\label{eq:v_ijrs_multilinear}
\vect{v}^{(ijrs)} = 
 (\widetilde{\matr{B}} \kron  \widetilde{\matr{A}} \kron \widetilde{\matr{B}} \kron  \widetilde{\matr{A}})  
\left(\theta^{(ijrs)}_1 \vecl{{\vect{e}}_i \tensp {\vect{e}}_j \tensp {\vect{e}}_r \tensp {\vect{e}}_s}+ 
\theta^{(ijrs)}_2 \vecl{{\vect{e}}_i \tensp {\vect{e}}_s \tensp {\vect{e}}_r \tensp {\vect{e}}_j}\right).
\end{equation}
\item Thanks to \eqref{eq:v_ijrs_multilinear}, we have that $\sigma_2(\vect{v}^{(ijrs)}) =  (\widetilde{\matr{B}} \kron  \widetilde{\matr{A}} \kron \widetilde{\matr{B}} \kron  \widetilde{\matr{A}})  \vect{w}^{(ijrs)}$, where
\begin{align*}
\vect{w}^{(ijrs)} =& 
\frac{\theta^{(ijrs)}_1}{2} (\vecl{{\vect{e}}_i \tensp {\vect{e}}_j \tensp {\vect{e}}_r \tensp {\vect{e}}_s}+ 
\vecl{{\vect{e}}_r \tensp {\vect{e}}_s \tensp {\vect{e}}_i \tensp {\vect{e}}_j}) + \\
&\frac{\theta^{(ijrs)}_2}{2} (\vecl{{\vect{e}}_i \tensp {\vect{e}}_s \tensp {\vect{e}}_r \tensp {\vect{e}}_j}+ \vecl{{\vect{e}}_r \tensp {\vect{e}}_j \tensp {\vect{e}}_i \tensp {\vect{e}}_s}).
\end{align*}
We note that  we take only unique pairs of indices $1 \le i < r \le R$, $1 \le j < s \le S$,
hence for another two pairs $1 \le i' < r' \le R$, $1 \le j' < s' \le S$, such that $(i,r) \neq (i',r')$ or $(j,s) \neq (j',s')$, we have  $\langle\vect{w}^{(ijrs)},\vect{w}^{(i'j'r's')} \rangle = 0$.
Therefore, the $\vect{w}^{(ijrs)}$ are mutually ortogonal and the vectors $\{\vect{v}^{(ijrs)}\}_{\substack{1 \le i < r \le R\\1 \le j < s \le S}}$ are linearly independent.
\item Consider $\pi$ applied to $\vect{w}^{(ijrs)}$, which then becomes
\begin{align*}
\pi(\vect{w}^{(ijrs)}) = 
\frac{\theta^{(ijrs)}_1 +\theta^{(ijrs)}_2}{4}\Big(&
({\vect{e}}_i \kron {\vect{e}}_r)({\vect{e}}_j \kron {\vect{e}}_s)^{\T}+ 
({\vect{e}}_r \kron {\vect{e}}_i)( {\vect{e}}_s \kron {\vect{e}}_j )^{\T}+  \\
&({\vect{e}}_i\kron {\vect{e}}_r)({\vect{e}}_s \kron {\vect{e}}_j)^{\T}+
({\vect{e}}_r  \kron {\vect{e}}_i)(  {\vect{e}}_j\kron {\vect{e}}_s )^{\T}\Big) \\
 =\,& \sigma_2({\vect{e}}_i \kron {\vect{e}}_r) (\sigma_2({\vect{e}}_j \kron {\vect{e}}_s))^{\T}.
\end{align*}
By multilinearity we have 
\[
\pi((\widetilde{\matr{B}} \kron  \widetilde{\matr{A}} \kron \widetilde{\matr{B}} \kron  \widetilde{\matr{A}})   \vect{w}^{(ijrs)}) = (\widetilde{\matr{A}} \kron  \widetilde{\matr{A}}) \pi(\vect{w}^{(ijrs)}) (\widetilde{\matr{B}} \kron \widetilde{\matr{B}})^{\T},
\]
which completes the proof.
\end{enumerate}
\end{proof}

\subsection{Extracting ${A}$ and ${B}$ from the kernel}
We begin by proving the following lemma.
\begin{lemma}\label{lem:Phi_rank_GGHH}
For $\tens{T} = \ptd{\matr{A}}{\matr{B}}{\matr{F}}{\matr{G}}{\matr{H}} \in \RR^{R\times S\times \edim{3}}$ with $\matr{A} \in \RR^{R\times R}$ and $\matr{B} \in \RR^{S\times S}$ we have:
\begin{equation}\label{eq:Phi_rank_GGHH}
\rank{\matr{\Phi}(\tens{T})} \le \rank{\matr{\CoreVert} \kr \matr{\CoreVert} \kr \matr{\CoreHoriz}  \kr  \matr{\CoreHoriz}}.
\end{equation}
\end{lemma}
\begin{proof}
By the properties of vectorizations,  $\vecl{\tens{T}_{:,:,k}} = (\matr{B} \kron \matr{A}) \Diag{\vecl{\matr{F}}} (\vect{h}_k \kron \vect{g}_k)$,
hence the columns $(\matr{\Phi}(\tens{T}))_{:,k} = \vecl{\tens{T}_{:,:,k}} \kron \vecl{\tens{T}_{:,:,k}}$ can be expressed as
\[
(\matr{\Phi}(\tens{T}))_{:,k} = (\matr{B} \kron \matr{A} \kron\matr{B} \kron \matr{A}) \Diag{\vecl{\matr{F}} \kron  \vecl{\matr{F}}} (\vect{h}_k \kron \vect{g}_k \kron \vect{h}_k \kron \vect{g}_k),   
\]
which is equivalent to saying that
\begin{equation}\label{eq:Phi_explicit_form}
 \matr{\Phi}(\tens{T}) = (\matr{B} \kron \matr{A} \kron\matr{B} \kron \matr{A})
\Diag{\vecl{\matr{F}} \kron \vecl{\matr{F}}} (\matr{\CoreHoriz} \kr \matr{\CoreVert} \kr \matr{\CoreHoriz}  \kr \matr{\CoreVert}).
\end{equation}
Then the statement follows since ${(\matr{B} \kron \matr{A} \kron\matr{B} \kron \matr{A}) \Diag{\vecl{\matr{F}} \kron  \vecl{\matr{F}}}}$ is a square matrix and
\[
\rank{\matr{\CoreHoriz} \kr \matr{\CoreVert} \kr \matr{\CoreHoriz}  \kr \matr{\CoreVert}} = \rank{\matr{\CoreVert} \kr \matr{\CoreVert} \kr \matr{\CoreHoriz}  \kr \matr{\CoreHoriz}}.
\]
\end{proof}

Lemma \Cref{lem:Phi_rank_GGHH} serves as a base for proving \Cref{lem:Phi_rank_bound_nonsymmetric}.
\begin{proof}[Proof of \Cref{lem:Phi_rank_bound_nonsymmetric}]
\begin{enumerate}
\item  The columns of $\matr{\CoreVert} \kr \matr{\CoreVert} \kr \matr{\CoreHoriz}  \kr  \matr{\CoreHoriz}$ are vectorizations of tensors $(\vect{g}_k \tensp \vect{g}_k)\tensp ( \vect{h}_k \tensp \vect{h}_k)$,
hence 
\[
\range{\matr{\CoreVert} \kr \matr{\CoreVert} \kr \matr{\CoreHoriz}  \kr  \matr{\CoreHoriz}}\subseteq \vecl{S^2(\RR^{R}) \tensp S^2(\RR^{S}) }.
\]
This implies by \Cref{lem:Phi_rank_GGHH},
\begin{equation}\label{eq:inequalities_proof_Phi_rank}
\rank{\matr{\Phi}(\tens{T})} \le \rank{{\matr{\CoreVert} \kr \matr{\CoreVert} \kr \matr{\CoreHoriz}  \kr  \matr{\CoreHoriz}}} \le \dim\left\{S^2(\RR^{R}) \tensp S^2(\RR^{S})\right\} = \binom{R+1}{2} \binom{S+1}{2},
\end{equation}
\item Let \Cref{as:AB_full_rank,as:nonzero_F} be satisfied.
Then, by inspecting the proof of \Cref{lem:Phi_rank_GGHH}, we see that in \eqref{eq:Phi_explicit_form},   $(\matr{B} \kron \matr{A} \kron\matr{B} \kron \matr{A}) \Diag{\vecl{\matr{F}} \kron  \vecl{\matr{F}}}$ is an invertible matrix and hence equality  
$\rank{\matr{\Phi}(\tens{T})} = \rank{{\matr{\CoreVert} \kr \matr{\CoreVert} \kr \matr{\CoreHoriz}  \kr  \matr{\CoreHoriz}}}$ is achieved.
Therefore \eqref{eq:inequalities_proof_Phi_rank} becomes an equality if, in addition, \Cref{as:full_rank} holds true.

The dimension of the symmetrized kernel follows from a simple count of dimensions, as
\begin{align*}
\dim(\symlker{\matr{\Phi}(\tens{T})}) &= \dim(S^{2}(\RR^{RS})) - \rank{\matr{\Phi}(\tens{T})}  \\
&= \binom{RS+1}{2} -   \binom{R+1}{2} \binom{S+1}{2} =   \binom{R}{2} \binom{S}{2}.
\end{align*}

\item
Assume,  that \Cref{as:nonzero_F,as:full_rank} are not satisfied simultaneously. We have the following cases:
\begin{itemize}
\item If \Cref{as:full_rank} is not satisfied, then the strict inequality in \eqref{eq:inequalities_proof_Phi_rank} follows.

\item If \Cref{as:nonzero_F} is not satisfied,  there is $i,j$ such that $F_{i,j} = 0$, hence the (only) diagonal element in  $\Diag{\vecl{\matr{F}} \kron  \vecl{\matr{F}}}$  corresponding to $F_{i,j}^2$ is zero (this diagonal element occurs at the position of $1$ in the vector $\vect{v} = \vect{e}_{j} \kron \vect{e}_{i} \kron \vect{e}_{j} \kron \vect{e}_{i}$).
Denote $\Pi = \matr{I}_{(RS)^2} - \vect{v}\vect{v}^{\T}$ (the identity matrix without the corresponding diagonal element).
Then
\begin{align*}
\rank{\matr{\Phi}(\tens{T})} & \le \rank{\Diag{\vecl{\matr{F}} \kron \vecl{\matr{F}}} (\matr{\CoreHoriz} \kr \matr{\CoreVert} \kr \matr{\CoreHoriz}  \kr \matr{\CoreVert})} \\
& \le \rank{\matr{\Pi}(\matr{\CoreHoriz} \kr \matr{\CoreVert} \kr \matr{\CoreHoriz}  \kr \matr{\CoreVert})} < \rank{(\matr{\CoreHoriz} \kr \matr{\CoreVert} \kr \matr{\CoreHoriz}  \kr \matr{\CoreVert})},
\end{align*}
where the strict inequality follows from the fact  that $\matr{\Pi} ( (\vect{h} \kron \vect{g} \kron \vect{h} \kron \vect{g}))$ does not contain the element $h^2_i g^2_i$ which appears only once in  $\vect{h} \kron \vect{g} \kron \vect{h} \kron \vect{g}$, for all  $\vect{g} \in \RR^{R}$, $\vect{h} \in \RR^{S}$.
\end{itemize}
\end{enumerate}
\end{proof}

\begin{proof}[Proof of \Cref{thm:pt2d_nonsymmetric}]
\begin{enumerate}
\item
Let $\matr{P}_1, \ldots, \matr{P}_M$  be as in  \Cref{thm:pt2d_nonsymmetric}. Then we have that, by \Cref{prop:lker_unstructured}
\[
\Span{\{\matr{P}_k\}^M_{k=1}} = \Span{\left\{\left(\det{\matr{F}_{(i,r),(j,s)}} \right) \cdot
\sigma_2(\widetilde{\vect{a}}_i \kron \widetilde{\vect{a}}_r) (\sigma_2(\widetilde{\vect{b}}_j \kron \widetilde{\vect{b}}_s) )^{\T}\right\}_{\substack{1 \le i \le r \le R\\1 \le j \le s \le S}}}.
\]
Therefore, we have that
\begin{equation}\label{eq:subspace_symmetric_A}
\range{\matr{P}_A}  =\Span{\{\range{\matr{P}_k}\}^M_{k=1}} \subseteq \Span{ \{\sigma_2(\widetilde{\vect{a}}_i \kron \widetilde{\vect{a}}_r)\}_{{1 \le i \le r \le R}}},
\end{equation}
and, by the dimension count,
\begin{equation}\label{eq:subspace_symmetric_A_rank}
\dim \left( \symlker{\matr{P}_A} \right) = \dim S^2(\RR^R) - \rank{\matr{P}_A} \ge \binom{R+1}{2} - \binom{R}{2} = R.
\end{equation}
If $\krank{\matr{F}^{\T}} \ge 2$, for any pair $(i,r)$ there exist $(j,s)$ such that $\det{\matr{F}_{(i,r),(j,s)}} \neq 0$, then this implies that the equality holds in \eqref{eq:subspace_symmetric_A} and \eqref{eq:subspace_symmetric_A_rank}.
On the other hand, if $\krank{\matr{F}^{\T}} =1$, then there exists a pair $(i,j)$, such that $\sigma_2(\widetilde{\vect{a}}_i \kron \widetilde{\vect{a}}_r) \not\in \range{\matr{P}_A} $, hence  $\dim \left( \symlker{\matr{P}_A} \right) > R$ in this case.
Finally we note that for any $i,r,\ell$
\[
\langle \sigma_2(\widetilde{\vect{a}}_i \kron \widetilde{\vect{a}}_r), \vect{a}_{\ell} \kron \vect{a}_{\ell} \rangle  = \langle\widetilde{\vect{a}}_i,\vect{a}_{\ell}\rangle \langle\widetilde{\vect{a}}_r,\vect{a}_{\ell}\rangle  = 0.
\]
Therefore,
\[
\Span{\{ \vect{a}_{\ell} \kron \vect{a}_{\ell}\}_{\ell=1}^R} \subseteq \symlker{\matr{P}_A},
\]
where the equality holds if and only if $\krank{\matr{F}^{\T}} \ge 2$.

\item
The proof is analogous to the one of the previous statement.
We observe that
\begin{equation}\label{eq:subspace_symmetric_B}
\range{\matr{P}_B} \subseteq \Span{ \{\sigma_2(\widetilde{\vect{b}}_j \kron \widetilde{\vect{b}}_s)\}_{{1 \le j \le s \le S}}}.
\end{equation}
and 
\[
\Span{\{ \vect{b}_{\ell} \kron \vect{b}_{\ell}\}_{\ell=1}^R} \subseteq \symlker{\matr{P}_B},
\]
where the equalities hold if and only if $\krank{\matr{F}} = 2$, by invoking the dimension count argument.
\end{enumerate}
\end{proof}

\subsection{On the case of collinear columns or rows of F}\label{sec:collinear_jbd}
We will make some comments on the case of collinear rows or columns in $\matr{F}$ (see  \Cref{sec:collinearities}).
In fact, if \Cref{as:krank_F} is not satisfied, then we will show that $\matr{A}$ and $\matr{B}$ can be recovered subject to additional ambiguities  outlined in \Cref{sec:collinearities}.
To do this, we use the following proposition
\begin{proposition}
Let  \cref{as:AB_full_rank,as:nonzero_F,as:full_rank} be satisfied, $\matr{F}$ be as in \eqref{eq:F_repeating_columns} (see \Cref{lem:F_repeating_columns}).
Then for $\matr{P}_A$ and $\matr{P}_B$ as in \Cref{thm:pt2d_nonsymmetric}, it holds that
\begin{align}
\dim \left( \symlker{\matr{P}_A} \right)  & =
\Span{(\{ \matr{A} \Diag{\matr{Y}_{1} ,\ldots,\matr{Y}_{R'}} \matr{A}^{\T}: \matr{Y}_k  \text{ is }m_k \times m_k \text{ symmetric}\})}, \label{eq:P1_kernel_blocks}\\
\dim \left( \symlker{\matr{P}_B} \right) & =
\Span{(\{ \matr{B} \Diag{\matr{Z}_{1} ,\ldots,\matr{Z}_{S'}} \matr{B}^{\T}: \matr{Z}_k  \text{ is }n_k \times n_k \text{ symmetric}\})}.\label{eq:P2_kernel_blocks}
\end{align}
\end{proposition}
\begin{proof}
We will prove the statement just for $\matr{A}$.
We define $\set{I}_k$, $1 \le k \le R'$, such that
\begin{align*}
\set{I}_1 & = \{1,\ldots ,m_1\} \\
\set{I}_2 & = \{m_1+1,\ldots ,m_1+m_2\} \\
& \vdots \\
\set{I}_{R'} & = \{\sum\limits_{\ell=1}^{R'} m_{\ell}+1,\ldots ,R\}. 
\end{align*}
The sets of indices $\set{I}_k$ correspond to blocks of rows in  \eqref{eq:F_repeating_columns}.
Then, as in the proof of \Cref{thm:pt2d_nonsymmetric}, $\sigma_2(\widetilde{\vect{a}}_i \kron \widetilde{\vect{a}}_r)$ belongs to $\range{\matr{P}_A}$ if and only if $\rank{\matr{F}_{(i,r),:}} = 2$,
which can only happen in $i$ and $r$ belong to different blocks of rows.
Therefore, we have that 
\begin{equation}\label{eq:subspace_symmetric_A_blocks}
\range{\matr{P}_A} = \Span{ \{\sigma_2(\widetilde{\vect{a}}_i \kron \widetilde{\vect{a}}_r)\}_{\substack{i \in \set{I}_k, r \in \set{I}_{\ell} \\ k \neq \ell}}}.
\end{equation}
Note that the dimension of this space is equal to 
\[
\dim( \range{\matr{P}_A})  = \binom{R+1}{2} - \left(\sum\limits_{\ell=1}^{S'} \binom{m_{\ell}+1}{2}\right).
\]
as we only forbid the pairs belonging to the same blocks.
Next, we observe that for any $(i',r')$ and $(i,r)$ 
\begin{equation}\label{eq:crossprod}
\langle\sigma_2({\vect{a}}_{i'} \kron {\vect{a}}_{r'}),\sigma_2(\widetilde{\vect{a}}_{i} \kron \widetilde{\vect{a}}_{r})\rangle =
 \frac{1}{2} \left(\langle\vect{a}_{i'},\widetilde{\vect{a}}_{i}\rangle\langle\vect{a}_{r'},\widetilde{\vect{a}}_{r}\rangle + \langle\vect{a}_{i'},\widetilde{\vect{a}}_{r}\rangle\langle\vect{a}_{r'},\widetilde{\vect{a}}_{i}\rangle\right).
\end{equation}
Now, if $i',r' \in \set{I}_t$ and  $i \in \set{I}_k, r \in \set{I}_{\ell}$, $k \neq \ell$, then we have that
\eqref{eq:crossprod} vanishes because $(i',r') \not\in \{(i,r),(r,i)\}$ due to the constraints.
This implies that for a fixed $1 \le t \le R$, vectorization of any matrix
\[
\matr{X} = \sum\limits_{i,r\in \set{I}_t} c_{i,j} \vect{a}_{i}\vect{a}_{r}^{\T}.
\]
belongs to the left kernel:  $\vecl{\matr{X}} \in \lker{\matr{P}_A}$.
Thus any matrix of the form as in \eqref{eq:P1_kernel_blocks} belongs to $\symlker{\matr{P}_A}$ and the two sets in \eqref{eq:P1_kernel_blocks} coincide by the dimension count.
\end{proof}

\section{Symmetric case: proof of the main theorem}\label{sec:symmetric_proofs}

\subsection{Basis of the left kernel}
Let $\vect{a}_k$ and $\widetilde{\vect{a}}_k$ be in \Cref{sec:nonsymmetric_proofs}.
In the symmetric case we have   $\vect{a}_k = \vect{b}_k$ and $\widetilde{\vect{a}}_k = \widetilde{\vect{b}}_k$.
Note that the vector $\vect{v}^{(ijrs)}$ belongs to the left kernel of $\matr{\Phi}(\tens{T})$ and in this case becomes
\begin{equation}\label{eq:lker_unstructured_summetric}
\vect{v}^{(ijrs)}:= \left(\theta^{(ijrs)}_1 \vecl{\widetilde{\vect{a}}_i \tensp \widetilde{\vect{a}}_j \tensp \widetilde{\vect{a}}_r \tensp \widetilde{\vect{a}}_s}+ 
\theta^{(ijrs)}_2 \vecl{\widetilde{\vect{a}}_i \tensp \widetilde{\vect{a}}_s \tensp \widetilde{\vect{a}}_r \tensp \widetilde{\vect{a}}_j}\right).
\end{equation}
Let $\sigma_{(2,2)}(\cdot)$ be the orthogonal projection on $\vecl{S^2(\vecl{S^2(\RR^S)})}$.
Then we have that 
\begin{align*}
\vect{u}^{(ijrs)} := \sigma_{(2,2)}(\vect{v}^{(ijrs)}) &=
\theta^{(ijrs)}_1 \sigma_2(\sigma_2(\widetilde{\vect{a}}_i \kron \widetilde{\vect{a}}_j) \kron \sigma_2(\widetilde{\vect{a}}_r \kron \widetilde{\vect{a}}_s))+ 
\theta^{(ijrs)}_2 \sigma_2(\sigma_2(\widetilde{\vect{a}}_i \kron \widetilde{\vect{a}}_s) \kron \sigma_2(\widetilde{\vect{a}}_r \tensp \widetilde{\vect{a}}_j)).
\end{align*}
We have the following proposition.
\begin{proposition}\label{prop:lker_unstructured_symmetric}
Let $\tens{T}$ have a symmetric PT2D with invertible $\matr{A}=\matr{B}$.
Then we have that:
\begin{enumerate}
\item Under  \Cref{as:AB_full_rank,as:nonzero_F,as:full_rank_symmetric},  the symmetric left kernel is spanned by the vectors $\vect{u}^{(ijrs)} $
\[
\symlkerd{\matr{\Phi}(\tens{T})} =\Span{ \{ \vect{u}^{(ijrs)} \}_{\substack{1 \le i < r \le R \\1 \le j < s \le R}}}
\]
\item The full symmetrization  applied to $\vect{u}^{(ijrs)}$ gives the  following
\[
\sigma_4 (\vect{u}^{(ijrs)} ) = \sigma_4(\vect{v}^{(ijrs)} ) = \left(\det{\matr{F}_{(i,r),(j,s)}} \right) \cdot
\sigma_4(\widetilde{\vect{a}}_i \kron \widetilde{\vect{a}}_r\kron \widetilde{\vect{a}}_j \kron \widetilde{\vect{a}}_s) .
\]
\end{enumerate}
\end{proposition}
\begin{proof}
\begin{enumerate}
\item  Note that we have
\[
\Span{ \{ \vect{u}^{(ijrs)} \}_{\substack{1 \le i < r \le R \\1 \le j < s \le R}}} \subseteq
\symlkerd{\matr{\Phi}(\tens{T})}.
\]
To show that the inequality holds, we will show that we can choose $N = \frac{R^2(R-1)(R+1)}{12}$ linearly independent vectors among $\vect{u}^{(ijrs)}$.
Consider the following sets of indices:
\begin{align}
\set{I}_2 &= \{ (k,k,\ell,\ell)   \}_{1\le k < \ell \le R}, \\
\set{I}_3 &= \{ (k,\ell,m,m), (k,\ell,\ell,m) , (k,k,\ell,m)  \}_{1 \le k <\ell<m \le R }, \\
\set{I}_{4} & = \{(k,\ell,m,n),(k,m,\ell,n)\}_{1\le k < \ell < m < n \le R}, 
\end{align}
Next, as in the proof \Cref{prop:lker_unstructured}, we express $\vect{u}^{(ijrs)} =  (\widetilde{\matr{A}} \kron  \widetilde{\matr{A}} \kron \widetilde{\matr{A}} \kron  \widetilde{\matr{A}})  \vect{x}^{(ijrs)}$, where
\begin{align*}
&\vect{x}^{(ijrs)} 
= \theta^{(ijrs)}_1 \sigma_2(\sigma_2({\vect{e}}_i \kron {\vect{e}}_j) \kron \sigma_2({\vect{e}}_r \kron {\vect{e}}_s))+ 
\theta^{(ijrs)}_2 \sigma_2(\sigma_2({\vect{e}}_i \kron {\vect{e}}_s) \kron \sigma_2({\vect{e}}_r \tensp {\vect{e}}_j)) \\
&\quad=  \frac{\theta^{(ijrs)}_1}{4} \left(({\vect{e}}_i \kron {\vect{e}}_j + {\vect{e}}_j \kron {\vect{e}}_i) \kron ({\vect{e}}_r \kron {\vect{e}}_s+{\vect{e}}_s \kron {\vect{e}}_r) +
 ({\vect{e}}_r \kron {\vect{e}}_s+{\vect{e}}_s \kron {\vect{e}}_r) \kron ({\vect{e}}_i \kron {\vect{e}}_j + {\vect{e}}_j \kron {\vect{e}}_i)  \right) \\
 &\quad=  \frac{\theta^{(ijrs)}_2}{4} \left(({\vect{e}}_i \kron {\vect{e}}_s + {\vect{e}}_s \kron {\vect{e}}_i) \kron ({\vect{e}}_r \kron {\vect{e}}_j+{\vect{e}}_j \kron {\vect{e}}_r) +
 ({\vect{e}}_r \kron {\vect{e}}_j+{\vect{e}}_j \kron {\vect{e}}_r) \kron ({\vect{e}}_i \kron {\vect{e}}_s + {\vect{e}}_s \kron {\vect{e}}_i)  \right). 
\end{align*}
We note that all the multiindices in $\set{I}_2,\set{I}_3,\set{I}_{4}$ correspond to different  choices  of $4$ numbers (with possible repetition) among $\{1,\ldots, R\}$ (except the set $\set{I}_4$ which is composed of couple of multiindices, corresponding to the same subset of numbers).
Therefore, for two distinct multi-indices $(i,j,r,s), (i',j',r',s') \in \set{I}_2 \cup \set{I}_3 \cup\set{I}_{4}$, we have $\langle \vect{x}^{(ijrs)},  \vect{x}^{(i'j'r's')}\rangle =0$ if 
\[
(i,j,r,s) \in \set{I}_t, (i',j',r',s') \in \set{I}_t', \text{ with } (t,t') \in \{ (2,2),(2,3),(2,4), (3,3), (3,4) \}.
\]
This proves that the vectors in $\{ \vect{x}^{(ijrs)} \}_{\{(i,j,r,s) \in \set{I}_2 \cup\set{I}_3  \}}$ are mutually orthogonal and that 
\[
\{ \vect{x}^{(ijrs)} \}_{(i,j,r,s) \in \set{I}_2 \cup\set{I}_3  } \bot \{ \vect{x}^{(ijrs)} \}_{{(i,j,r,s) \in \set{I}_4  }}.
\]
Next, we are going to show that $\{ \vect{x}^{(ijrs)} \}_{\{(i,j,r,s) \in \set{I}_4  \}}$ is linearly independent.
For this we note that for two multi-indices $(i,j,r,s), (i',j',r',s') \in \set{I}_4$ we have $\langle \vect{x}^{(ijrs)},  \vect{x}^{(i'j'r's')}\rangle =0$ as long as they are not linked by a permutation (i.e., $\{i,j,r,s\} \neq  \{i',j',r',s'\}$).
Therefore, it is left to show that $\vect{x}^{(k{\ell}mn)}$ is linearly independent of  $ \vect{x}^{(km{\ell}n)}$.
For simplicity, take $(k,{\ell},m,n) = (1,2,3,4)$  and note that  the term $\vect{e}_{1} \kron \vect{e}_{2} \kron \vect{e}_{3} \kron \vect{e}_{4}$ appears only in the expansion of $\vect{x}^{(1234)}$ while $\vect{e}_{1} \kron \vect{e}_{3} \kron \vect{e}_{2} \kron \vect{e}_{4}$ appears only in the expansion of  $\vect{x}^{(1324)}$.
This proves that $\vect{x}^{(1234)}$ and $\vect{x}^{(1324)}$ are linearly independent.

Finally, we note that the cardinality of the set
\[
\# \{\set{I}_2 \cup\set{I}_2 \cup\set{I}_3  \} = \binom{R}{2} + 3  \binom{R}{3} + 2  \binom{R}{4} = 
\frac{R^2(R-1)(R+1)}{12} =  \dim\left( \symlkerd{\matr{\Phi}(\tens{T})}\right).
\]
where the last equality holds from 
\Cref{lem:Phi_rank_bound_symmetric}.
This proves that  $\{ \vect{u}^{(ijrs)} \}_{\{(i,j,r,s) \in \set{I}_2 \cup\set{I}_3  \cup\set{I}_4  \}}$ forms the basis of  $\symlkerd{\matr{\Phi}(\tens{T})}$.
\item The proof is similar to the corresponding statement of \Cref{prop:lker_unstructured}. Note that $\vecl{S^{4}(\RR^{R})} \subset \vecl{S^2(\vecl{S^2(\RR^R)})}$ and therefore 
\[
\sigma_4 (\vect{u}^{(ijrs)} ) = \sigma_4(\sigma_{(2,2)}(\vect{v}^{(ijrs)}) ) = \sigma_4 (\vect{v}^{(ijrs)} ).
\]
Then it is easy to see that
\begin{align*}
\sigma_4(\vect{v}^{(ijrs)} ) &= \theta^{(ijrs)}_1 \sigma_4\left(\vecl{\widetilde{\vect{a}}_i \tensp \widetilde{\vect{a}}_j \tensp \widetilde{\vect{a}}_r \tensp \widetilde{\vect{a}}_s}\right)+ 
\theta^{(ijrs)}_2 \sigma_4 \left(\vecl{\widetilde{\vect{a}}_i \tensp \widetilde{\vect{a}}_s \tensp \widetilde{\vect{a}}_r \tensp \widetilde{\vect{a}}_j}\right) \\
&= ( \theta^{(ijrs)}_1 + \theta^{(ijrs)}_2)\sigma_4\left(\vecl{\widetilde{\vect{a}}_i \tensp \widetilde{\vect{a}}_j \tensp \widetilde{\vect{a}}_r \tensp \widetilde{\vect{a}}_s}\right) = 
\left(\det{\matr{F}_{(i,r),(j,s)}} \right)  \sigma_4\left({\widetilde{\vect{a}}_i \kron \widetilde{\vect{a}}_j \kron \widetilde{\vect{a}}_r \kron \widetilde{\vect{a}}_s}\right).
\end{align*}
\end{enumerate}
\end{proof}

\subsection{Extracting the vectors from the left kernel}

\begin{proof}[Proof of \Cref{lem:Phi_rank_bound_symmetric}]
\begin{enumerate}
 \item Similarly, the columns of $\matr{\CoreVert}^{\kr^4} = \matr{\CoreVert} \kr \matr{\CoreVert} \kr \matr{\CoreVert}  \kr  \matr{\CoreVert}$ are vectorizations of symmetric tensors $(\vect{g}^{\tensp 4})$,
hence by \Cref{lem:Phi_rank_GGHH} we have
\begin{equation}\label{eq:inequalities_proof_Phi_rank_symmetric}
\rank{\matr{\Phi}(\tens{T})} \le \rank{{\matr{\CoreVert}^{\kr 4} }} \le \dim\{S^4(\RR^{R})\} = \binom{R+3}{4},
\end{equation}
The second part of the statement follows from the dimension count:
\begin{align*}
& \dim\left( \symlkerd{\matr{\Phi}(\tens{T})}\right) = 
\dim\left(\vecl{S^2(\vecl{S^2(\RR^R)})}\right) - \rank{\matr{\Phi}(\tens{T})} \\
&= \frac{(\binom{R+1}{2} +1)\binom{R+1}{2}}{2} - \binom{R+3}{4} = \frac{R^2(R-1)(R+1)}{12}.
\end{align*}
\item As in the nonsymmetric case, both inequalities in \eqref{eq:inequalities_proof_Phi_rank_symmetric} become equalities if  all \Cref{as:AB_full_rank,as:nonzero_F,as:full_rank_symmetric} are satisfied.
\item
We proceed as in the nonsymmetric case. Assume  that \Cref{as:nonzero_F,as:full_rank_symmetric} are not satisfied simultaneously. We have the following cases:
\begin{itemize}
\item If \Cref{as:full_rank_symmetric} is not satisfied, then the strict inequality in \eqref{eq:inequalities_proof_Phi_rank_symmetric} follows.

\item If \Cref{as:nonzero_F} is not satisfied,  there is $i,j$ such that $F_{j,i} = F_{i,j} = 0$.
Then we have that $F_{j,i}F_{i,i} = F_{i,j}F_{i,i}  = 0$ and all the corresponding diagonal elements  in $\Diag{\vecl{\matr{F}} \kron  \vecl{\matr{F}}}$ (there are $4$ of them if $i\neq j$ and just one otherwise) are zero.
This implies that in the product
 $\Diag{\vecl{\matr{F}} \kron \vecl{\matr{F}}}  \vect{g}^{\kron 4}$ the monomial
$g^3_i g_j$ never appears and thus the corresponding rows in
\[
\Diag{\vecl{\matr{F}} \kron \vecl{\matr{F}}} \matr{\CoreVert}^{\kr 4}
\]
are zero. Therefore,
\begin{align*}
\rank{\matr{\Phi}(\tens{T})} & \le \rank{\Diag{\vecl{\matr{F}} \kron \vecl{\matr{F}}} (\matr{\CoreVert}^{\kr 4})} < \rank{(\matr{\CoreVert}^{\kr 4})},
\end{align*}
which completes the proof.
\end{itemize}
\end{enumerate}
\end{proof}

\begin{proof}[Proof of \Cref{thm:pt2d_symmetric}]
Let $\matr{p}_1, \ldots, \matr{p}_N$  be as in  \Cref{thm:pt2d_symmetric}. Then we have that, by \Cref{prop:lker_unstructured_symmetric}
\[
\Span{\{\matr{p}_k\}^N_{k=1}} = \Span{\left\{\left(\det{\matr{F}_{(i,r),(j,s)}} \right) \cdot
\sigma_4(\widetilde{\vect{a}}_i \kron \widetilde{\vect{a}}_r \kron \widetilde{\vect{a}}_j \kron \widetilde{\vect{a}}_s) )^{\T}\right\}_{\substack{1 \le i <r \le R\\1 \le j <s \le R}}}.
\]
Now, let $\matr{P}_k = \matricize{R^3}{R}{\matr{p}_k}$ and $\matr{P}_{sym}$ as in \eqref{eq:Psym}.
 Then we have that
\begin{equation}\label{eq:proof_range_Psym}
\range{\matr{P}_{sym}} =  \Span{\left\{\range{\matr{Z}^{(ijrs)} } \right\}_{\substack{(i,j,r,s) :\\ \det{\matr{F}_{(i,r),(j,s)}\neq 0}}}}
\end{equation}
where $\matr{Z}^{(ijrs)}  = \matricize{R^3}{R}{\sigma_4({\widetilde{\vect{a}}_i \kron \widetilde{\vect{a}}_j \kron \widetilde{\vect{a}}_r \kron \widetilde{\vect{a}}_s}) }$
Now let us expand $\sigma_4({\widetilde{\vect{a}}_i \kron \widetilde{\vect{a}}_j \kron \widetilde{\vect{a}}_r \kron \widetilde{\vect{a}}_s})$.
Since the order of $(i,j,r,s)$ does not not matter, we will consider just the three cases 
\[
(i,j,r,s) \in \{(1,1,2,2), (1,1,2,3), (1,2,3,4)\}.
\]
Direct calculation shows that 
\begin{align*}
\matr{Z}^{(1122)} &=
\frac{1}{2}
 \sigma_3(\widetilde{\vect{a}}_1 \kron \widetilde{\vect{a}}_1 \kron \widetilde{\vect{a}}_2) ({\widetilde{\vect{a}}_2})^{\T} +
 \frac{1}{2}  \sigma_3(\widetilde{\vect{a}}_1 \kron \widetilde{\vect{a}}_2 \kron \widetilde{\vect{a}}_2)  (\widetilde{\vect{a}}_1)^{\T} \\
 &=
\frac{1}{2} \begin{bmatrix} \sigma_3(\widetilde{\vect{a}}_1 \kron \widetilde{\vect{a}}_2 \kron \widetilde{\vect{a}}_2) &
\sigma_3(\widetilde{\vect{a}}_1 \kron \widetilde{\vect{a}}_1 \kron \widetilde{\vect{a}}_2)
 \end{bmatrix}
\begin{bmatrix} \widetilde{\vect{a}}_1 & \widetilde{\vect{a}}_2 \end{bmatrix}^{\T},
\end{align*}
hence $\range{\matr{Z}^{(1122)}} = \Span{\widetilde{\vect{a}}_1 \kron \widetilde{\vect{a}}_1 \kron \widetilde{\vect{a}}_2, \widetilde{\vect{a}}_1 \kron \widetilde{\vect{a}}_2 \kron \widetilde{\vect{a}}_2}$.
Similarly
\[
\matr{Z}^{(1123)} =
\frac{1}{4}
 \sigma_3(\widetilde{\vect{a}}_1 \kron \widetilde{\vect{a}}_1 \kron \widetilde{\vect{a}}_2) ({\widetilde{\vect{a}}_3})^{\T} +
 \frac{1}{4}  \sigma_3(\widetilde{\vect{a}}_1 \kron \widetilde{\vect{a}}_1 \kron \widetilde{\vect{a}}_3)  (\widetilde{\vect{a}}_2)^{\T} +
 \frac{1}{2}  \sigma_3(\widetilde{\vect{a}}_1 \kron \widetilde{\vect{a}}_2 \kron \widetilde{\vect{a}}_3)  (\widetilde{\vect{a}}_1)^{\T},
\]
hence $\range{\matr{Z}^{(1123)} } = \Span{\widetilde{\vect{a}}_1 \kron \widetilde{\vect{a}}_1 \kron \widetilde{\vect{a}}_2, \widetilde{\vect{a}}_1 \kron \widetilde{\vect{a}}_1 \kron \widetilde{\vect{a}}_3,\widetilde{\vect{a}}_1 \kron \widetilde{\vect{a}}_2 \kron \widetilde{\vect{a}}_3}$.
Finally,
\[
\matr{Z}^{(1234)} =
\frac{1}{4} \begin{bmatrix}
\sigma_3(\widetilde{\vect{a}}_2 \kron \widetilde{\vect{a}}_3 \kron \widetilde{\vect{a}}_4)\,
\sigma_3(\widetilde{\vect{a}}_1 \kron \widetilde{\vect{a}}_3 \kron \widetilde{\vect{a}}_4)\,
\sigma_3(\widetilde{\vect{a}}_1 \kron \widetilde{\vect{a}}_2 \kron \widetilde{\vect{a}}_4)\, 
\sigma_3(\widetilde{\vect{a}}_1 \kron \widetilde{\vect{a}}_2 \kron \widetilde{\vect{a}}_3)\,
\end{bmatrix}
\begin{bmatrix} \widetilde{\vect{a}}_1 & \widetilde{\vect{a}}_2 & \widetilde{\vect{a}}_3 & \widetilde{\vect{a}}_4  \end{bmatrix}^{\T},
\]
thus $\range{\matr{Z}^{(1234)}} = \Span{\sigma_3(\widetilde{\vect{a}}_1 \kron \widetilde{\vect{a}}_2 \kron \widetilde{\vect{a}}_3), \sigma_3(\widetilde{\vect{a}}_1 \kron \widetilde{\vect{a}}_2 \kron \widetilde{\vect{a}}_4), \sigma_3(\widetilde{\vect{a}}_1 \kron \widetilde{\vect{a}}_3 \kron \widetilde{\vect{a}}_4),\sigma_3(\widetilde{\vect{a}}_2 \kron \widetilde{\vect{a}}_3 \kron \widetilde{\vect{a}}_4)}$.

Thus, we have that
\begin{equation}\label{eq:subspace_antisymmetric_3}
\range{\matr{P}_{sym}} \subseteq \Span{\{\sigma_3(\widetilde{\vect{a}}_i \kron \widetilde{\vect{a}}_j \kron \widetilde{\vect{a}}_r) : 1 \le i < r \le R, 1 \le j \le R \}},
\end{equation}
Assume that $\krank{\matr{F}} =\krank{\matr{F}}^{\T} \ge 2$. Fix $1 \le i< r \le R$ and $1 \le j \le R$, the submatrix $\matr{F}_{(i,r),:}$ has rank $2$.
Note that the $j$-th column $\matr{F}_{(i,r),j}$  of that matrix is nonzero (by \Cref{as:nonzero_F}), and therefore there must exist at least one other column $\matr{F}_{(i,r),s}$ that is not a multiple of $\matr{F}_{(i,r),j}$, which implies $\det \matr{F}_{(i,r),(j,s)} \neq 0$.
Therefore, in view of \eqref{eq:proof_range_Psym}, we have $\sigma_3(\widetilde{\vect{a}}_i \kron \widetilde{\vect{a}}_r \kron \widetilde{\vect{a}}_j)$, and hence the two spaces in \eqref{eq:subspace_antisymmetric_3} are equal due to the dimension count (which also implies $\dim \range{\matr{P}_{sym}}  = \dim(S^3(\RR^R)) - R$).

Now assume that $\krank{F} = 1$. This means that there exists $1 \le k< \ell \le R$ such that the submatrix $\matr{F}_{(k,\ell),:}$ is rank $1$.
In particular, $\det \matr{F}_{(k,\ell),(k,m)} = 0$ for any $m$.
This implies that 
\[
\range{\matr{P}_{sym}} \subseteq \Span{\{\sigma_3(\widetilde{\vect{a}}_i \kron \widetilde{\vect{a}}_j \kron \widetilde{\vect{a}}_r) : 1 \le i < r \le R, 1 \le j \le R, (i,r,j) \neq (k,\ell,k) \}},
\]
and hence the two linear spaces are not equal in \eqref{eq:subspace_antisymmetric_3}, and thus $\dim \range{\matr{P}_{sym}}  < \dim(S^3(\RR^R)) - R$.

Finally, we note that in case when equality holds in  \eqref{eq:subspace_antisymmetric_3},
\[
(\vect{a}_k \kron\vect{a}_k \kron \vect{a}_k )^{\T}\matr{P}_{sym} = 0, \quad k=1,\ldots, R,
\]
which completes the proof.
\end{proof}

\section{Numerical experiments}
A preliminary implementation is available in MATLAB at \url{https://github.com/kdu/paratuck-2}, together with the reproducible experiments.
We use the following metrics to evaluate the results  of the algorithms:
\begin{itemize}
\item squared Frobenius  error: 
\[
MSE(\tens{T},\widehat{\tens{T}}) = \|\tens{T} -\widehat{\tens{T}}\|_F^2
\]
 (to measure the reconstruction error);
\item sum of squared sines of angles between columns (to measure the correct recovery of the factors)
\[
SSS(\matr{A},\widehat{\matr{A}}) = \sum\limits_{k=1}^R 1 - \left(\frac{\langle\vect{a}_k,\widehat{\vect{a}}_k\rangle}{\|\vect{a}_k\|\|\widehat{\vect{a}}_k\|}\right)^2
\]
where the columns of $\widehat{\matr{A}}$ are permuted to match the columns of $\matr{A}$ beforehand.
The permutation is computed based on Hungarian algorithm.
Note that $SSS(\matr{A},\widehat{\matr{A}})  = \frac{\|\matr{A}-\widehat{\matr{A}}\|^2_F}{2}$ if the columns of $\matr{A}$ and $\widehat{\matr{A}}$ are normalized.
\end{itemize}

We note that in the implementation we use a basic joint eigenvalue computation algorithm, and thus the performance can be improved by using generalized Schur decompositions or approximate joint diagonalization techniques.

\subsection{Nonsymmetric case: simplest  example}
We consider the following simple test example:
\[
\matr{\CoreFront} = \begin{bmatrix}1 &1\\ 2 &-1\end{bmatrix},\quad
\matr{A} =  \begin{bmatrix}1& 1\\ -1& 2\end{bmatrix},\quad
\matr{B} = \begin{bmatrix}1 & 1\\  1 & -3\end{bmatrix},
\]
\[
\matr{\CoreVert} = \begin{bmatrix} -5& -4 &-3 &-2& -1& 0& 1& 2 &3 &4 \\ 1 &0& 2& 1& -3 &2 &-2 & -1& 0& 1\end{bmatrix},\quad
\matr{\CoreHoriz} = \begin{bmatrix} -5& -4 &-3 &-2& -1& 0& 1& 2 &3 &4 \\ 1 &1& 1& 1& 1 &1 &1 & 1& 1& 1\end{bmatrix}.
\]
The squared Frobenius norm of $\tens{T} = \ptd{\matr{A}}{\matr{B}}{\matr{F}}{\matr{G}}{\matr{H}}$ is $13538$.

\begin{figure}[hbt!]
\[
\includegraphics[height=5cm]{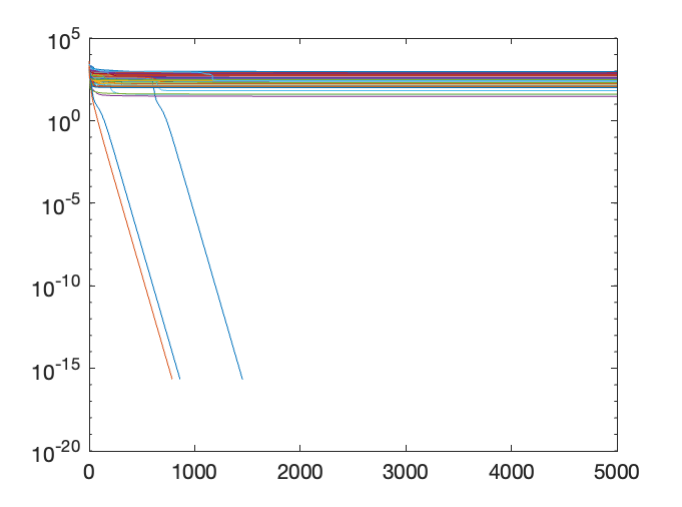} \quad \includegraphics[height=5cm]{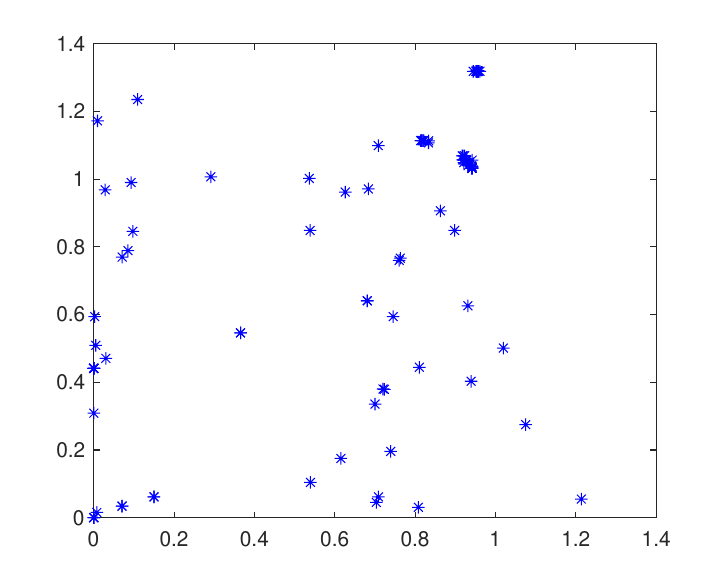}
\]
\caption{Performance of the plain ALS for the $2\times 2$ example. Left: convergence plots for $100$ initializations. Right: $SSS(\matr{A},\widehat{\matr{A}})$ versus $SSS(\matr{B},\widehat{\matr{B}})$  after 5000 ALS iterations.}
\label{fig:paratuck2_als}
\end{figure}

Algorithm~\ref{alg:paratuck_2x2_complete}  returns a Paratuck-2 approximation with  approximation error as $MSE(\tens{T},\widehat{\tens{T}})  = 1.7 \cdot 10^{-26}$ and factor distances 
$SSS(\matr{A},\widehat{\matr{A}}) =  6.6 \cdot 10^{-26}$ and $SSS(\matr{B},\widehat{\matr{B}}) =  2.2 \cdot 10^{-26}$

We also ran ALS(alternating least squares) \cite[pp. 68--71]{Bro98:thesis} with random initialization (factors drawn from the uniform distribution). With 5000 iterations, ALS gave an approximation error $MSE(\tens{T},\widehat{\tens{T}})$ below $10^{-10}$ only in $3\%$ of the cases (number of random initializations was 100).

\subsection{Random example}
\begin{figure}[hbt!]
\[
\includegraphics[height=4.8cm]{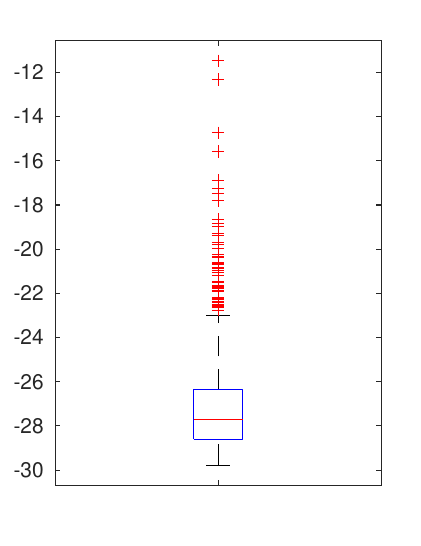} \quad \includegraphics[height=4.8cm]{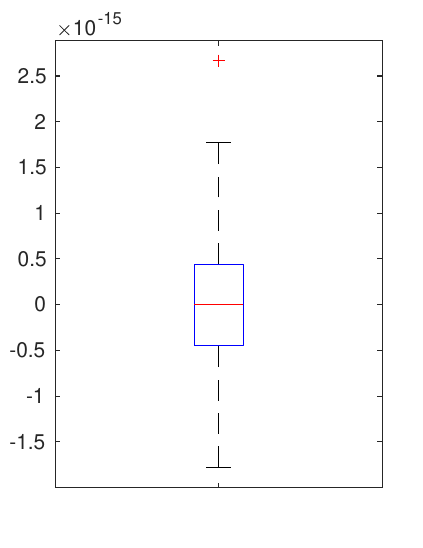} \quad \includegraphics[height=4.8cm]{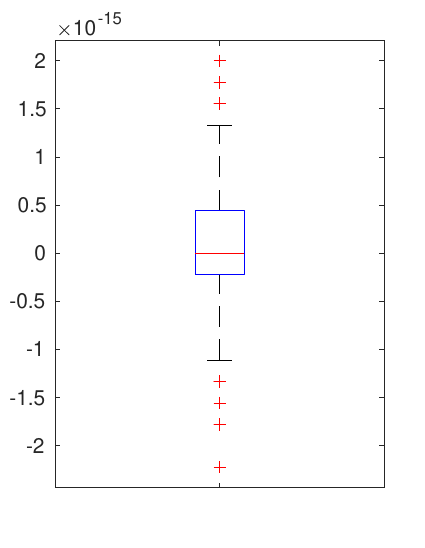}
\]
\caption{Performance of \Cref{alg:paratuck_2x2_complete} for random examples with $(R,S) = (4,4)$, boxplots for 1000 realizations. Left:  relative error in logarithmic scale $\log_{10}(MSE(\tens{T},\widehat{\tens{T}})/ \|\tens{T}\|^2)$, middle:  $SSS(\matr{A},\widehat{\matr{A}})$. Right: $SSS(\matr{B},\widehat{\matr{B}})$.
}
\label{fig:4x4_exact}
\end{figure}

We next consider random example with $I = J = 10$, $R=S=4$ and $K=100$.
We generate $\matr{A}$, $\matr{B}$, $\matr{H}$,  $\matr{G}$ having i.i.d. entries drawn randomly from $N(0,1)$,  and $\matr{F}$ from uniform distribution in $[0.5,1.5]$ (to guarantee \Cref{as:nonzero_F,as:krank_F}. ) without any perturbation explicitly added to the tensor.
For all $1000$ realizations of the factors, \Cref{alg:paratuck_2x2_complete} was able to accurately compute the factors of the PT2D, as shown in Fig.~\ref{fig:4x4_exact}.

\begin{figure}[hbt!]
$\includegraphics[height=5.8cm]{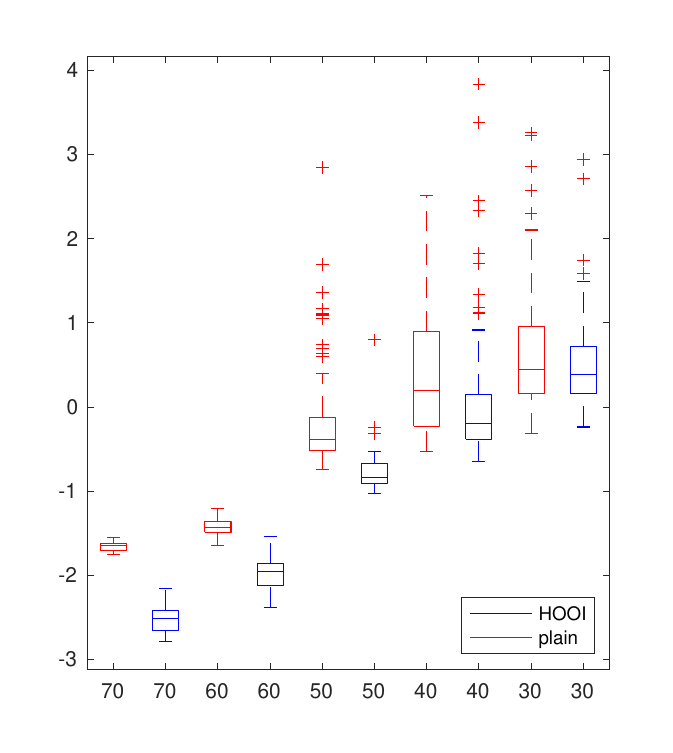} \quad 
\includegraphics[height=5.8cm]{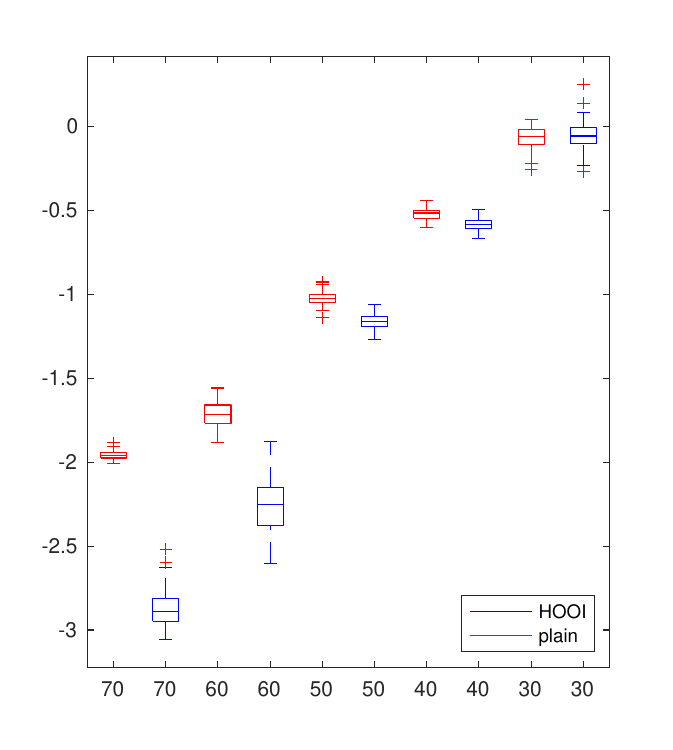} \includegraphics[height=5.8cm]{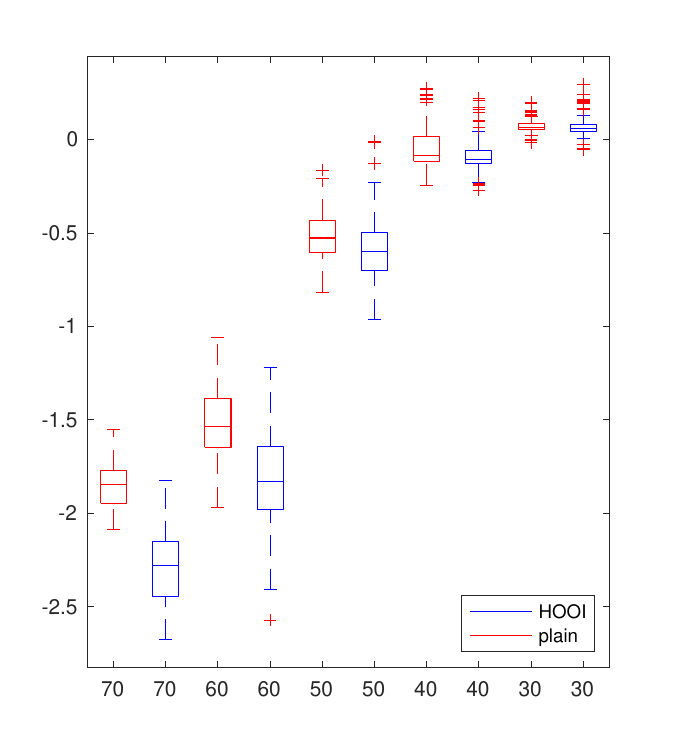}$
\caption{Performance of \Cref{alg:paratuck_2x2_complete} with respect to additive noise, for 100 realizations (plots in logarithmic scale for easier visualization). Left:  relative error $\log_{10}(MSE(\tens{T},\widehat{\tens{T}})/ \|\tens{T}\|^2)$. Middle:  $\log_{10}(SSS(\matr{A},\widehat{\matr{A}}))$. Right: $\log_{10}(SSS(\matr{B},\widehat{\matr{B}}))$.
For each of the dB levels, the pair of columns is given: without and with HOOI improvement (10 iterations).
}
\label{fig:4x4_noisy}
\end{figure}

\paragraph{Noisy scenario}
We also pick a single realization of factors from the previous experiments, and test the quality of the reconstruction of the factors for the case of additive noise $\tens{T} = \tens{T}_0 + \tens{E}$.
The level of noise is measured in dB, where $X$ dB  means that $10\log_{10} ({\|\tens{T}_0\|^2_F}/{\|\tens{E}\|^2_F}) = X$.
For example, $40$dB corresponds to the case $\|\tens{E}\|_F = 10^{-2} \|\tens{T}_0\|_F$.
In Fig.~\ref{fig:4x4_noisy}, we show the simulation results for random noise from $N(0,1)$, normalized to have achieve a prescribed noise level (essentially uniform distribution on a high-dimensional hypersphere $\|\tens{E}\|_F  = C$).

The results in Fig.~\ref{fig:4x4_noisy} show that \Cref{alg:paratuck_2x2_complete} still works for small perturbations, but the performance degrades for higher levels of noise.
We have tried an improvement based on takin into account the structure of the subspace space spanned by matrices $\matr{P}_k = \pi(\vect{p}_k)$ from \Cref{thm:pt2d_nonsymmetric}.
Note that from \Cref{prop:lker_unstructured} it follows that  the (noiseless)  $R^2 \times S^2 \times M$ tensor $\tens{P}$ build as $\tens{P}_{:,:,k} = \matr{P}_k$, must have the multilinear ranks equal to $(\binom{R}{2},\binom{S}{2}, \cdot)$.
Therefore, before performing steps $5$ of \Cref{alg:paratuck_2x2_complete}, we can compute the Tucker-2 approximation of  $\tens{P}$ by the higher-order orthogonal iteration (HOOI).
In our experiments (Fig.~\ref{fig:4x4_noisy}) HOOI gives a visible improvement, and we only need a few iterations of HOOI.

\begin{figure}[hbt!]
$\includegraphics[height=4.8cm]{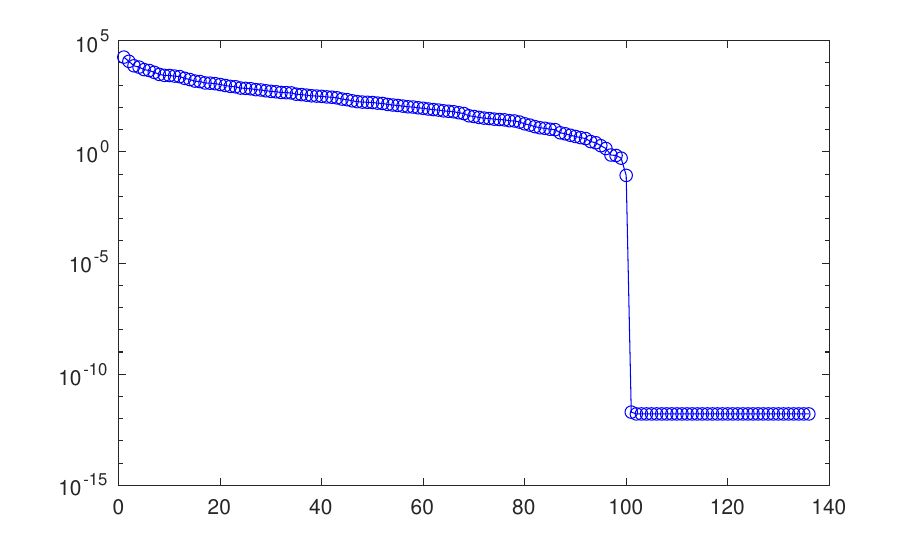} \includegraphics[height=4.8cm]{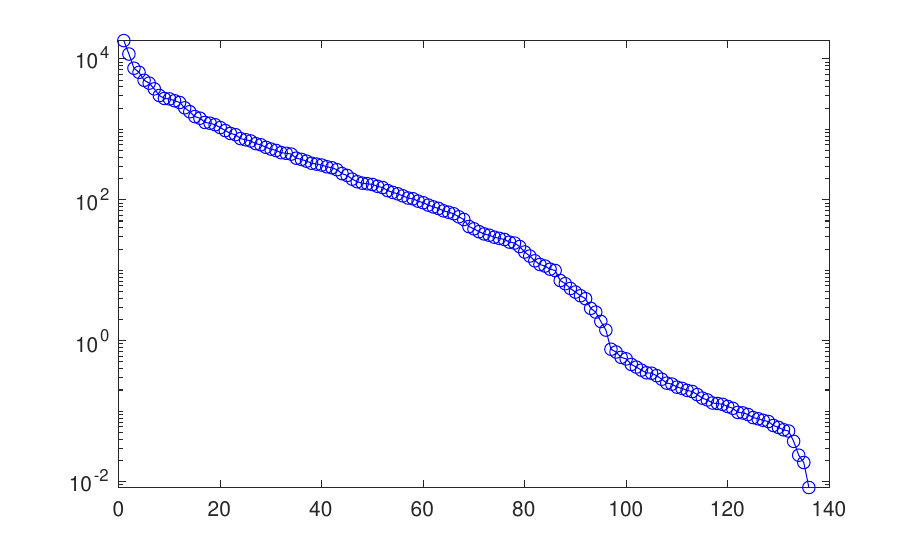}$
\caption{Singular values of $\matr{\Phi}(\tens{T})$. Left: noiseless scenario. Right: noise level $60$dB.}
\label{fig:sig_Phi}
\end{figure}

\paragraph{Singular values}
However, there is a fundamental limitation of \Cref{alg:paratuck_2x2_complete}, which is illustrated by the plots of singular values in Fig.~\ref{fig:sig_Phi}.
We see that $\matr{\Phi}(\tens{T})$ has $\binom{RS+1}{2}$ singular values in total, out of which we are looking for $\binom{R}{2}\binom{S}{2}$ smallest singular values.
When $\tens{T}$ is perturbed by the additive noise, this also results in a nonlinear perturbation of $\matr{\Phi}{(\tens{T})}$.

\subsection{Symmetric case}
Finally, to illustrate the performance of the algorithm for the DEDICOM (symmetric PT2D), we consider the example of PARAFAC-2 decomposition as in  \Cref{ex:PT2_PF2}.
We take $I=30$, $J_k = 40$, $K = 20$ and $R = 3$, where the factors
$\matr{A}$ is generated randomly from the uniform distribution in $[1,2]$ and $\matr{B}_k$ (of the same size), $\matr{D}_k$ are with elements drawn from the uniform distribution in $[0,1]$.
This setup was taken from the manual  of the TensorLy package \cite{tensorlydoc2024},
and we generate a single example, stacking all the 
Note that in this case the columns of the matrix $\matr{A}\in \RR^{30\times 3}$ are highly correlated, with the correlation matrix of the columns given by
\[
\approx \begin{bmatrix}
 1  &  0.97 &    0.95 \\    0.97   & 1   &  0.96\\   0.95  &  0.96&    1
\end{bmatrix}.
\]

\begin{figure}[hbt!]
\[
\includegraphics[height=4.8cm]{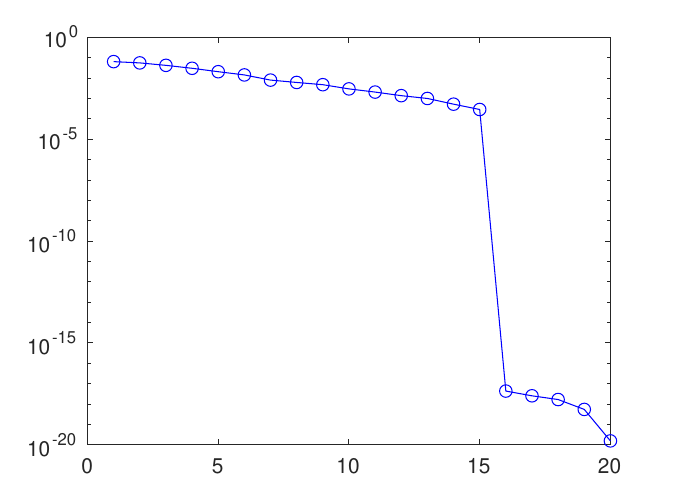}
\]
\caption{Singular values of $\matr{\Phi}(\widetilde{\tens{M}})$.}
\label{fig:parafac2_svalues}
\end{figure}

We proceed by adding additive noise  $\tens{T} = \tens{T}_0 + \tens{E}$ for different noise level and use the method in \Cref{ex:PT2_PF2} to reduce to the DEDICOM decomposition problem, with
dimension reduction using the SVD of the first unfolding of $\tens{T}$:
\begin{itemize}
\item the SVD  of the first unfolding $\matr{T}^{(1)} = \widetilde{\matr{U}} \widetilde{\Sigma} \widetilde{\matr{V}}^{\T}$ is computed;
\item the matrices $\matr{X}_k$ are projected on the leading $R$ left singular vectors: $\widetilde{X}_{k} = (\Sigma_{1:R,1:R})^{-1} (\matr{U}_{1:R,:}^\T)\matr{X}_k$.
\end{itemize} 
The algorithm \Cref{alg:paratuck_2x2_symmetric} is then used to compute the DEDICOM approximation of a tensor $\widetilde{\tens{M}} \in \RR^{3 \times 3 \times 20}$ with slices $\widetilde{\tens{M}}_{:,:,k} = \widetilde{\matr{X}}_{k}(\widetilde{\matr{X}}_{k})^{\T}$.
For reference, we plot in \Cref{fig:parafac2_svalues} the singular values of the matrix $\matr{\Phi}(\widetilde{\tens{M}})$.

\begin{figure}[hbt!]
\[
\includegraphics[height=5.5cm]{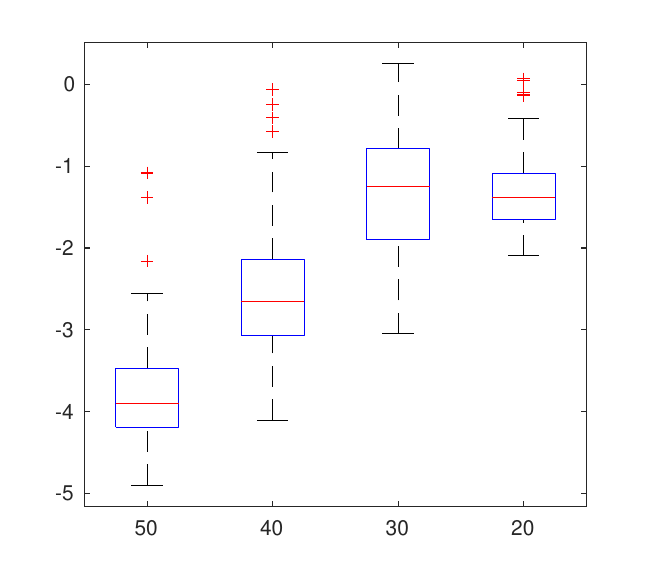} \quad\includegraphics[height=5.5cm]{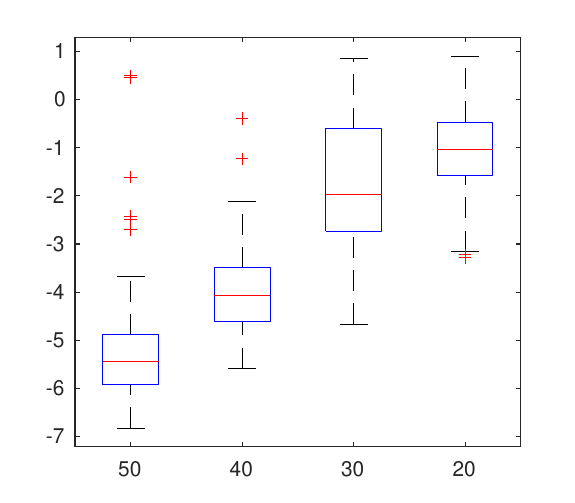}
\]
\caption{Reconstructions of the factor and of the tensor, for different noise levels (in dB), for 100 noise realizations.
  Left:  box plots of $\log_{10}(SSS(\matr{A},\widehat{\matr{A}}))$. Right: relative errors  $\log_{10}(MSE(\tens{M},\widehat{\tens{M}})/ \|\tens{M}\|^2)$.}
\label{fig:parafac2}
\end{figure}

In \Cref{fig:parafac2}, we report the results of approximation computed by DEDICOM. We see that the algorithm gives low error on the difference between the true and estimated factor $\matr{A}$.
We also plot the relative errors for estimation of the true (uncompressed) $\tens{M} \in \RR^{I\times I\times K}$, and see that the algorithm works for higher levels of noise than in the nonsymmetric case.

In terms of computational speed, the whole computation of a compression of matrices $\matr{X}_k$ together with DEDICOM  took just a few milliseconds in MATLAB R2024a on a MacBook laptop with M3 Pro chipset.
As a comparison, running the PARAFAC-2 decomposition for the same example with an iterative algorithm from \cite{tensorlydoc2024} took $\approx 5-6$ seconds for 500 iterations.

\section*{Conclusion}
In this paper we proposed a lifting approach to ParaTuck-2 decomposition and DEDICOM. This approach clarifies and relaxes the existing uniqueness conditions and leads to simple algebraic algorithms. The algebraic algorithms, that use just nullspace and eigenvalue computations, are the first (except the case of $R=2$ in DEDICOM), up to the author's knowledge, that are guaranteed to compute the exact decomposition and are also  are applicable in the approximation scenario, as shown by the numerical experiments.

The main drawback of the approach is its computational complexity, as it requires finding an SVD of a matrix with $\binom{R+1}{2}\binom{S+1}{2}$  matrix, thus the complexity grows very quickly as $R,S$ increase.
Another issue is that we need $\binom{R}{2}\binom{S}{2}$ last singular vectors, a number that also grows with $R,S$ and can affect the robustness to noise and numerical stability of the algorithm for high $R$, $S$.
However, we believe, that there are plenty of possible improvements of the basic algebraic algorithms presented in this paper, and that lifting approach can serve as the base for reliable algorithms to compute ParaTuck-2 decompositions.

\section*{Acknowledgements}
%
I would  like to thank many colleagues for stimulating  discussions.

\bibliographystyle{siamplain}
\bibliography{paratuck}

\end{document}

%% file: preamble.tex
\usepackage[utf8]{inputenc}
\usepackage{hyperref}
\usepackage[scale=0.78]{geometry}
\usepackage{graphicx}
\usepackage{amsmath,amssymb}
\usepackage{mathtools,mathdots,mathrsfs}
\usepackage{cleveref}
\usepackage{bm}
\usepackage{subcaption}
\usepackage{cancel}

\usepackage{algorithmic} 
\Crefname{ALC@unique}{Line}{Lines} 

\newcommand{\newmythm}[2]{
  \theoremstyle{plain}
  \theoremheaderfont{\normalfont\sc}
  \theorembodyfont{\normalfont\itshape}
  \theoremseparator{.}
  \theoremsymbol{}
  \newtheorem{#1}{#2}
}
\newmythm{assumption}{Assumption}

 \newsiamremark{remark}{Remark}
  \newsiamremark{example}{Example}
\crefname{assumption}{Assumption}{Assumptions}

\usepackage[textsize=footnotesize]{todonotes}

\usepackage{tikz}
\usetikzlibrary{shapes.geometric}

%% file: macros.tex
\definecolor{darkgreen}{rgb}{0,0.6,0}
\definecolor{lightred}{rgb}{1,0.2,0.2}

\newcommand{\tens}[1]{\boldsymbol{\mathcal{#1}}}
\newcommand{\tenselem}[1]{\mathcal{#1}}

\newcommand{\matr}[1]{\boldsymbol{#1}}

\newcommand{\bmx}{\begin{bmatrix}}
\newcommand{\emx}{\end{bmatrix}}
\newcommand{\bsm}{\left[\begin{smallmatrix}}
\newcommand{\esm}{\end{smallmatrix}\right]}

\newcommand{\vect}[1]{\boldsymbol{#1}}

\newcommand{\cpd}[3]{[\![#1, #2,#3]\!]}
\newcommand{\cpdfour}[4]{[\![#1, #2,#3,#4]\!]}
\makeatletter
\newcommand{\triprod}[3]{\mathop{\operator@font TP}(#1, #2,#3)}
\makeatother

\newcommand{\set}[1]{\mathscr{#1}}
\newcommand{\tensp}{\mathop{\otimes}}      
\newcommand{\kr}{\odot}     
\newcommand{\hadam}{\circledast}    
\newcommand{\kron}{\mathop{\boxtimes}}   
\newcommand{\con}{\mathop{\bullet}}   
\newcommand{\T}{{\sf T}}        

\makeatletter
\newcommand{\tensorize}[4]{\mathop{\operator@font tens}_{#1,#2,#3}\{#4\}}
\newcommand{\matricize}[3]{\mathop{\operator@font matr}_{#1,#2}\{#3\}}
\newcommand{\vecl}[1]{\mathop{\operator@font vec}\{#1\}}
\newcommand{\rank}[1]{\mathop{\operator@font rank}\{#1\}}
\newcommand{\colrank}[1]{\mathop{\operator@font colrank}\{#1\}}
\newcommand{\krank}[1]{\mathop{\operator@font krank}\{#1\}}
\newcommand{\trace}[1]{\mathop{\operator@font trace}\{#1\}}
\newcommand{\Diag}[1]{\mathop{\operator@font Diag}\{#1\}}    
\newcommand{\diag}[1]{\mathop{\operator@font diag}\{#1\}}    
\newcommand{\Span}[1]{\mathop{\operator@font Span}\{#1\}}    
\newcommand{\argmin}{\mathop{\operator@font argmin}}

\newcommand{\cond}[1]{\mathop{\operator@font cond}\{#1\}}

\newcommand{\lker}[1]{\mathop{\operator@font lker}\{#1\}}
\newcommand{\ptd}[5]{\mathop{\operator@font PT2D}(#1, #2,#3,#4,#5)}
\newcommand{\ptdsym}[3]{\mathop{\operator@font DEDICOM}(#1, #2,#3)}
\newcommand{\symlker}[1]{\sigma_2(\lker{#1})}
\newcommand{\symlkerd}[1]{\sigma_{(2,2)}(\lker{#1})}
\newcommand{\symlkert}[1]{\sigma_{3}(\lker{#1})}

\newcommand{\range}[1]{\mathop{\operator@font range}\{#1\}}

\makeatother

\newcommand{\RR}{\mathbb{R}}
\newcommand{\CC}{\mathbb{C}}

\newcommand{\FF}{\mathbb{F}}

\newcommand{\unfold}[2]{\matr{#1}^{(#2)}}
\newcommand{\contr}[1]{\con_{#1}}

\newcommand{\edim}[1]{%
\def\FirstArg{#1}%
\def\One{1}%
\def\Two{2}%
  \ifx\FirstArg\One%
    \ensuremath{I}%
    \else%
  \ifx\FirstArg\Two%
    \ensuremath{J}%
    \else%
    \ensuremath{K}
  \fi%
 \fi%
}

\newcommand{\rkRow}{R}
\newcommand{\rkCol}{S}

\newcommand{\DiagA}[1]{\matr{D}^{(#1)}_{\matr{\CoreVert}}}
\newcommand{\DiagB}[1]{\matr{D}^{(#1)}_{\matr{\CoreHoriz}}}

\newcommand{\Core}{C}
\newcommand{\CoreVert}{G}
\newcommand{\CoreHoriz}{H}
\newcommand{\CoreFront}{F}

\newcommand{\Mreshaped}[1]{\matr{P}_{#1}}

\newcommand{\symglobal}{\sigma}
\newcommand{\symperm}{\pi}

%% file: paratuck2_algo.bbl
\begin{thebibliography}{10}

\bibitem{BadeHK07:temporal}
{\sc B.~W. Bader, R.~A. Harshman, and T.~G. Kolda}, {\em Temporal analysis of
  semantic graphs using {ASALSAN}}, in Seventh IEEE international conference on
  data mining (ICDM 2007), IEEE, 2007, pp.~33--42.

\bibitem{Bro98:thesis}
{\sc R.~Bro}, {\em Multi-way Analysis in the Food Industry}, PhD thesis, Vrije
  Universiteit Brussel (VUB), 1998.

\bibitem{AlmeFX13:MIMO}
{\sc A.~L.~F. de~Almeida, G.~Favier, and L.~R. Ximenes}, {\em
  Space-time-frequency {(STF) MIMO} communication systems with blind receiver
  based on a generalized {PARATUCK2} model}, IEEE Transactions on Signal
  Processing, 61 (2013), pp.~1895--1909.

\bibitem{DeJoUDI23:CAMSAP}
{\sc J.~De~Jonghe, K.~Usevich, P.~Dreesen, and M.~Ishteva}, {\em Compressing
  neural networks with two-layer decoupling}, in 2023 IEEE 9th International
  Workshop on Computational Advances in Multi-Sensor Adaptive Processing
  (CAMSAP), 2023, pp.~226--230,
  \url{https://doi.org/10.1109/CAMSAP58249.2023.10403509}.

\bibitem{OlivFFB19:mimo}
{\sc P.~M.~R. de~Oliveira, C.~A.~R. Fernandes, G.~Favier, and R.~Boyer}, {\em
  {PARATUCK Semi-Blind Receivers for Relaying Multi-Hop MIMO Systems}},
  {Digital Signal Processing}, 92 (2019), pp.~127--138,
  \url{https://doi.org/10.1016/j.dsp.2019.05.011}.

\bibitem{DiacS98:aos}
{\sc P.~Diaconis and B.~Sturmfels}, {\em {Algebraic algorithms for sampling
  from conditional distributions}}, The Annals of Statistics, 26 (1998),
  pp.~363 -- 397, \url{https://doi.org/10.1214/aos/1030563990}.

\bibitem{DomaL14:gevd}
{\sc I.~Domanov and L.~D. Lathauwer}, {\em Canonical polyadic decomposition of
  third-order tensors: Reduction to generalized eigenvalue decomposition}, SIAM
  Journal on Matrix Analysis and Applications, 35 (2014), pp.~636--660,
  \url{https://doi.org/10.1137/130916084}.

\bibitem{DreeIS15:decoupling}
{\sc P.~Dreesen, M.~Ishteva, and J.~Schoukens}, {\em Decoupling multivariate
  polynomials using first-order information and tensor decompositions}, SIAM
  Journal on Matrix Analysis and Applications, 36 (2015), pp.~864--879,
  \url{https://doi.org/10.1137/140991546}.

\bibitem{FaviA14:constrained}
{\sc G.~Favier and A.~L. de~Almeida}, {\em Overview of constrained {PARAFAC}
  models}, EURASIP Journal on Advances in Signal Processing, 2014 (2014),
  p.~142.

\bibitem{GarcSS05:AG}
{\sc L.~D. Garcia, M.~Stillman, and B.~Sturmfels}, {\em Algebraic geometry of
  {Bayesian} networks}, Journal of Symbolic Computation, 39 (2005),
  pp.~331--355, \url{https://doi.org/10.1016/j.jsc.2004.11.007}.
\newblock Special issue on the occasion of MEGA 2003.

\bibitem{Hars78:dedicom}
{\sc R.~A. Harshman}, {\em {Models for analysis of asymmetrical relationships
  among n objects or stimuli}}, in {First Joint Meeting of the Psychometric
  Society and the Society for Mathematical Psychology}, McMaster University,
  Hamilton, Ontario, Aug. 1978.

\bibitem{HarsL96:uniqueness}
{\sc R.~A. Harshman and M.~E. Lundy}, {\em Uniqueness proof for a family of
  models sharing features of {Tucker's} three-mode factor analysis and
  parafac/candecomp}, Psychometrika, 61 (1996), pp.~133--154.

\bibitem{Kier89:ALS}
{\sc H.~A. Kiers}, {\em An alternating least squares algorithm for fitting the
  two-and three-way {DEDICOM} model and the idioscal model}, Psychometrika, 54
  (1989), pp.~515--521.

\bibitem{Kier93:ALS}
{\sc H.~A. Kiers}, {\em An alternating least squares algorithm for {PARAFAC2}
  and three-way {DEDICOM}}, Computational Statistics \& Data Analysis, 16
  (1993), pp.~103--118, \url{https://doi.org/10.1016/0167-9473(93)90247-Q}.

\bibitem{KierTB99:direct}
{\sc H.~A. Kiers, J.~M. Ten~Berge, and R.~Bro}, {\em {PARAFAC2—Part I.} a
  direct fitting algorithm for the {PARAFAC2} model}, Journal of Chemometrics:
  A Journal of the Chemometrics Society, 13 (1999), pp.~275--294.

\bibitem{KoldB09:tensor}
{\sc T.~G. Kolda and B.~W. Bader}, {\em Tensor decompositions and
  applications}, SIAM review, 51 (2009), pp.~455--500.

\bibitem{Nask20:thesis}
{\sc K.~Naskovska}, {\em Advanced tensor based signal processing techniques for
  wireless communication systems and biomedical signal processing}, PhD thesis,
  Ilmenau, Jan 2020.
\newblock Dissertation, Technische Universit{\"a}t Ilmenau, 2019.

\bibitem{RoalSCABCA22:AOADMM}
{\sc M.~Roald, C.~Schenker, V.~Calhoun, T.~Adali, R.~Bro, J.~Cohen, and
  E.~Ataman}, {\em An {AO-ADMM} approach to constraining {PARAFAC2} on all
  modes}, SIAM Journal on Mathematics of Data Science, 4 (2022),
  pp.~1191--1222, \url{https://doi.org/10.1137/21M1450033}.

\bibitem{BergK96:uniqueness}
{\sc J.~M. ten Berge and H.~A. Kiers}, {\em Some uniqueness results for
  {PARAFAC2}}, Psychometrika, 61 (1996), pp.~123--132.

\bibitem{tensorlydoc2024}
{\sc {TensorLy developers}}, {\em Demonstration of parafac2 (tensorly package
  manual)}, 2024.
\newblock
  \url{tensorly.org/stable/auto_examples/decomposition/plot_parafac2.html}
  [Accessed: 13 November 2024].

\bibitem{UsevZIDA23:EUSIPCO}
{\sc K.~Usevich, Y.~Zniyed, M.~Ishteva, P.~Dreesen, and A.~de~Almeida}, {\em
  {Tensor-based two-layer decoupling of multivariate polynomial maps}}, in
  {31st European Signal Processing Conference, EUSIPCO 2023}, Helsinki,
  Finland, Sept. 2023, {European Association for Signal Processing},
  \url{https://doi.org/10.23919/EUSIPCO58844.2023.10289900}.

\bibitem{ZniyA25:stochastic}
{\sc Y.~Zniyed and A.~L. de~Almeida}, {\em A stochastic algorithm for the
  {ParaTuck} decomposition}, Digital Signal Processing, 156 (2025), p.~104767.

\end{thebibliography}
